\documentclass[11pt,a4paper]{article} 
\usepackage[utf8]{inputenc}
\usepackage[english]{babel}
\usepackage[T1]{fontenc}
\usepackage{amsmath}
\usepackage{amsthm}
\usepackage{comment}
\usepackage{amsfonts}
\usepackage{amssymb}
\usepackage{multicol}

\usepackage{url}
\usepackage{pdfpages}

\usepackage{graphicx}
\usepackage{tikz-qtree}
\usepackage[totoc]{idxlayout}

\usepackage[intoc,refpage]{nomencl}

\usepackage{graphicx}
\graphicspath{ {images/} }
\usepackage[all,cmtip]{xy}

\usepackage[toc,page]{appendix}
\usepackage[left=2cm,right=2cm,top=3cm,bottom=3cm]{geometry}
\usepackage{pgf,tikz}
\usepackage[font=small,skip=0pt]{caption}

\usetikzlibrary{arrows}

\usepackage{etoolbox}
\renewcommand\nomgroup[1]{%
  \item[\bfseries
  \ifstrequal{#1}{P}{Notation for properties of Henselian valued fields}{%
  \ifstrequal{#1}{L}{Languages}{%
  \ifstrequal{#1}{SE}{Stable embeddedness}{%
  \ifstrequal{#1}{C}{Classification theory}{%
  \ifstrequal{#1}{A}{Adler's convention and more}{%
  \ifstrequal{#1}{O}{Other Symbols}{}}}}}}%
]}

\makeatletter
\newcommand\ackname{Acknowledgements}
\if@titlepage
   \newenvironment{acknowledgements}{%
       \titlepage
       \null\vfil
       \@beginparpenalty\@lowpenalty
       \begin{center}%
         \bfseries \ackname
         \@endparpenalty\@M
       \end{center}}%
      {\par\vfil\null\endtitlepage}
\else
   \newenvironment{acknowledgements}{%
       \if@twocolumn
         \section*{\abstractname}%
       \else
         \small
         \begin{center}%
           {\bfseries \ackname\vspace{-.5em}\vspace{\z@}}%
         \end{center}%
         \quotation
       \fi}
       {\if@twocolumn\else\endquotation\fi}
\fi
\makeatother

\def\restriction#1#2{\mathchoice
              {\setbox1\hbox{${\displaystyle #1}_{\scriptstyle #2}$}
              \restrictionaux{#1}{#2}}
              {\setbox1\hbox{${\textstyle #1}_{\scriptstyle #2}$}
              \restrictionaux{#1}{#2}}
              {\setbox1\hbox{${\scriptstyle #1}_{\scriptscriptstyle #2}$}
              \restrictionaux{#1}{#2}}
              {\setbox1\hbox{${\scriptscriptstyle #1}_{\scriptscriptstyle #2}$}
              \restrictionaux{#1}{#2}}}
\def\restrictionaux#1#2{{#1\,\smash{\vrule height .8\ht1 depth .85\dp1}}_{\,#2}} 

\newcommand*{\upwedge}{\mathbin{\raisebox{0ex}{\textasciicircum}}}

\DeclareMathOperator{\acl}{acl}

\DeclareMathOperator{\bdn}{bdn}
\DeclareMathOperator{\TP}{TP}
\DeclareMathOperator{\NTP}{NTP}

\DeclareMathOperator{\tp}{tp}
\DeclareMathOperator{\qftp}{qftp}
\DeclareMathOperator{\val}{val}
\DeclareMathOperator{\res}{res}
\DeclareMathOperator{\Res}{Res}
\DeclareMathOperator{\ac}{ac}
\DeclareMathOperator{\WD}{WD}
\DeclareMathOperator{\rv}{rv}
\DeclareMathOperator{\RV}{RV}
\DeclareMathOperator{\inp}{inp}
\DeclareMathOperator{\Th}{Th}

\DeclareMathOperator{\cf}{cf}
\DeclareMathOperator{\act}{act}
\DeclareMathOperator{\Card}{Card}

\DeclareMathOperator{\supp}{supp}

\author{Pierre Touchard  
\footnote{The author was funded
by DFG through SFB 878 and under Germany 's Excellence Strategy – EXC 2044 – 390685587, Mathematics Münster: Dynamics – Geometry - Structure. This research has been also supported by the DAAD through the `Kurzstipendien für Doktoranden 2020/21' programme and by the MFO through the `Oberwolfach Leibniz Fellows 2020' programme.}}
\title{Burden in Henselian Valued Fields}

\theoremstyle{plain}
\newtheorem{theorem}{Theorem}[section]
\newtheorem{lemma}[theorem]{Lemma}
\newtheorem{proposition}[theorem]{Proposition}
\newtheorem{claim}{Claim}
\newtheorem{remark}[theorem]{Remark}
\newtheorem{remarks}[theorem]{Remarks}
\newtheorem{observation}[theorem]{Observation}
\newtheorem{fact}[theorem]{Fact}

\theoremstyle{definition}
\newtheorem{definition}[theorem]{Definition}
\newtheorem{corollary}[theorem]{Corollary}

\newtheorem{heuristic}[theorem]{Heuristic}

\theoremstyle{remark}
\newtheorem*{remark*}{Remark}

\newtheorem*{note}{Note}

\newtheorem*{notation}{Notation}

\definecolor{black}{rgb}{0,0,0}
\newtheorem*{example}{Example}
\newtheorem*{examples}{Examples}

\makeindex

\begin{document}

\date{\today}
\maketitle
\vspace{-0.5cm}
\begin{abstract}
    In the spirit of the Ax-Kochen-Ershov principle, we show that in certain cases the burden of a Henselian valued field can be computed in terms of the burden of its residue field and that of its value group. To do so, we first see that the burden of such a field is equal to the burden of its leading term structure. These results are generalisations of Chernikov and Simon's work in \cite{CS19}.
\end{abstract}

\tableofcontents
\newpage

\section*{Introduction}
Since the work of Shelah \cite{She90}, model theorists have defined and studied combinatorial
configurations that first order theories may encode, in order to measure their relative tameness. Theories which are able
to encode complex configurations are considered less tame. Arising from
this classification  a complex hierarchy of first order theories. The most important class is undoubtedly
that of stable theories; the classes of NIP (Non-Independance Property), 
simple and $\text{NTP}_2$ theories (Non-Tree Property of the second kind) are
generalisations of stability which have been extensively studied in the last
20 years. The abstract study of these classes leads to a better understanding
of algebraic structures, as some algebraic phenomena may or not occur,
depending of the complexity of the theory (e.g. \cite{KSW11}). Locating a given concrete first-order theory in this hierarchy is often interesting and challenging, and many methods are known: one can show stability (resp. NIP) by counting types over small sets (resp. coheirs over models). Simple theories are exactly theories with some
abstract independence relation (\cite{KP97}). The situation seems to be different when
we look at quantitative versions of these notions, such the burden. This is a notion of dimension for $\text{NTP}_2$ theories, and in order to compute it, a concrete good understanding of formulas is required.

Since Ax, Kochen and Ershov's  work, Henselian valued fields appear to be perfect playgrounds for this kind of consideration, and many model theoretic questions on valued fields have been reduced to the residue field and value group. The theorem of Delon \cite{Del81} is one of the first instances: a Henselian valued field of equicharacteristic 0 is NIP if and only if both of its residue field and its value group are NIP\footnote{It was later show by Gurevich and Schmitt that any pure ordered abelian group is NIP.}. A more quantitative transfer in NIP Henselian valued fields was then showed by Shelah in \cite{She14}: a Henselian valued field of equicharacteristic 0 is strongly dependant if and only if both its residue field and its value group are strongly dependant. Both results were generalised to $\text{NTP}_2$ theories by Chernikov in \cite{Che14}: a Henselian valued field of equicharacteristic 0 is $\text{NTP}_2$ (resp. strong, of finite burden) if and only if both of its residue field and its value group are $\text{NTP}_2$ (resp. strong, of finite burden). This approach is based on the faith that the study of the residue field and the value group might be enough to classify some rather nice valued fields. In the case of equicharacteristic $0$ Henselian valued fields, this faith is justified by the theorem of Pas: a Henselian valued field of equicharacteristic zero equipped with an angular component (also called $\ac$-map), eliminates quantifiers relative to the value group and the residue field. It is indeed an important tool for producing transfer principles. However, the theorem of Pas has some limits: considering an $\ac$-map has the disadvantage of adding new definable sets to the structure of valued fields. From the point of view of complexity, it has an impact: any ultraproduct of p-adic fields over a non-principal utrafilter on prime numbers is inp-minimal, i.e. of burden 1 (see \cite{CS19}), but it is of burden 2 when endowed with an ac-map. Another approach initiated by Basarab and Kuhlmann, is to consider another interpretable sort capturing both information from the value group and the residue field. This lead to the definition of the \textit{leading term structure} (see Paragraph \ref{SusubsectionRVSort}), also called the $\RV$-sort. Unlike the ac-map, it is always interpretable in the standard languages of valued fields, and adding it to the language does not add definable sets. The study of valued fields together with its $\RV$-sort offers an additional point of view: let us cite Hrushovski and Kazdan's work in motivic integration, where RV-sort structures are used, as opposed to Denef, Loeser and Cluckers' work, where ac-maps are used. More relevant to this paper, Chernikov and Simon prove in \cite{CS19} that, under some hypothesis, an equicharacteristic zero Henselian valued field is inp-minimal if only if both its residue field and its value group are inp-minimal, going via an intermediate step: they first reduce the question to the RV-sort.

The main aim of this paper is to give a general transfer principle for burden in certain Henselian valued fields.  In particular, we provide a full answer to \cite[Problems 4.3 \& 4.4]{CS19}. Here is an overview of the paper:\\
The first section consists of preliminaries on pure model theory and on model theory of valued fields and groups. We define the burden of a theory, indiscernible sequences and few related lemmas. Then we define the sort $\RV$ and recall some relative quantifier elimination results.\\

Section 2 is dedicated to the proof of the following theorem:
    {
    \renewcommand{\thetheorem}{\ref{ThmBdnExSeq}}
    \begin{theorem}
        Consider an $\{A\}$-$\{C\}$-enrichement of a pure exact sequence $\mathcal{M}$ of abelian groups
    \[ \xymatrix{0 \ar[r] & A \ar[r]^{\iota}& B \ar[r]^{\nu} & C\ar[r] &0}, \]
in a language $\mathrm{L}$. Let $\mathcal{D}=\mathcal{D}(x)$ be the set of formulas in the pure language of groups which are conjunction of formulas of the form $\exists y\  nx=my$ for $n,m \in \mathbb{N}$.
For $D(x)\in \mathcal{D}$ and $A$ an abelian group, $D(A)$ is an subgroup of $A$, and we have 
        \[ \bdn\mathcal{M} =\max_{D\in \mathcal{D}}(\bdn(A/D(A))+\bdn(D(C))).\]
        
        In particular:
        \begin{itemize}
            \item If $A/nA$ is finite for all $n\geq 1$, then 
        \[\bdn{\mathcal{M}}= \max_{k\in \mathbb{N}}(\bdn(kA), \bdn(C_{[k]})),\]
        where $C_{[k]}:= \{c\in C \ \vert  \ kc=0\}$ is the subgroup of $k$-torsion.
            \item If $C$ has finite $k$-torsion of all $k\geq 1$, then 
        \[\bdn{\mathcal{M}}= \max_{n\in \mathbb{N}}(\bdn(A/nA) + \bdn(nC)).\] 
            \item If $C$ has finite $n$-torsion and $A/nA$ is finite for all $n\geq 1$, then 
            \[\bdn{\mathcal{M}}= \max(\bdn(A), \bdn(C)).\]
        
        \end{itemize}
    \end{theorem}
    \addtocounter{theorem}{-1}
    }   
     We also give an easy generalisation, notably for short exact sequences of ordered abelian groups.\\

In Section 3, we prove our main theorems. The first concerns the following partial theories of Henselian valued fields, that we gather here under the name of \emph{benign} theories of valued fields:
    \begin{enumerate}
        \item Henselian valued fields of characteristic $(0,0)$,
        \item algebraically closed valued fields, 
        \item algebraically maximal Kaplansky valued fields.
    \end{enumerate}

{
    \renewcommand{\thetheorem}{\ref{ThmBdnHenValFieCha00}}
    \begin{theorem}
        Let $\mathcal{K}=(K, \Gamma, k)$ be a benign Henselian valued field, with value group $\Gamma$ and residue field $k$. Then: 
        \[\bdn(\mathcal{K})= \max_{n\geq 0} \left(\bdn(k^\star/{k^\star}^n) + \bdn (n\Gamma) \right).\]
    \end{theorem}
    \addtocounter{theorem}{-1}
}   
The second concerns unramified mixed characteristic Henselian valued fields with perfect residue field.  
    {
    \renewcommand{\thetheorem}{\ref{theoremmixedchar}}
    \begin{theorem}
        Let $\mathcal{K}=(K,k,\Gamma)$ be an unramified mixed characteristic Henselian valued field  with value group $\Gamma$ and residue field $k$. We denote by $\mathcal{K}_{\ac_{<\omega}}=(K,k,\Gamma,\ac_n,n<\omega)$ the structure $\mathcal{K}$ endowed with compatible $\ac$-maps. Assume that the residue field $k$ is perfect. One has 
        \[\bdn(\mathcal{K})= \bdn(\mathcal{K}_{\ac_{<\omega}})= \max(\aleph_0 \cdot \bdn(k), \bdn(\Gamma)).\]
    \end{theorem}
    \addtocounter{theorem}{-1}
    }
 For both proofs, and as in \cite{CS19}, we proceed in two steps: we show first that the burden of these Henselian valued fields is equal to the burden of their $\RV$-sorts (see Theorem \ref{ThmHensValuedFieldReductionRV}) and use Section 2 to conclude. Finally, we produce in Appendix \ref{LexicographicProduct} a similar transfer principle for burden in Lexicographic products and in Appendix \ref{AxiomatisationRV}, we give an axiomatisation of $\RV$-structures and show how to construct a canonical Henselian valued field from any given $\RV$-structure.

Let us conclude this introduction with some generalities about transfer principles. We summarise the strategy presented above by formalising the reduction in valued fields and pure short exact sequences of abelian groups. We introduce \emph{reduction diagrams}. It is nothing else than a concise way to picture relative quantifier elimination and by extension, the strategy for proving reduction principles. 
    
        \begin{heuristic}\label{DefRedDiag}
    A \emph{reduction diagram} of a structure $\mathcal{M}$ is a rooted tree such that: 
    \begin{itemize}
        \item all nodes are pure sorts of $\mathcal{M}$ (in some $\emptyset$-interpretable language) endowed with their full structure;
        \item the root is $\mathcal{M}$;
        \item any node admits relative resplendent quantifier elimination (in some $\emptyset$-interpretable language) to the set of its children;
        \item any two sorts in two different branches are orthogonal. 
    \end{itemize}
    The idea is that one might be able to reduce certain questions on the structure $\mathcal{M}$ to the set of its leaves. Every node describes then an intermediate step.    Reduction to a node would also have the advantage of being generalised to any enrichment of structure below the node.
    \end{heuristic}  
    In this text, we compute the burden (Definition \ref{DefBurden1}) of the following examples in terms of the burden of the leaves.  In \cite{Tou20b} we also characterise stable embeddedness of elementary pairs of models in terms of stable embeddedness of elementary pairs of structures in the leaves. 
    \begin{example} 
    \begin{enumerate}
     \item If $\mathcal{M}_0,\mathcal{M}_1$ are arbitrary structures, both the direct product $\mathcal{M}_0 \times \mathcal{M}_1$ and the disjoint union $\mathcal{M}_0 \cup \mathcal{M}_1$ reduce to $\mathcal{M}_0$ and $\mathcal{M}_1$ (Fact \ref{FactQEProd}):
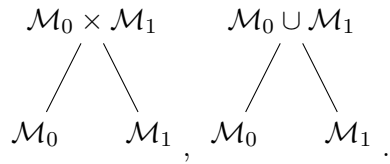
    \begin{center}
        \begin{tikzpicture}
        \node{$\mathcal{M}_0 \times \mathcal{M}_1$}
    child {  node {$\mathcal{M}_0$}}
    child { node {$\mathcal{M}_1$} };
    \end{tikzpicture},
    \begin{tikzpicture}
        \node{$\mathcal{M}_0 \cup \mathcal{M}_1$}
    child {  node {$\mathcal{M}_0$}}
    child { node {$\mathcal{M}_1$} };
    \end{tikzpicture}.
    \end{center}
    We can of course keep going:
    
    
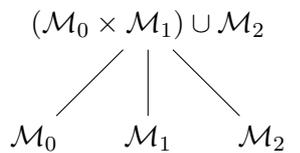
\begin{center}
    \begin{tikzpicture}
        \node{$\left(\mathcal{M}_0 \times \mathcal{M}_1 \right) \cup \mathcal{M}_2$}
    child {node {$\mathcal{M}_0$}}
    child {node {$\mathcal{M}_1$}}
    child {node {$\mathcal{M}_2$}};
    \end{tikzpicture}.
    \end{center}
        \item Let $\mathcal{K}_{ac} =\{K,\Gamma,k,\ac:K \rightarrow k \}$ be a Henselian valued field of equicharacteristic $0$ of valued group $\Gamma$, residue field $k$, and angular component $\ac$. It admits the following reduction diagrams (Theorem of Pas):
            
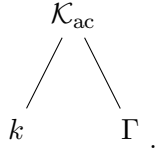
\begin{center}
    \begin{tikzpicture}
        \node{$\mathcal{K}_{\ac}$}
    child {  node {$k$}}
    child { node {$\Gamma$} };
    \end{tikzpicture}.
    \end{center}
    
        \item Let $\mathcal{M}=\{A,B,C, \iota, \nu\}$ be a short  exact sequence of abelian groups  \[ \xymatrix{0 \ar[r] & A \ar[r]^{\iota}& B \ar[r]^{\nu} & C\ar[r] &0}, \]
        seen as a three-sorted structure. Assume $A$ is a pure subgroup of $B$. It admits the following reduction diagram (Fact \ref{FactACGZ}):
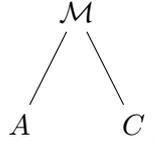
    \begin{center}
        \begin{tikzpicture}
            \node{$\mathcal{M}$}
        child {  node {$A$}}
        child { node {$C$} };
        \end{tikzpicture}.
    \end{center}
        To get relative quantifier elimination, one has to consider interpretable maps from $B$ to $A/nA$, $n\geq 0$. The sort $A/nA$ are understood to be part of the induced structure on $A$.
        \item Let $\mathcal{K} =\{K,\Gamma,k,\RV(K) \}$ be a Henselian valued field of equicharacteristic $0$, value group $\Gamma$, residue field $k$ and RV-sort $\RV(K)$ (definition in Paragraph \ref{SusubsectionRVSort}). It admits the following reduction diagram (Fact \ref{FactRelQERV} and Fact \ref{FactACGZ}):
        
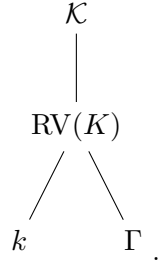
    \begin{center}
    \begin{tikzpicture}

        \node{$\mathcal{K}$}
    child {  node {$\RV(K)$}
        child { node {$k$} }
        child { node {$\Gamma$} }
    };
    \end{tikzpicture}.
    \end{center}
    
        \item If $\mathcal{K} =\{K,\Gamma,k\}$ is a Henselian valued field of equicharacteristic $0$, where moreover the residue field $k$ is endowed with a structure $(k,\Gamma',k')$ of Henselian valued field of equicaracteristic $0$, and $\Gamma$ is endowed with a predicate for a convex subgroup $\Delta$. Then by Corollary \ref{CoroPureOrth} (and resplendence), we have the following reduction diagram:
                    \begin{center}
    \begin{tikzpicture}
        \node{$\mathcal{K}$}
    child{ node {$\RV(K)$}
        child {  node {$k$}
            child{node{$\RV(k)$} 
                child{node {$k'$}}
                child{node {$\Gamma'$}}
            }
        }
        child[missing]
        child { node {$\Gamma$} 
            child{node{$\Delta$}}
            child{node{$\Gamma/\Delta$}}
        }
    };
    \end{tikzpicture}.
    \end{center}
    \end{enumerate}
    \end{example}
    \newpage

\section{Preliminaries} \label{sectionPreliminaries}

\subsection{On pure model theory}\label{SectionOnPureModelTheory}
 We will assume the reader to be familiar with basic model theory concepts, and in particular with standard notations. One can refer to \cite{TZ12}. Symbols $x,y,z,\dots$ will usually refer to tuples of variables, $a,b,c,\dots$ to parameters. Capital letters $K,L,M,N,\dots$ will refer to sets, and calligraphic letters $\mathcal{K},\mathcal{L},\mathcal{M},\mathcal{N},\dots$ will refer to structures with respective base sets $K,L,M,N,\dots$. If there is no ambiguity, we may respectively name a very saturated elementary extension with blackboard bold letters $\mathbb{K},\mathbb{L},\mathbb{M},\dots$ Languages will be denoted will a roman character $\mathrm{L},\mathrm{L}',\mathrm{L}_{Rings}, \mathrm{L}_{\Gamma,k}$ etc.

 In this section, we will consider any (possibly multi-sorted) first order language $\mathrm{L}$, and an arbitrary $\mathrm{L}$-structure $\mathcal{M}$. 
 
\subsubsection{Relative quantifier elimination and resplendence}\label{SubsectionRelativeQuantifierElimination}
    We will use freely Rideau-Kikuchi's terminology about enrichment that we briefly recall now. The reader can refer to \cite[Annexe A]{Rid17} for a more detailed exposition.
    This will allow us to give effortless generalisations to richer structures, or to simplify the notation, by reducing the language to the strict necessity for producing transfer principles.
    
    First, let us recall two notions of relative quantifier elimination. 
    
    \begin{definition}
        Let $\mathcal{M}$ be a multisorted structure in a language $\mathrm{L}$, and consider $\Pi \cup \Sigma$ a partition of the set of sorts. We denote by $\restriction{\mathrm{L}}{\Sigma}$ the language of all function symbols and relation symbols in $\mathrm{L}$ involving only sorts in $\Sigma$.  Then, we say that
        \begin{itemize}
            \item $\mathcal{M}$ eliminates $\Pi$-quantifiers if every formula $\phi(x)$ is equivalent to a formula without quantifier in a sort in $\Pi$.
            \item $\mathcal{M}$ eliminates quantifiers relatively to $\Sigma$ if the theory of $\mathcal{M}^{\Sigma-\text{Mor}}$ -- obtained by naming all $\restriction{\mathrm{L}}{\Sigma}$-definable sets (without parameters) with a new predicate-- eliminates quantifiers.  \index{Relative quantifier elimination}
        \end{itemize}
    \end{definition}
    
    As observed in \cite[Annexe A]{Rid17}, $M$ eliminates quantifiers relatively to $\Sigma$, then it eliminates $\Pi$-quantifiers. 
    
    \begin{definition}
        Let $\mathcal{M}$ be a multi-sorted structure in a language $\mathrm{L}$, and let $\Sigma$ be a set of sorts in $\mathrm{L}$.
        \begin{itemize}
            \item a language $\mathrm{L}_e$ containing $\mathrm{L}$ is said to be a $\Sigma$\emph{-enrichment} of $\mathrm{L}$ if all new function symbols and relation symbols only involve the sorts in $\Sigma$ and the new sorts $\Sigma_e$ in  $\mathrm{L}_e \setminus \mathrm{L}$. An expansion $\mathcal{M}_e$ of $\mathcal{M}$ to $\mathrm{L}_e$ is called a $\Sigma$\emph{-enrichment} of $\mathcal{M}$.
            \item $\Sigma$ is said to be \emph{closed} if any relation symbol involving a sort in $\Sigma$ or any function symbol with a domain involving a sort in $\Sigma$ only involves sorts in $\Sigma$.
        \end{itemize}
    \end{definition}
    
    \begin{fact}[\cite{Rid17}]
        Let $\mathcal{M}$ be a multisorted structure, and consider $\Pi \cup \Sigma$ a partition of the set of sorts.
        If $\Sigma$ is a closed set of sorts, then $\mathcal{M}$ eliminates $\Pi$-quantifiers if and only if $\mathcal{M}$ eliminates quantifiers relatively to $\Sigma$.
    \end{fact}
    
    In the context of this text, these two notions of quantifier elimination will be often equivalent. Another consequence of closedness is the automatic resplendence of relative quantifier elimination:

    \begin{definition}
      Let $\mathcal{M}$ be a multi-sorted structure in a language $\mathrm{L}$, and let $\Sigma$ be a set of sorts in $\mathrm{L}$.
      We say that $\mathcal{M}$ \emph{eliminates quantifiers resplendently relatively}\index{Relative quantifier elimination! resplendent ~} to $\Sigma$ if for any $\Sigma$-enrichment  $\mathcal{M}_e$ of $\mathcal{M}$, $\Th(\mathcal{M}_e)$ eliminates quantifiers relatively to $\Sigma \cup \Sigma_e$ (where $\Sigma_e$ is the set of new sorts in $\mathrm{M}_e$ ).
    \end{definition}
    \begin{fact}[{\cite[Proposition A.9]{Rid17}}]\label{FactClosedSortResplQE}
        Let $\mathcal{M}$ be a multi-sorted structure in a language $\mathrm{L}$, and assume that $\Th(\mathcal{M})$ eliminates quantifiers relative to a closed set of sorts $\Sigma$. Then $\Th(\mathcal{M})$ eliminates quantifiers resplendently relatively to $\Sigma$.  
    \end{fact}

    Notice however that closedness does not characterise resplendence of relative quantifier elimination, as we will see later with pure short exact sequence of abelian groups. Let us introduce the notion of stable embedded definable sets and of pure sorts. 
    \begin{definition}
        \begin{itemize}
            \item A definable subset $D$ of $\mathcal{M}$ is called \emph{stably embedded} if all definable subsets of $D^n$, $n\in \mathbb{N}$ can be defined with parameters in $D$. \index{Stable embeddedness}
            \item  Two definable subsets $D$ and $D'$ of $\mathcal{M}$ are called \textit{orthogonal} if for all formulas \index{Orthogonality}
        \[\phi(x_0,\ldots, x_{n-1}; x_0', \ldots, x_{n-1}',a)\]
        with parameters $a$ in $\mathcal{M}$, there is finitely many formulas $\theta_i(x_0,\ldots, x_{n-1},a_i)$  and $\theta_i'(x_0', \ldots, x_{n-1}',a_i'),$  with $i<k$ and parameters  $a_0,\ldots,a_{n-1}$, $a_0',\ldots,a_{n-1}'$ in $\mathcal{M},$ such that 
        \[ \phi(D^n,{D'}^n,a)= \cup_{i<k} \theta_i(D^n,a_i)\times \theta'_i({D'}^n,a_i). \]
        \end{itemize}
    \end{definition}
    

    If $S$ is a sort, we use the following terminology in order to say that definable sets in $S$ can be given by formulas with parameters in $S$ and function/predicate symbols \emph{ contained} in $S$.
    \begin{definition}\label{DefinitionPureSort}
        A sort $S$ in an $\mathrm{L}$-structure $\mathcal{M}$ is called \emph{pure}  \index{Pure! sort (or unenriched sort)}or \emph{unenriched} if definable subsets of $S$ (with parameters) are given by $\restriction{\mathrm{L}}{S}(S)$-formulas where $\restriction{\mathrm{L}}{S}$ is the language restricted to function/predicated symbols which only involve $S$.
    \end{definition}
    A pure sort  $S$ can be seen as an $\restriction{\mathrm{L}}{S}$-structure on its owns. In particular, it is stably embedded. Purity of a sort $S$ is usually a simple corollary of quantifier elimination relative to $S$ and closedness of $S$ (see Fact \ref{FactClosedImpliesPureWithControlOfParameters}).

    Remark that the notion of closedness is syntactic, which is not ideal. One may indeed use another bi-interpretable language, where the sort is no longer closed, but where resplendent relative quantifier elimination still holds\footnote{This happens for instance with the residue field $k$ of an equicharacteristic 0 Henselian valued field: it is closed in the traditional 3-sorted language of valued fields but it is also natural to interpret in addition a short exact sequences of abelian groups $1\rightarrow k^\times \rightarrow \RV^\star \rightarrow \Gamma \rightarrow 0$. See Paragraph \ref{SusubsectionRVSort}}. Here is a sightly improved version of purity which can replace the notion of closedness. It is, in a certain sens, less dependent of the language.

\begin{definition}\label{DefinitionPureWithControlOfParameters}
Consider $\mathcal{M}$ a structure. An imaginary sort $\mathcal{S}=(S,\ldots)$ endowed with an interpretable structure in a language $\mathrm{L}_S$ is called \textit{pure with control of parameters}\index{Pure! sort (or unenriched sort)!with control of parameters} if every formula  $\phi_i(x_S,b)$ where $x_S$ is a tuple of $S$-variables and $b$ is a tuple of parameters in M, is equivalent to a formula $\phi_S(x_S,\pi_S(t(b)))$ where $\pi_S$ is the canonical projection onto $S$, $\phi_S$ is an $\mathrm{L}_S$-formula and $t(x)$ is a tuple of $\mathrm{L}$-terms.
\end{definition}

The following is immediate: 
\begin{fact}\label{FactClosedImpliesPureWithControlOfParameters}
    Consider $\mathcal{M}$ an $\mathrm{L}$-structure. If $S$ is a closed sort and $\mathcal{M}$ eliminates quantifier relative to $S$, then $S$ endowed with its induced structure in $\restriction{\mathrm{L}}{S}$ is pure with control of parameters. In particular, it is pure and stably embedded.
\end{fact}


\begin{proposition}\label{PropositionEquivPureRelativeQuantifierElimination}
    Let $\mathcal{M}$ be an $\mathrm{L}$-structure, and $\mathcal{S}$ an imaginary sort of arity $n$ with some interpretable structure. Assume that $\mathcal{M}$ has quantifier elimination and that $\mathcal{S}$ is a pure imaginary structure with control of parameters. 
    The (multisorted) structure $\{\mathcal{M},\mathcal{S}, \pi_S:M^n \rightarrow S\}$ in the language $ \mathrm{L}_2:=\mathrm{L} \cup \mathrm{L}_S \cup \{\pi_S:M^n \rightarrow S\}$ admits quantifier elimination relatively to $\mathcal{S}$.
\end{proposition}

As $\mathcal{S}$ is by definition a closed sort in the language $\mathrm{L}_2$, this is in fact a characterisation of purity with control of parameters:

\begin{corollary}
    An interpretable structure $\mathcal{S}$ is pure with control of parameters if and only if $\{\mathcal{M},\mathcal{S}, \pi_S:M^n \rightarrow S\}$ admits quantifier elimination relative to $\mathcal{S}$.
\end{corollary}
 
 \begin{proof}
    We use the usual back-and-forth argument. Let $\mathcal{N}$ be a $\vert M\vert$-saturated model of the theory of $\mathcal{M}$ in the language $\mathrm{L}_2$. Let $f:(A,S_A) \rightarrow (B,S_B)$ be an isomorphism between a substructure $(A,S_A)$ of $\mathcal{M}$ and a substructure $(B,S_B)$ of $\mathcal{N}$. Assume that the restriction $\restriction{f}{S}$ to $\mathcal{S}$ is elementary. We want to extend $f$ to an embedding of $\mathcal{M}$ into $\mathcal{N}$.\\
    \textbf{Step 0:} We may assume that $S_A=S_M$. \\
    Indeed, by elementarity of $\restriction{f}{S}$, there exists an isomophism $\restriction{\tilde{f}}{S}:S_M \rightarrow \restriction{\tilde{f}}{S}(S_M) \subset S_N$ extending $\restriction{f}{S}$. The union $f \cup \restriction{\tilde{f}}{S}$ is a partial isomorphism as the sort $S$ is closed. Indeed, every quantifier-free formula $\phi(a,s)$ with parameters in $(A,S_M)$ can be written of the form:
    
    \[\bigvee \phi_{\mathrm{L}}(a)\wedge \phi_{\mathcal{S}}(s,\pi_S(t(a))), \]
    
    where $ \phi_{\mathrm{L}}$ is an $\mathrm{L}$-formula, $\phi_{\mathcal{S}}$ is an $\mathrm{L}_S$-formula and $t$ is a tuples of $\mathrm{L}$-terms. As $A$ is a structure, all terms $t(a)$ are elements of $A$. It follows that $f \cup \restriction{\tilde{f}}{S}$ preserves these formulas. 
    \\
    \textbf{Step 1:} We may assume that $A=M$ and thus conclude the proof.\\
    Indeed, let $a \in M \setminus A$. We denote by $p(x)$ the quantifier free type of $a$ over $A$. We want an appropriate answer for $a$, \textit{i.e.} an element $\tilde{f}(a)$ of $\mathcal{N}$ satisfying the set of formulas:
    \[\{ \phi(x,f(b), f(s)) \ \vert \ \phi(x,b,s) \in p(x), b\in A, \ s\in S_M   \}\]
    By compactness, we need to show that it is finitely consistent.
    Consider a formula
    \[ \phi(x,b,s) \in p(x), \]
    where $b\in A$ and $s\in S_M$.
    As $S$ is pure with control of parameters, the formula 
    \[\exists x \ \phi(x,b,y_S)\]
    is equivalent to an $\mathrm{L}_S(S_M)$-formula $\psi_S(t(b),y_S)$  (with a tuple of $\mathrm{L}$-terms $t(y)$).

    The formula  
    \[  \theta(y)= \forall y_S \ \psi_S(t(y),y_S)\Leftrightarrow \exists x \ \phi(x,y,y_S) \]
    is interpreted in the language $\mathrm{L}$ by a formula $\Sigma(y)$. As $\mathcal{M}$ has quantifier elimination in the language $\mathrm{L}$, we may assume that $\Sigma(y)$ is quantifier-free. We have:
    \begin{align*}\mathcal{M} & \models \Sigma(b),\\
        \mathcal{M} & \models  \psi_S(t(b),s).
    \end{align*}
    As $f$ respects quantifier free-formula and $\restriction{f}{S_M}$ respects $\mathrm{L}_S(S_M)$-formula, we have 
    \begin{align*}\mathcal{N} \models \Sigma(f(b)),\\
        \mathcal{N} \models  \psi_S(f(t(b)),f(s)).
    \end{align*}
    
    Of course, $f(t(b))=t(f(b))$. We get: $\mathcal{N} \models \exists x \ \phi(x,f(b),f(s))$. This concludes our proof.
 \end{proof}

 As an example, we treat the question of quantifier elimination in the field of $p$-adics in a two-sorted language of valued fields. This is a well known result, but we are not aware of a reference.     
\begin{example}
    Consider the theory $T$ of the $p$-adics $\mathbb{Q}_p$ for some $p$. By Macintyre's theorem \cite{Mac76}, it admits quantifier elimination in the language $\mathrm{L}_{Mac}:= \mathrm{L}_{rings} \cup \{P_n\}_{n<\omega}$ where the predicate $P_n$ interprets the $n^{th}$-powers. 
    \begin{itemize}
        \item The value group $\Gamma$, simply considered as a set, is not pure. Indeed, the theory $T$ in the language $\mathrm{L}_{Mac} \cup \{\Gamma \}\cup \{\val \}$ (no structure on $\Gamma$) does not eliminate the quantifiers in the formula encoding the addition:
        \[ \phi(x_\Gamma,y_\Gamma,z_\Gamma) \equiv \exists x,y \in K \ \val(x)=x_\Gamma \wedge \val(y)=y_\Gamma \wedge \val(xy)=z_\Gamma,\]
        where $x_\Gamma,y_\Gamma,z_\Gamma$ are variables in $\Gamma$.
        \item By Bélair's theorem \cite[Theorem 5.1]{Bel99}, the structure $\{\mathbb{Q}_p, \mathcal{O}_n,\Gamma,\val:\mathbb{Q}_p \rightarrow \mathbb{Z},\ac_n:\mathbb{Q}_p\rightarrow \mathcal{O}_n\}$ enriched with angular components (see Paragraph \ref{Preliminaries Unramified mixed characteristic Henselian valued fields}) eliminates quantifiers in the sort for $\mathbb{Q}_p$. It results that the value group $(\Gamma,+,0,<,\infty)$ --as an imaginary sort of $\mathbb{Q}_p$ in the ring language-- is pure with control of parameters. Then, by the theorem, $T$ eliminates quantifiers relatives to $\Gamma$ in the language
        \[\mathrm{L}_2:=\mathrm{L}_{Mac} \cup \{\Gamma,<,+,0,\infty \} \cup \{\val \} \]
        \item To get full elimination of quantifiers, one only needs to eliminates quantifiers in $\{\Gamma,<,+,0 ,\infty \}$. So the theory $T$ eliminates quantifiers in the language $\mathrm{L}_{Mac} \cup \{\Gamma,<,+,P_{\Gamma,n},0,1,\infty \} \cup \{\val \}$ where $P_{\Gamma,n}$ interprets the set of values divisible by $n$. 
        
    \end{itemize}
\end{example}

\subsubsection{Burden of a theory}\label{SubsectionClassificationtheory}

In \cite{She90}, Shelah defined the notion of burden as an invariant cardinal $\kappa_{inp}$ and implicitly defined the \emph{tree property of the second kind}. A theory which does not satisfy it is called $\text{NTP}_2$. Interest in the class of $\text{NTP}_2$ theories grew after the success of stability theory and with the necessity of extending methods to unstable contexts. In \cite{CK12}, Chernikov and Kaplan studied the forking relation in $\text{NTP}_2$ theories, establishing notably that types over models fork if and only if they divide. In \cite{Che14}, Chernikov continued the study of $\text{NTP}_2$ theories, establishing in particular a criterion with indiscernible sequences and the sub-multiplicativity of the burden. 

We recall here a definition of burden, some of the results cited above and give some important lemmas required for the proof of Theorem \ref{ThmBdnHenValFieCha00}.  We will give a second definition ( slightly different) of the burden in order to formalise a convention due to Adler \cite{Adl07}.

	\begin{definition}\label{DefBurden1}
		 Let $\lambda$ be  a cardinal. For all $i<\lambda$, $\phi_i(x,y_i)$ is $\mathrm{L}$-formula where $x$ is a common tuple of free variables, $b_{i,j}$ are elements of $\mathbb{M}$ of size $\vert  y_i \vert$ and $k_i$ is a positive natural number. Finally, let $p(x)$ be a partial type. We say that $\lbrace \phi_i(x,y_i),(b_{i,j})_{j \in \omega},k_i \rbrace_{i <\lambda}$ is \textit{an inp-pattern of depth} $\lambda$ \index{Inp-pattern} in $p(x)$ if: 
		\begin{enumerate}
		\item \textbf{for all $i<\lambda$, the $i^{\text{th}}$ row is $k_i$-inconsistent}: any conjunction $\bigwedge_{l=1}^{k_i}\phi_i(x,b_{i,j_l})$ with ${j_1<\cdots< j_{k_i}<\omega}$, is inconsistent.
		\item \textbf{all (vertical) paths are consistent}: for every $f: \lambda \rightarrow \omega$, the set $\lbrace \phi_i(x,b_{i,f(i)}) \rbrace_{i <\lambda} \cup p(x)$ is consistent.
		\end{enumerate}		   
	\end{definition}	
		Most of the time, we will not mention the $k_i$'s and only say that the rows are finitely inconsistent.
		\begin{definition}\label{DefInpMin}
		    
			\begin{itemize}
			     \item Let $p(x)$ be a partial type. \textit{The burden} \index{Burden!of a type} of $p(x)$, denoted by $\bdn(p(x))$\nomenclature[C]{$\bdn(p(x))$,$\bdn(T)$}{}, is the cardinal defined as the supremum of the depths of inp-patterns in $p(x)$. If $C$ is a small set of parameters, we write $\bdn(a/C)$ instead of $\bdn(\tp(a/C))$.
				\item The cardinal $\sup_{S \in \mathcal{S}}\bdn( \lbrace x_S=x_S \rbrace)$ where $x_S$ is a single variable from the sort $S$ and $\mathcal{S}$ is the set of sorts, is called the \textit{burden}\index{Burden!of a theory} of the theory $T$, and it is denoted by $\kappa^1_{\inp}(T)$ or by $\bdn(T)$. The theory $T$ is said to be \textit{inp-minimal}\index{Inp-minimality} if $\kappa^1_{\inp}(T)=1$.
				\item More generally, for $\lambda$ a cardinal, we denote by $\kappa^{\lambda}_{\inp}(T)$ \nomenclature[C]{$\kappa^{\lambda}_{\inp}(T)$}{} the supremum of $\bdn(\lbrace x=x\rbrace)$ where $\vert x\vert= \lambda$ and variables run in all sorts $S\in \mathcal{S}$. We always have $\kappa^{\lambda}_{\inp}(T)\geq \lambda \cdot \kappa^{1}_{inp}(T)$. In particular, if models of $T$ are infinite, $\kappa^{\lambda}_{\inp}(T)\geq \lambda$.
				\item A formula $\phi(x,y)$ has $\TP_2$ \index{Three property of the second kind} \nomenclature[C]{$\TP_2$}{Three property of the second kind} if there is an inp-pattern of the form $\lbrace \phi(x,y),(b_{i,j})_{j < \omega},k_i \rbrace_{i < \omega}$. Otherwise, we say that $\phi(x,y)$ is $\NTP_2$\index{Non Three Property of the second kind ($\NTP_2$)}\nomenclature[C]{$\NTP_2$}{Non three property of the second kind}.
				\item The theory $T$ is said $\NTP_2$ if $\kappa^1_{\inp}(T)<\infty$. Equivalently, $T$ is $\NTP_2$ if and only if there is no $\TP_2$ formula. (See \cite[Remark 3.3]{Che14})

			\end{itemize}		
		\end{definition}	
		In all the notation above, we may replace $T$ by $\mathcal{M}$.
		In \cite{Che14}, Chernikov proves the following:
		\begin{fact}[Sub-multiplicativity]
		    Let $a_1,a_2 \in \mathbb{M}$. If there is an inp-pattern of depth $\kappa_1 \times \kappa_2$ in $\tp(a_1a_2/C)$, then either there is an inp-pattern of depth $\kappa_1$ in $\tp(a_1/C)$ or there is an inp-pattern of depth $\kappa_2$ in $\tp(a_2/a_1C)$.
		\end{fact}
		
		As a corollary, for $n<\omega$, we have $\kappa_{inp}^n(T)+1 \leq (\kappa_{inp}^1(T)+1)^n$ and then $\kappa_{inp}^n(T)=\kappa_{inp}^1(T)=\bdn(T)$ as soon as one of these cardinal is infinite.
		
		If the reader knows the notion of dp-rank of a theory $T$, usually denoted by $\text{dp-rank}(T)$, let us say the following: it admits as well a similar definition in term of depth of \emph{ict-patterns} and it has been showed that a theory $T$ is NIP if and only if the depth of $ict$-pattern is bounded by some cardinal. In this paper, the reader only needs to know that the notions of $\text{dp-rank}$ and burden coincide in NIP theories. If they are not familiar with the notion of $\text{dp-rank}$, they may take it as a definition.
		\begin{fact}[{\cite[Proposition 10]{Adl07}}]
		    Let $T$ be an NIP theory, and $p(x)$ a partial type. Then $\text{dp-rank}(p(x))=\bdn(p(x))$.
		\end{fact}
		The previous fact is only stated with partial type $p(x)=\{x=x\}$, but the proof is the same.
		
		\begin{example}
		    \begin{itemize}
		        
		        \item Any quasi-o-minimal theory is inp-minimal (see e.g. \cite[Theorem A.16]{Sim15}). In particular, $\{\mathbb{Z}, 0, +, <\}$ is inp-minimal.
		        \item Let $\mathrm{L}=\{R,B\}$ be the language with two binary predicates, and let $\mathcal{M}$ be a set with two cross-cutting equivalence relations with infinitely many infinite classes (for all $a$ and $b$, there is infinitely many $c$ such that $aRc$ and $bBc$). 
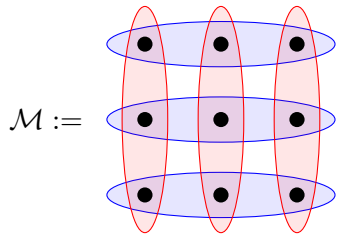
\begin{center} 
\begin{tikzpicture}\clip(-2,-0.7) rectangle (5,2.7);
\draw (-1.3,1) node {$\mathcal{M}:=$};

\draw[blue, fill=blue!10] (1,0) ellipse (1.5 and 0.3);
\draw[blue, fill=blue!10] (1,1) ellipse (1.5 and 0.3);
\draw[blue, fill=blue!10] (1,2) ellipse (1.5 and 0.3);

\fill[red,opacity=0.1] (0,1) ellipse (0.3 and 1.5);
\draw[red] (0,1) ellipse (0.3 and 1.5);
\fill[red,opacity=0.1] (1,1) ellipse (0.3 and 1.5);
\draw[red] (1,1) ellipse (0.3 and 1.5);
\fill[red,opacity=0.1] (2,1) ellipse (0.3 and 1.5);
\draw[red] (2,1) ellipse (0.3 and 1.5);

\fill (0,0) circle (0.1);
\fill (1,0) circle (0.1);
\fill (2,0) circle (0.1);
\fill (0,1) circle (0.1);
\fill (1,1) circle (0.1);
\fill (2,1) circle (0.1);
\fill (0,2) circle (0.1);
\fill (1,2) circle (0.1);
\fill (2,2) circle (0.1);
\end{tikzpicture}
\end{center}
One proves easily that $\bdn(\mathcal{M})=2$. For $\lambda$ a cardinal, one can consider $\lambda$-many cross cutting equivalence relations, and shows that the structure is of burden $\lambda$.
		            
		    \end{itemize} 
		\end{example}

	The definition of burden of a theory, as many other notion of complexity, gives to unary sets an important role. But one has to notice that the notion of unary set is syntactic, and is not preserved under bi-interpretability:
	
	\begin{remark}\label{k to the n kappa n}
        For $n\in \mathbb{N}$, one can consider the multisorted structure $(\mathcal{M}^n, \mathcal{M}, p_i,i< n )$ where $p_i: \mathcal{M}^n \rightarrow \mathcal{M}, (a_0,\ldots, a_{n-1}) \mapsto a_i$ is the projection to the $i^{\text{th}}$ coordinate. If we denote its theory by $T^n$, then we clearly have $\kappa_{inp}^n(T)=\kappa_{inp}^1(T^n)$. 
    \end{remark}

	To clarify, let us introduce the following terminology:
	\begin{definition}
        Let $\mathcal{M}$ and $\mathcal{N}$ be two structures. We say that $\mathcal{N}$ is \emph{interpretable on a unary set} \index{Interpretability on a unary set} in $\mathcal{M}$ if there is a bijection $f:N \rightarrow D/ \sim$ where $D$ is a unary definable set in $\mathcal{M}$, $\sim$ is a definable equivalence relation, and the pull-back in $\mathcal{M}$ of any graph of function and relation of $\mathcal{N}$ is definable. The structures  $\mathcal{M}$ and $\mathcal{N}$ are said \emph{bi-interpretable on unary sets} if $\mathcal{N}$ is interpretable on a unary set in $\mathcal{M}$ and $\mathcal{M}$ is interpretable on a unary set in $\mathcal{N}$.
    \end{definition}

    We will work up to bi-interpretability on unary sets, meaning in particular that the main results of this text will only depend on the structure that we want to consider and not on the language.  
    
    \begin{fact}\label{FactBdnInterpretUnarySet}
        Let $\mathcal{M}$ and $\mathcal{N}$ be two structures, and assume that $\mathcal{N}$ is interpretable on a unary set in $\mathcal{M}$, then $\bdn(\mathcal{N}) \leq \bdn(\mathcal{M})$. In particular, if $\mathcal{M}$ and $\mathcal{N}$ are bi-interpretable on unary sets, then $\bdn(\mathcal{M})=\bdn(\mathcal{N})$.
    \end{fact}
    For example, $\{\mathbb{Z},0,+,< \}$ does not interpret 
    \[\{\mathbb{Z}\times \mathbb{Z},(\mathbb{Z}, 0, +, <), \pi_1:\mathbb{Z}\times \mathbb{Z} \rightarrow \mathbb{Z} ,\pi_2:\mathbb{Z}\times \mathbb{Z} \rightarrow \mathbb{Z} \}\]
    on a unary set (the first being of burden 1, the second being of burden at least 2).
    However, if $k$ is an imperfect field, we will see that $k$ interprets $\{k\times k,(k,0,1,+,\cdot), \pi_1:k\times k \rightarrow k ,\pi_2:k \times k \rightarrow k \}$ on a unary set.

	\paragraph{Lemmas on inp-patterns}
	Let $\mathrm{L}$ be any first order language, $\mathcal{M}$ a $\mathrm{L}$-structure of base set $M$ and let $\lambda$ be a cardinal. 
	
	\begin{definition}
	     \begin{itemize}
	        \item A sequence $(b_j)_{j\in\lambda}$ of (tuples of) elements of $M$ is indiscernible over a subset $A\subset M$ if for every $n\in \mathbb{N}$ and every formula $\phi(x_0,\ldots, x_{n-1}, a)$ with parameters $a\in A$, we have
	        \[\mathcal{M} \models \phi(b_0, \ldots,b_{n-1},a) \Leftrightarrow \phi(b_{j_0}, \ldots,b_{j_{n-1}},a)\]
	        for every $j_0<\cdots<j_{n-1} \in \lambda$.
            \item An array $(b_{i,j})_{i\in \lambda,j \in \omega}$ is mutually indiscernible if every line $(b_{i,j})_{j \in \omega}$ is indiscernible over $\{b_{k,j}\}_{k\neq i, k\in \lambda, j<\omega}$.
        \end{itemize}

	\end{definition}
	
	We will intensively use the following fact:
	\begin{fact}[{\cite[Lemma 2.2]{Che14}}]\label{LemmaWMAArrMutInd}
	    If $p(x)$ is a partial type and if $\lbrace \phi_i(x,y_i),(b_{i,j})_{j \in \omega},k_i \rbrace_{i <\lambda}$ is an inp-pattern in $p(x)$, there is an inp-pattern $\lbrace \phi_i(x,y_i),(\tilde{b}_{i,j})_{j \in \omega},k_i \rbrace_{i <\lambda}$ in $p(x)$ with a mutually indiscernible array $(\tilde{b}_{i,j})_{i< \lambda,j < \omega}$.
	\end{fact}
	
	We will now present some easy lemmas, which we will later use. They give us tools to `transform' inp-patterns into simpler ones which are easier to analyse.
		
    \begin{lemma} \label{INPequiv}
    Let $\lbrace \phi_i(x,y_{i}), (a_{i,j})_{j<\omega},k_i\rbrace_{i<\lambda}$ be an inp-pattern with $(a_{i,j})_{i<\lambda,j<\omega}$ mutually indiscernible. Assume for every $i<\lambda$, $\phi_i(x,a_{i,0})$ is equivalent to some formula $\psi_{i}(x,b_{i,0})$ with parameter $b_{i,0}$. Then we may extend $(b_{i,0})_{i <\lambda}$ to an mutually indiscernible array $(b_{i,j})_{i < \lambda, j < \omega}$ such that 
    $$\lbrace \psi_{i}(x,y_{i}), (b_{i,j})_{j<\omega},k_i\rbrace_{i<\lambda},$$
    is an inp-pattern.
    \end{lemma}
    \begin{proof}
    By 1-indiscernibility, we find $b_{i,j}$ such that $\phi_i(\mathbb{M},a_{i,j}) = \psi_i(\mathbb{M},b_{i,j})$. Then, the statement is clear.
    \end{proof}
    
    \begin{remark}\label{rkSEDS}
        Let $D$ be a stably embedded definable set in $\mathbb{M}$, and $\lbrace \phi_i(x,y_{i,j}), (a_{i,j})_{j<\omega},k_i\rbrace_{i<\lambda}$ an inp-pattern in $D$. This in particular implies that solutions of paths can be found in $D$ but the parameters $(a_{i,j})$ may not belong to $D$. Using the previous lemma, we may actually assume that this is the case. It follows that $D$ endowed with the induced structure is at least of burden $\lambda$.  
    \end{remark}
    
    The next lemma shows that one can `eliminate' disjunction symbols in inp-patterns. A direct consequence is that if the theory has quantifier elimination, then we may assume that formulas of inp-patterns are conjunctions of atomic and negation of atomic formulas. 
    
    \begin{lemma}\label{INPdisj}
    Let $\lbrace \phi_i(x,y_{i,j}), (a_{i,j})_{j<\omega},k_i\rbrace_{i<\lambda}$ be an inp-pattern with $(a_{i,j})_{i<\lambda,j<\omega}$ mutually indiscernible. Assume  that $\phi_i(x,y_{i,j})= \bigvee_{l\leq n_i} \psi_{l,i}(x,y_{i,j})$. Then there exists a sequence of natural numbers $(l_i)_{i < \lambda}$ such that $l_i \leq n_i$ and  
    $$\lbrace \psi_{l_i,i}(x,y_{i,j}), (a_{i,j})_{j<\omega},k_i\rbrace_{i<\lambda}$$
    is an inp-pattern.
    \end{lemma}

    \begin{proof}
     Let $d \models \lbrace \phi_i(x,a_{i,0})\rbrace_{j<\lambda}$. For every $i< \lambda$, let $l_i \leq n_i$ be such that $d \models \psi_{l_i,i}(x,a_{i,0})$. By the mutual-indiscernibility of $(a_{i,j})_{i<\lambda,j<\omega}$, every path of the pattern $\lbrace \psi_{l_i,i}(x,y_{i,j}), (a_{i,j})_{j<\omega},k_i\rbrace_{i<\lambda}$ is consistent. The inconsistency of the rows follows immediately from the inconsistency of the rows of the initial pattern.
    \end{proof}
    
    `Elimination' of conjunction symbols may happen in more specific context. Notably:
   \begin{proposition}\label{sumburden}
    Let $\mathcal{K}$ and $\mathcal{H}$ be two structures, and consider the multisorted structure $\mathcal{G}$:
    \[\mathcal{G}=\{ K \times H, \mathcal{K},\mathcal{H}, \pi_K: K\times H \rightarrow K, \pi_H: K\times H \rightarrow H \},\]
    called the \emph{direct product structure} (where $\pi_K$ and $\pi_H$ are the natural projections). Then $\mathcal{G}$ eliminates quantifiers relative to $\mathcal{K}$ and $\mathcal{H}$, and $\mathcal{K}$ and $\mathcal{H}$ are orthogonal and stably embedded within $\mathcal{G}$.\index{Product of structures} 
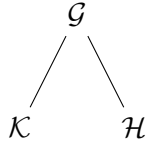
    \begin{center}
        \begin{tikzpicture}
        \node{$\mathcal{G}$ }
    child {  node {$\mathcal{K}$}}
    child { node {$\mathcal{H}$} };
    \end{tikzpicture}
    \end{center}
    
    We have $$\bdn(\mathcal{G}) = \bdn(\mathcal{K})+ \bdn(\mathcal{H}).$$
\end{proposition}
       
We prove an obvious generalisation for product of more than two structures in the next paragraph.  

\begin{proof}
Relative quantifier elimination, stable embeddedness and orthogonality are rather obvious.
The inequality $\bdn(\mathcal{G}) \geq \bdn(\mathcal{K}) + \bdn(\mathcal{H})$ is easy but we give a detailed proof. Let $\lbrace \phi_i(x_K,y_i),(a_{i,j})_{j<\omega} \rbrace_{i \in \lambda_1}$ be an inp-pattern in $\mathcal{K}$ and $\lbrace \psi_i(x_H,y_i),(b_{i,j})_{j<\omega} \rbrace_{i\in \lambda_2}$ an inp-pattern in $\mathcal{H}$. Then 
$$\lbrace \phi_i(\pi_K(x_{K}, x_{H}),y_i),(a_{i,j})_{j<\omega} \rbrace_{i \in \lambda_1} \cup \lbrace \psi_i(\pi_H(x_{K}, x_{H}),y_i),(b_{i,j})_{j<\omega} \rbrace_{i\in \lambda_2}$$
is an inp-pattern in $\mathcal{G}$ of depth $\lambda_1  + \lambda_2$. Indeed, first notice that inconsistency of each rows is clear. Secondly, take a path ${f: \lambda_1 \sqcup \lambda_2 \rightarrow \omega} $. There is an element $d_{K}\in K$ satisfying $\lbrace \phi_i(x_K,a_{i,f(i)})\rbrace_{i\in \lambda_1}$ and an element $d_H \in H$ satisfying $\lbrace \psi_i(x_H,b_{i,f(i)})\rbrace_{i\in \lambda_2}$. Then, the element $d=(d_K,d_H)$ of $\mathcal{G}$ is a solution of the pattern along the path $f$.

For the other inequality, let $\lbrace \theta_i(x,y_{i,j}), (c_{i,j})_{j<\omega},k_i\rbrace_{i<\lambda}$ be an inp-pattern in $\mathcal{G}$, with $(c_{i,j})_{i<\lambda,j<\omega}$ mutually indiscernible.
We may assume $\theta_i(x,c_{i,j})$ is of the form $\phi_{i}(x_K, a_{i,j}) \wedge \psi_{i}(x_H, b_{i,j})$ where $x_K= \pi_K(x)$, $x_H=\pi(x_H)$, $c_{i,j}=a_{i,j}\upwedge b_{i,j}$, $\phi_i(x_K,a_{i,j})$ is a $\mathcal{K}$-formula and $\psi_i(x_H,b_{i,j})$ is a $\mathcal{H}$-formula. Indeed, let $d \models \lbrace \theta_i(x,c_{i,0})\rbrace_{j<\lambda}$, by orthogonality, $\theta_i(x,c_{i,0})$ is equivalent to a formula of the form:
$$ \bigvee\limits_{k<n_i}\phi_{i,k}(x_K, a_{i,0}) \wedge \psi_{i,k}(x_H, b_{i,0}).$$
Then we conclude by using Lemmas \ref{INPequiv} and \ref{INPdisj}. For every $i$, at least one of the sets $\lbrace \phi_i(x_K, a_{i,j})\rbrace_{j<\omega}$ and $\lbrace \psi_i(x_H, b_{i,j})\rbrace_{j<\omega}$ is $k_i$-inconsistent (by indiscernibility of $(c_{i,j})_j$). We may "eliminate" the conjunction as well and assume that every line is an $\mathrm{L}_K$-formula or an $\mathrm{L}_H$-formula. We conclude that $\lambda \leq \bdn(\mathcal{K}) + \bdn(\mathcal{H})$.
 \end{proof}

   Together with Fact \ref{FactBdnInterpretUnarySet}, we get more generally:
    \begin{fact}\label{FactSumBurdenOrthSorts}
        Let $M=(A,C,\ldots )$ be a many-sorted structure. Assume that $A$ and $C$ are orthogonal and stably embedded in $M$. Then we have 
        $\bdn(A \times C) = \bdn(A)+ \bdn(C).$
    \end{fact}

    
    Let us finish this paragraph with one more lemma:

    \begin{lemma}\label{Ft1Bdn}
Let $D$ and $D'$ two type-definable sets respectively given by the partial types $p(x)$ and $p'(x)$ and let $f: D \rightarrow D'$ be a surjective finite to one type-definable function.  Then we have $\bdn(D):= \bdn(p(x))= \bdn(p'(x)) =: \bdn(D')$. 
    \end{lemma}

    \begin{proof}
We may assume that $D$ and $D'$ are definable, the general case can be similarly deduced. Let $\lbrace \phi_i'(x',y_i),(a_{i,j})_{j< \omega},k_i \rbrace_{i < \lambda}$ be an inp-pattern in $D'$. Clearly, $\lbrace \phi_i'(f(x),y_i),(a_{i,j})_{j< \omega},k_i \rbrace_{i <\lambda}$ is an inp-pattern in $D$. Hence $\bdn(D) \geq \lambda$. Conversely, let $\lbrace \phi_i(x,y_i),(a_{i,j})_{j< \omega},k_i \rbrace_{i < \lambda}$ be an inp-pattern of depth $\lambda$ in $D$. Consider the pattern
$$\lbrace \phi_i'(x',y_{i,j}),(a_{i,j})_{j< \omega}\rbrace_{i < \lambda},$$ 
where 
$$\phi_i'(x',a_{i,j}) \quad \equiv \quad \exists x  \ x \in D \wedge x'=f(x) \wedge \phi(x,a_{i,j}).$$ 

Clearly every path is consistent. Assume for some $i <\lambda $, the row $\lbrace \phi_i'(x',a_{i,j}) \rbrace_{j< \omega}$ is consistent, witnessed by some $h'$. Note that $h'$ is in $D'$. By the pigeonhole principle, there is $h \in D$ and an infinite subset $J$ of $\omega$ such that $f(h)= h'$ and $h \models \lbrace \phi_i(x,a_{i,j})\rbrace_{j \in J}$, contradiction. It follows that $\lbrace \phi_i(x',y_{i,j}), (a_{i,j})_{j< \omega} \rbrace_{i <\lambda}$ is an inp-pattern in $D'$. We conclude that $\bdn(D') \geq \lambda $.
    \end{proof}

\subsubsection{More on burden and strength}\label{More on burden and strongness}
	
     We will formally introduce a well known convention with respect to the burden, which consists of writing $\bdn(\mathcal{M})= \lambda_-$ for a limit cardinal $\lambda$ if $\mathcal{M}$ admits inp-patterns of depth $\mu$ for all $\mu<\lambda$, but no inp-pattern of depth $\lambda$. It has been introduced in \cite{Adl07}, and has the advantage to emphasising a relevant distinction. If the reader is not interested by such subtleties, they may move to the next subsection. 
     Proposition \ref{PropInfProd} might be interesting on its own, as it corresponds to the `baby case' for the difficulty that we will encounter for mixed-characteristic Henselian valued fields. One can refer to \cite{Adl07} for this paragraph.
	\begin{definition}
        We define the ordered class $(\Card^\star, <)$\nomenclature[A]{$(\Card^\star, <)$}{} as the linear order obtained from the ordered class of cardinals $(\Card, <)$ by adding for any limit cardinal $\lambda$ a new element $\lambda_-$ (called `lambda minus'). This new element comes immediately before $\lambda$: $\lambda_- < \lambda$ and if $\mu \in \Card^\star$ with $\mu < \lambda$, then $\mu \leq \lambda_-$. In addition to the natural injection $\Card \hookrightarrow \Card^\star$, we define the \textit{actualisation map} $\act:\Card^\star \rightarrow \Card$ \nomenclature[A]{$\act(\lambda)$ for $\lambda \in \Card^\star$ }{Actualisation map} as the map such that $\act(\lambda_-)=\lambda$ for every limit cardinal $\lambda$, and $\act(\kappa)=\kappa$ for any cardinal $\kappa\in \Card$. It will be convenient to also set $\kappa_-=\lambda$ when $\kappa= \lambda^+\in \Card$ is a successor cardinal.
        \end{definition}
        
        If $\lambda$ is a limit cardinal, one should think $\lambda$ as an `actual' lambda and $\lambda_-$ as a `potential' lambda. We don't change our notion of cardinality of a set. As we will see, this definition of $\Card^\star$ is motivated by the burden, \textit{i.e.} by a notion of dimension. It also motivates to (partially) extend the arithmetic operations of $\Card$ to $\Card^\star$. We will have to answer any question of the form: should the cardinal $\aleph_0\cdot \aleph_{\omega-}$ be $\aleph_{\omega-}$ or $\aleph_{\omega}$? As the definitions themselves appear to be a bit technical, we prefer to first give intuition to the reader with a small digression on graphs.
    
    \paragraph{Graphs and cliques}
        We consider symmetric graphs in the language $\mathrm{L}=\{R\}$. We denote by $K_\kappa$ the complete graph on $\kappa$-many vertices, for $\kappa$ a cardinal in $\Card$. Given a graph $\mathcal{G}$, we denote by $C(\mathcal{G})$ the cardinal in $\Card^\star$:
    \[C(\mathcal{G})=\begin{cases}
        \kappa \in \Card \text{ if $K_{\kappa}$ embeds in $\mathcal{G}$ and $K_{\kappa^+}$ does not,}\\
        \kappa_- \in \Card^\star \text{ if $K_\lambda$ embeds in $\mathcal{G}$ for all cardinal $\lambda< \kappa$ and $K_{\kappa}$ does not.} 
        \end{cases}\]
        
    \begin{example}
    Let $\mathcal{G}$ be the disjoint union of graphs $\cup_{n<\aleph_0} K_n$:
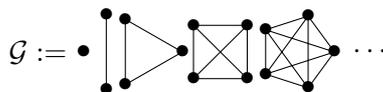
    \begin{center}
        
    \begin{tikzpicture} [scale= 0.5]

           \draw (-0.8,0) node{$\mathcal{G}:=$};
        \draw (0.4,0) node{$\bullet$};
        
        \begin{scope}[shift = {(1,0)}]
            \draw (0,1) node{$\bullet$} -- (0,-1) node{$\bullet$};
        \end{scope}
        
        \begin{scope}[shift = {(2,0)}]
            \draw (0:1) node{$\bullet$} -- (120:1) node{$\bullet$} -- (240:1) node{$\bullet$} -- cycle;
        \end{scope}
        \begin{scope}[shift = {(4,0)}]
            \draw (45:1) node{$\bullet$} -- (135:1) node{$\bullet$} -- (225:1) node{$\bullet$} -- (315:1) node{$\bullet$}-- cycle;
            \draw (45:1) -- (225:1);
            \draw (135:1) -- (315:1); 
        \end{scope}
        
        \begin{scope}[shift ={ (6,0)}]
        \draw (0:1) node{$\bullet$} -- (72:1) node{$\bullet$} -- (144:1) node{$\bullet$} -- (216:1) node{$\bullet$} -- (288:1) node{$\bullet$} -- (0:1);
        \draw (0:1) -- (144:1) -- (288:1) -- (72:1) -- (216:1) -- (0:1)  ;
        \end{scope}
        \draw (8,0) node{$\cdots$};
    \end{tikzpicture}
    
    \end{center}
    By definition, we have $C(\mathcal{G})= \aleph_{0-}$.
    \end{example}
    In addition to the union of graphs, we want to consider another natural operation:
    \begin{definition}
        We define the lexicographic product of graphs \index{Lexicographic product! of graphs} $\mathcal{G}$ and $\mathcal{F}$ as the graph $\mathcal{G}[\mathcal{F}]$ with set of vertices $G\times F$ and a symmetric relation  given by:
        \[(g_0,f_0) R^{\mathcal{G}[\mathcal{F}]} (g_1,f_1) \ \Leftrightarrow 
        \begin{cases}
        g_0=g_1 \ \text{ and } f_0 R^{\mathcal{F}} f_1, \\
        g_0 \neq g_1 \ \text{ and } g_0 R^{\mathcal{G}} g_1.
        
        \end{cases}\]
    \end{definition}

\begin{example}
    Consider the lexicographic product of $K_4$ and $K_3$. We simply obtain $K_{12}$:
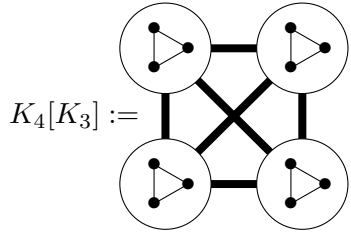
\begin{center} 
\begin{tikzpicture}[scale=0.6]\clip(-4,-1) rectangle (10,5);
\draw[line width=3pt] (0,0)--(0,3);
\draw[line width=3pt] (0,0)--(3,0);
\draw[line width=3pt] (0,0)--(3,3);
\draw[line width=3pt] (3,0)--(3,3);
\draw[line width=3pt] (0,3) -- (3,3);
\draw[line width=3pt] (0,3) -- (3,0);
\draw (-2,1.5) node {$
K_4 [ K_3]:=$};
\begin{scope}[shift={(0,0)}] 
\filldraw[fill=white] (0,0) circle (1);
 \draw (0:0.5) node{$\bullet$} -- (120:0.5) node{$\bullet$} -- (240:0.5) node{$\bullet$} -- cycle;

\end{scope}
\begin{scope}[shift={(0,3)}]
\filldraw[fill=white] (0,0) circle (1);
 \draw (0:0.5) node{$\bullet$} -- (120:0.5) node{$\bullet$} -- (240:0.5) node{$\bullet$} -- cycle;
 \end{scope}
\begin{scope}[shift={(3,0)}]
\filldraw[fill=white] (0,0) circle (1);
 \draw (0:0.5) node{$\bullet$} -- (120:0.5) node{$\bullet$} -- (240:0.5) node{$\bullet$} -- cycle;
\end{scope}
\begin{scope}[shift={(3,3)}]
\filldraw[fill=white] (0,0) circle (1);
 \draw (0:0.5) node{$\bullet$} -- (120:0.5) node{$\bullet$} -- (240:0.5) node{$\bullet$} -- cycle;
\end{scope}
\end{tikzpicture}
\end{center}

\end{example}

    If $C(\mathcal{G}), C(\mathcal{F}) \in \Card$ are cardinals greater or equal to 2, we have by pigeonhole principle that \[C(\mathcal{G}[\mathcal{F}])= C(\mathcal{G}) \times C(\mathcal{F}).\] 
    This gives us the intuition of how one can define the product of cardinals in $\Card^\star$. Let us look at two examples:
\begin{example}
    Consider $\mathcal{G} = \cup_{n<\omega} K_{\aleph_n}$ and  $\mathcal{F}= \cup_{\alpha<\omega_1} K_{\aleph_\alpha}$. Then
    we have $C( \mathcal{G})= \aleph_{\omega-}$ \nomenclature[A]{$\aleph_{0-}$, $\aleph_{\omega-}$, etc}{} and  $C( \mathcal{F})= \aleph_{\omega_1-}$. If we consider the lexicographic product with $K_{\aleph_0}$, we obtain:
    \begin{itemize}
        \item $C( K_{\aleph_0}[\mathcal{G}])= \aleph_\omega \ $ ,
        \item $C( K_{\aleph_0}[\mathcal{F}])= \aleph_{\omega_1-} \ $ .
    \end{itemize}
 We leave the proof to the reader, with the following picture for the intuition:
 
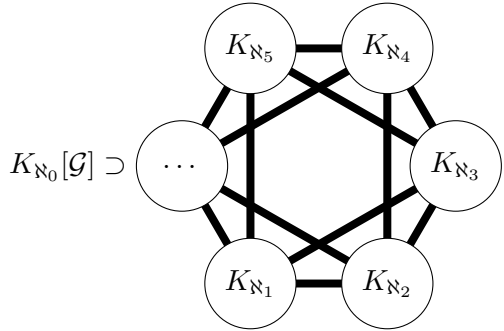
\begin{center} 
\begin{tikzpicture}[scale=0.6]\clip(-7,-4) rectangle (9,4);
\draw[line width=3pt] (0:3)--(60:3) -- (120:3)-- (180:3) -- (240:3) -- (300 :3) -- cycle;
\draw[line width=3pt] (0:3)--(120:3)--(240:3) -- cycle;
\draw[line width=3pt] (60:3)-- (180:3) -- (300 :3) -- cycle;
\draw (-5.5,0) node {$
K_{\aleph_0}[\mathcal{G}] \supset$};
\begin{scope}[shift={(0:3)}] 
\filldraw[fill=white] (0,0) circle (1);
 \draw (0:0) node{$K_{\aleph_3}$};
\end{scope}
\begin{scope}[shift={(60:3)}] 
\filldraw[fill=white] (0,0) circle (1);
 \draw (0:0) node{$K_{\aleph_4}$};
\end{scope}
\begin{scope}[shift={(120:3)}] 
\filldraw[fill=white] (0,0) circle (1);
 \draw (0:0) node{$K_{\aleph_5}$};
\end{scope}
\begin{scope}[shift={(240:3)}] 
\filldraw[fill=white] (0,0) circle (1);
 \draw (0:0) node{$K_{\aleph_1}$};
\end{scope}
\begin{scope}[shift={(300:3)}] 
\filldraw[fill=white] (0,0) circle (1);
 \draw (0:0) node{$K_{\aleph_2}$};
\end{scope}

\begin{scope}[shift={(180:3)}] 
\filldraw[fill=white] (0,0) circle (1);
 \draw (0:0) node{$\dots$} ;
\end{scope}

\end{tikzpicture}
\end{center}
 
 As a consequence, one might be tempted to write $\aleph_{0} \cdot \aleph_{\omega-}= \aleph_{\omega}$ and $\aleph_{0} \cdot \aleph_{\omega_1-}= \aleph_{\omega_1-}$. This is what we want to define now.

\end{example}

    \paragraph{Arithmetic on $\Card^\star$}
    We first define the cofinality of a cardinal $\lambda$ in $\Card^\star$ as the cofinality of $\act(\lambda)$, denoted by $\cf(\lambda)$. Secondly, we define the following operations:

    \begin{definition}\label{DefOpeCarEto}
         Let $\Lambda=(\lambda_i)_{i\in I}$ be a sequence in $\Card^\star$. Let $\lambda=\sup_{i\in I}(\act(\lambda_i))  \in \Card$ be the supremum in the usual sense, and $\supp(\Lambda)=\lbrace i\in I \ \vert \lambda_i\neq 0\rbrace$.
         We find a partition $I_1\cup I_2 \cup I_3$ of $I$ such that:
        \[\Lambda = (\lambda_i)_{i\in I_1}\cup (\lambda_-)_{i\in I_2} \cup (\lambda)_{i\in I_3},\]
        where $\lambda_i < \lambda_-$ for $i\in I_1$.
        We define $\sup^\star$ \nomenclature[A]{$\sup^\star$}{} as follows:
        \begin{itemize}
    
            \item  $\sup^\star_{i\in I}(\lambda_i)= 
            \begin{cases}
                \lambda \text{ if }\vert I_3\vert\neq \emptyset.\\
                \lambda_- \text{ otherwise}.
            \end{cases}$
        \end{itemize}
        If $\vert I \vert$ and $\lambda$ are finite, the definition of the sum $\sum^\star$\nomenclature[A]{$\sum^\star$}{} in $\Card^\star$ is the sum in the usual sense: \[{\sum_{i\in I}}^\star\lambda_i=\sum_{i\in I}\act(\lambda_i)=\sum_{i\in I}\lambda_i.\] Otherwise, we set:
        \begin{itemize}
            \item $\sum^\star_{i\in I} \lambda_i =
            \begin{cases}
                \vert \supp(\Lambda) \vert \text{ if $\vert \supp(\Lambda) \vert \geq \lambda$},\\
                \begin{cases}
                    \lambda \text{ if $I_3 \neq \emptyset$}, \\
                    \lambda \text{ if $\vert I_2 \vert \geq \cf(\lambda)$},\\
                    \lambda \text{ if $\sup_{i\in I_1}(\act(\lambda_i))=\lambda$}, \\
                    \lambda_- \text{ otherwise.}
                \end{cases} \text{ if $\vert \supp(\Lambda) \vert < \lambda$}.
            \end{cases}$
        \end{itemize}
        For $\lambda, \mu \in \Card$, $\lambda$ limit cardinal, we define the product $\cdot^\star: \Card\times \Card^\star \rightarrow \Card^\star$ in terms of sum: 
            \begin{itemize}
                \item $\sum_{\mu}^\star \lambda_- = \mu\cdot^\star \lambda_- = 
                    \begin{cases}
                        \lambda_- \text{ if $\mu < \cf(\lambda)$,}\\
                        \mu \cdot \lambda \text{ if $\mu \geq \cf(\lambda)$.}
                    \end{cases}$
            \end{itemize}

    \end{definition}
    We see in particular that, under these definition, $\sup^\star$ and $\sum^\star$ do not necessary coincide anymore when there are infinite. However, it is clear that we recover the usual definition via the actualisation map:  
    \[\xymatrix{
{\Card^\star}^{\vert I \vert} \ar[d]^{\act} \ar[r]^{\sum^\star} & \Card^\star \ar[d]^{\act}\\
\Card^{\vert I \vert} \ar[r]^{\sum}          &\Card}
\qquad \qquad
\xymatrix{
{\Card^\star}^{\vert I \vert} \ar[d]^{\act} \ar[r]^{\sup^\star} & \Card^\star \ar[d]^{\act}\\
\Card^{\vert I \vert} \ar[r]^{\sup}          &\Card}    \]

    Here are the promised examples: 
    \begin{examples}
        \begin{itemize}
            \item Consider the sequence $\Lambda_1=\aleph_{\omega-},1,2,3,\ldots$. We have  $\sup^\star \Lambda_1 =\sum^\star \Lambda_1= \aleph_{\omega-}$.
            \item Consider the sequences $\Lambda_2=(\aleph_{\omega-})_{i< \omega}=\aleph_{\omega-}, \aleph_{\omega-}, \ldots $ and 
        $\Lambda_3=(\aleph_i)_{i<\omega}=\aleph_{0}, \aleph_{1}, \ldots $. We have $\sup^\star \Lambda_2=\sup^\star \Lambda_3 = \aleph_{\omega-}$ and $\sum^\star \Lambda_2=\sum^\star \Lambda_3 = \aleph_\omega$.
        \item Consider $\Lambda_4= (\aleph_i)_{i<\omega}\cup (\aleph_{2\omega-})$. Then $\sup^\star \Lambda_4=\sum^\star \Lambda_4 =\aleph_{2\omega-}$.
        \item We have $\aleph_0 \cdot \aleph_{\omega-}= \aleph_{\omega}$, $\aleph_0 \cdot \aleph_{\omega_{1}-}=\aleph
        _{\omega_{1}-}$ and $\aleph_1 \cdot \aleph_{\omega_{1}-}=\aleph_{\omega_{1}}$.
        \end{itemize}
    \end{examples}

    Now, we go back to the burden. 
    \paragraph{Burden, strength and $\Card^\star$}
    In Definition \ref{DefInpMin}, the burden of the complete theory $T$ is the supremum (in $\Card \cup \{\infty\}$) of depth of inp-patterns in $T$. However this supremum is not necessarily attained by an actual inp-pattern. This distinction is in particular motivated by the following definition: 
    \begin{definition}[\cite{Adl07}]
        A complete theory is called \textit{strong} if there is no inp-pattern of infinite depth in $T$.
    \end{definition} 
    
    One sees that, paradoxically, some strong theories have burden $\aleph_0$ and some theories of burden $\aleph_0$ are not strong (see examples below). In other words, the definition of burden we gave failed to characterize strength.  We will follow Adler's convention (see \cite{Adl07}) which gives a solution to this problem: burden will now take values in $\Card^\star \cup \lbrace \infty \rbrace$. 
    \begin{definition}(second definition of burden)\label{DefBurden2}
        Let $T$ be a complete theory. We denote by $\mathcal{S}$ the set of sorts.
		\begin{itemize}
		     \item The burden $\bdn(\pi(x))$ of a partial type $\pi(x)$ is the supremum in $\Card^\star\cup \lbrace \infty \rbrace$ of the depths of inp-patterns in $p(x)$.\index{Burden!of a type}
			\item The cardinal $\sup_{S \in \mathcal{S}}^\star\bdn( \lbrace x_S=x_S \rbrace)$ where $x_S$ is a single variable from the sort $S$, is called the \textit{burden}\index{Burden!of a theory} of the theory $T$, and it is denoted by $\kappa^1_{\inp}(T)$ or by $\bdn(T)$.
		\end{itemize}
	\end{definition}

    In other words, if the supremum $\lambda \in \Card$ of depth of inp-patterns is attained, the burden is equal to $\lambda$ . Otherwise, it is equal to $\lambda_-$. In particular, strong theories are exactly theories of burden at most $\aleph_{0-}$.
    One can check that every lemma in the previous paragraph -and its proof- still hold. Let us give a formal definition:
    
    \begin{definition}
        Let $\mathcal{M}_i=(M_i,\ldots)$ be a structure in a language $\mathrm{L}_i$, for $i\in I$ a set of indices. We define the following multisorted structure:
        \begin{itemize}
            \item The \textit{disjoint union}  
                    \[\bigcup_i \mathcal{M}_i  = \lbrace (M_i,\ldots) \rbrace_{i\in I}.\] \nomenclature[]{$\bigcup_i \mathcal{M}_i$}{Disjoint union of structures $\mathcal{M}_i$}
                    with a sort for each $\mathcal{M}_i$'s.
            \item The \textit{direct}  \textit{product}         \[\prod_{i \in I} \mathcal{M}_i = \lbrace \prod_{i \in I} M_i, (M_i,\dots)_{i\in I}, (\pi_i: \prod_{j \in I} M_j \rightarrow M_i)_{i\in I}\rbrace,\] \nomenclature[]{$\prod_i \mathcal{M}_i$}{Direct product of structures $\mathcal{M}_i$}
        with a sort for each $\mathcal{M}_i$'s and a sort for the product and where $\pi_i: \prod_{j \in I} M_j \rightarrow M_i$ is the natural projection.
        \end{itemize}
    \end{definition}
        We have the following fact: 
        \begin{fact}\label{FactQEProd}
        \begin{itemize}
            \item  The sorts $\mathcal{M}_i$ in the union $\cup_{i\in I}\mathcal{M}_i$ are stably embedded and pairwise orthogonal.
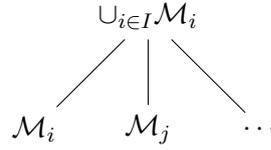
            \begin{center}
                \begin{tikzpicture}
                    \node{$\cup_{i\in I}\mathcal{M}_i$ }
                        child {  node {$\mathcal{M}_i$}}
                        child { node {$\mathcal{M}_j$}} 
                        child { node {$\cdots$}
                        };
                \end{tikzpicture}
            \end{center}

            \item The direct product $\prod_{i \in I} \mathcal{M}_i$ eliminates quantifiers relative to the sorts $\mathcal{M}_i$. In particular, the sorts $\mathcal{M}_i$ are stably embedded, and pairwise orthogonal.
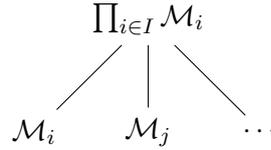
            \begin{center}
                \begin{tikzpicture}
                    \node{$\prod_{i\in I}\mathcal{M}_i$ }
                        child {  node {$\mathcal{M}_i$}}
                        child { node {$\mathcal{M}_j$}} 
                        child { node {$\cdots$}
                        };
                \end{tikzpicture}
            \end{center}
        \end{itemize}
        
    \end{fact}
    
    \begin{proof}
        The first point is easy and is solved by simple inspection on formulas. For the second point, we leave to the reader to prove quantifier elimination. 
        Stable embeddedness is clear by inspection: a formula $\phi(x_i)$ with variable $x_i\in \mathcal{M}_i$ without $M$-sorted quantifiers is a finite Boolean combination of formulas of the form 
        \[\phi_i(x_i, a_i, \pi_i(a)) \cup \bigcup_{j\in I \setminus \{i\} } \phi_j(a_j,\pi_j(a)) \cup \phi(a),\]
        where $\phi_i$ is an $\mathrm{L}_i$-formula, $\phi_j(a_j,\pi_j(a))$ are closed $\mathrm{L}_j$-formula $j\in I\setminus \{i\} $ and $\phi(a)$ is a closed formula in the empty language. It is clearly equivalent to an $\mathrm{L}_i$-formula with parameters in $\mathcal{M}_i$, as a closed formula is true or false and can be replaced either by $x_i=x_i$ or by $x_i\neq x_i$ . Same argument for orthogonality. 
        
    \end{proof}

    Naturally, we have a generalisation of Proposition \ref{sumburden} for infinite products. 
     \begin{proposition} \label{PropInfProd}
        Let $\mathcal{M}_i=(M_i,\ldots)$ be a structure in a language $\mathrm{L}_i$, for $i\in I$ a set of indices. Assume they are not all finite. One has:
            \begin{itemize}
                \item $\bdn(\bigcup_{i\in I}{\mathcal{M}_i})= \sup^\star_{i\in I}\bdn(M_i),$
                \item $\bdn(\prod_{i\in I} \mathcal{M}_i) = \sum^\star_{i\in I} \bdn(M_i).$
            \end{itemize}
    \end{proposition}
    
    \begin{remark}
        If all structures $\mathcal{M}_i$ are finite, there are two cases: either $\#\lbrace i\in I \ \vert \ \vert M_I\vert >1 \rbrace$ is infinite and $\bdn(\prod_{i\in I}\mathcal{M}_i)=1$, or $\#\lbrace i\in I \ \vert \ \vert M_I\vert >1 \rbrace$ is finite and $\bdn(\prod_{i\in I}\mathcal{M}_i)=0$. As the condition $\vert M_i\vert>1$ cannot be seen in terms of burden, this case is told apart.
    \end{remark}
    
    \begin{proof}
    The first point is clear: an inp-pattern $P(x)$ in $\bigcup_i \mathcal{M}_i$ has to `choose' in which sort $\mathcal{M}_i$ its variable $x$ lives. This sort, say $\mathcal{M}_{i_0}$, is stably embedded by Fact \ref{FactQEProd}. By Remark \ref{rkSEDS}, the depth of $P(x)$ is bounded by $\bdn(\mathcal{M}_{i_0})$. Going to the supremum, one sees that definitions match: 
    \[\bdn(\bigcup_{i\in I}{\mathcal{M}_i})= {\sup}^\star_{i\in I}\bdn(M_i).\]

    The second point is more subtle: if $Q(x)$ is an inp-pattern of depth $\mu$ in $\prod_{i \in I} \mathcal{M}_i$, with the variable $x$ in the main sort, then the pattern refers to the sorts $\mathcal{M}_i$ simultaneously. 
    
    \begin{claim}
    Assume that $\prod_{i \in I} \mathcal{M}_i$ admits an inp-pattern of depth $\mu$.
    Then there is an inp-pattern of depth $\mu$ in $\prod_{i \in I} \mathcal{M}_i$ of the following form:
    \[\{ \phi_{\alpha}(\pi_{f(\alpha)}(x),y_{f(\alpha)}), (a_{f(\alpha),j})_{j< \omega}\}_{\alpha<\mu},\]
    for some function $f:\mu \rightarrow I$ and where $\phi_{\alpha}(x_{f(\alpha)},y_{f(\alpha)})$ is a $\mathcal{M}_{f(\alpha)}$-formula. In other word, we may assume that a line $\alpha$ "mentions" only one structure $\mathcal{M}_i$.
    \end{claim}
    \begin{proof}
 Let us assume that $\prod\mathcal{M}_i$ admits an inp-pattern $Q(x)=\lbrace \psi_\alpha(x,\bar{y}_\alpha), (\bar{a}_{\alpha,j})_{j<\omega}, k_{\alpha} \rbrace_{\alpha<\mu}$ of depth $\mu\geq 2$. We assume the array $(\bar{a}_{\alpha,j})_{\alpha<\mu,j<\omega}$ to be mutually indiscernible. To simplify the notation, a generic line of $Q(x)$ is denoted by $\lbrace \psi(x,\bar{y}), (\bar{a}_j)_{j<\omega},k\rbrace$ (we drop the index $\alpha$). By relative quantifier elimination and by Lemma \ref{INPdisj}, we may assume that formulas $\psi(x,\bar{y})$ in $Q(x)$ are of the form 
    \[ \bigwedge_{n<N} x\neq y_n \wedge x=y \wedge \bigwedge \phi_i(\pi_i(x), y_i), \]
    where $\phi_i(x,y_i)$ are $\mathrm{L}_i$-formulas, $N\in \mathbb{N}$ and where $\bar{y}=(y_1,\ldots, y_N, y) \cup (y_i)_{i\in I}$ and $\bar{a}_j=(a_{1j},\ldots, a_{Nj}, a_j) \cup (a_{ij})_{i\in I}$ for $j<\omega$.  If the atomic formula $x=y$ does occur, for example in the first row, then consistency of paths contradicts $k_2$-inconsistency of the second row. Thus, we knows that formulas in $Q(x)$ are of the form 
    \[ \bigwedge_{n<N} x\neq y_n \wedge \bigwedge \phi_i(\pi_i(x), y_i), \]
    Now, the formula $\bigwedge_{n<N} x\neq y_n$ is co-finite. This implies that 
    \[ \lbrace \bigwedge \phi_i(\pi_i(x), a_{i,j})\rbrace_{j<\omega}\]
    is $k+1$-inconsistent. Indeed, otherwise, for one (equivalently for all) $k+1$-increasing tuple $j_0<\cdots<j_k<\omega$, the set
    \[\left\lbrace\bigwedge \phi_i(\pi_i(x), a_{i,j_1}), \ldots, \bigwedge \phi_i(\pi_i(x), a_{i,j_k}) \right\rbrace\]
    is satisfied by $a_{n,j_l}$ for some $n<N$ and $l\leq k$. Without loss of generality, assume that $n=N-1$ and $l=k$. Then, by mutual indiscernibility, $(a_{N-1,j})_{j\geq k}$ are solutions of
    \[\left\lbrace\bigwedge \phi_i(\pi_i(x), a_{i,0}), \ldots, \bigwedge \phi_i(\pi_i(x), a_{i,k-1}) \right\rbrace\]
    This contradicts the $k$-inconsistency of the line
    \[ \left\{\bigwedge_{n<N} x\neq a_{n,j} \wedge \bigwedge \phi_i(\pi_i(x), a_{i,j})\right\}_{j<\omega}, \]
    unless $(a_{N-1,j})_{j<\omega}$ is constant. In that case, this parameter can be ignore:
    replace the formula by 
    \[ \bigwedge_{n<N-1} x\neq y_n \wedge \bigwedge \phi_i(\pi_i(x), y_i), \]
    and we still have an inp-pattern. We get our contradiction by induction on $N$.
    Hence, we may assume that formulas $\psi(x,\bar{y})$ in $Q(x)$ are of the form 
    \[ \bigwedge \phi_i(\pi_i(x), y_i). \]
    We may now conclude using mutual indiscernibility that for at least one $i=:f(\alpha)$, the set
    \[\left\lbrace \phi_i(\pi_i(x), a_{i,j})\right\rbrace_{j<\omega}\]
    is $k$-inconsistant. We may replace the formula $\bigwedge \phi_i(\pi_i(x), y_i)$ by $\phi_{f(\alpha)}(\pi_{f(\alpha)},y_{f(\alpha)})$. In other word, we may assume that only the index $i=f(\alpha)$ occurs in the formula of the line $\alpha$. We found an inp-pattern of the desired form. 
    \end{proof} 
    
    We denote $\bdn(\mathcal{M}_i)$ by $\lambda_i$ and $\sup_{i\in I}\act(\lambda_i)\in \Card$ by $\lambda$. One immediate corollary is that 
    
    \[\bdn(\prod\limits_{i\in I}\mathcal{M}_i)= \bdn(\prod\limits_{\substack{i\in I \\ \lambda_i\neq 0}}\mathcal{M}_i) \geq 1,\]
    
    (notice that we used that some $\mathcal{M}_i$ is infinite). We may assume that $I=\supp(\lambda_i)_{i\in I}$. Now, the proof is straight forward and is just a case study. We distinguish six cases:

    \textbf{First case}: the cardinals $\vert I \vert$ and $\lambda$ are finite. Then, this is immediate from the previous claim: $\bdn(\prod\mathcal{M}_i)=\sum \lambda_i={\sum}^\star \lambda_i$.

    \textbf{Second case}: we have $\vert I \vert \geq \lambda $ and $\vert I \vert \geq \aleph_0$. Then, let $(b_{i,j})_{j<\omega}$ be a sequence of pairwise distinct elements of $\mathcal{M}_i$. Let $x$ be a variable in the main sort. Then, $\lbrace \pi_i(x)=y_i, (b_{i,j})_{j<\omega} \rbrace_{i\in I}$ is an inp-pattern of depth $\vert I \vert$. We have $\bdn(\prod M_i) \geq \vert I \vert$. 
    Reciprocally, assume $\prod\mathcal{M}_i$ admits an inp-pattern $Q(x)$ of depth $\mu> \vert I \vert $. By the previous claim and pigeonhole principle, we find an inp-pattern of depth $\mu$ in some $\mathcal{M}_i$, which is a contradiction with $\lambda \leq \vert I \vert <\mu$.
    We get $\bdn(\prod\mathcal{M}_i)=\sum^\star_{i\in I} \lambda_i = \vert I \vert$.

    \textbf{Third case}: we have $\vert I \vert < \lambda$ and $\lambda_i=\lambda \geq \aleph_0$ for some $i\in I$. Then clearly $\bdn(\prod \mathcal{M}_i) \geq \lambda$. Again, by pigeonhole principle, one gets $\bdn(\prod \mathcal{M}_i) \leq \lambda$.

    \textbf{Fourth case}: we have $\vert I \vert < \lambda$ and $\cf(\lambda) \leq \# \lbrace i\in I \ \vert \ \lambda_i=\lambda_- \rbrace$. Then, choose any sequence of cardinals $(\mu_\alpha)_{\alpha<\cf(\lambda)}$ with  supremum $\lambda$ (in the usual sense) and $\mu_\alpha<\lambda$ for all $\alpha$. We can assume that $I=\cf(\lambda)$ and that we have an inp-pattern $Q_i(x_i)$ in $\mathcal{M}_i$ of depth $\mu_i$. The inp-pattern $Q(x) = \cup_{i\in I} Q_i(\pi_i(x))$ is of depth $\lambda$. We get $\bdn(\prod\mathcal{M}_i)={\sum_{i\in I}}^\star \lambda_i =\lambda.$

    \textbf{Fifth case}: we have $\sup \lbrace \act(\lambda_i) \ \vert \ \lambda_i \notin \lbrace \lambda_-, \lambda \rbrace \rbrace = \lambda$. We conclude as in the previous case that $\bdn(\prod\mathcal{M}_i)=\sum_{i\in I} \lambda_i =\lambda$.

    \textbf{Last case}: we are not in the above cases. Then, by the previous claim, there is no inp-pattern of depth $\lambda$ in $\prod\mathcal{M}_i$. We have then
    \[\bdn(\prod\mathcal{M}_i) = {\sum_{i\in I}}^\star \lambda_i = \lambda_-. \qedhere\]
    \end{proof}
    
    In a supersimple theory, the burden of a complete type is always finite (see \cite{Adl07}). Hence, supersimple theories are examples of strong theories.

    \begin{example}
        The following structures have burden $\aleph_{0-}$:
        \begin{itemize}
            \item Any union structure $M= \bigcup_n M_n$,  where for every $n\in \mathbb{N}$, $M_n$ is a structure of burden $n$.
            \item Any model of $\text{ACFA}$, the model companion of the theory of algebraically closed fields with an automorphism.
            \item Any model of $\text{DCF}_0$, the theory of differentially closed fields.
        \end{itemize}
            
    \end{example}
    The first example is clear by the previous discussion but could look artificial. It will naturally appear when we will discuss the burden of the $\RV_{<\omega}$-sort in mixed characteristic (see Section \ref{SectionUnrimified}). The last one is already given in \cite{CH14}. The fact that the last two examples are of burden $\aleph_{0-}$  follows from the fact they are super-simple and from the next remark, once we notice that such fields are infinite dimensional vector spaces over respectively their fixed field and their constant field:
    
    \begin{remark}\cite[Remark 5.3]{CH14}
        Let $T$ be a simple theory and assume there is a $n$-dimensional type-definable vector space $V$ over a type-definable infinite field $F$. Then there is a type in $V$ of burden $\geq n$.
    \end{remark}
    
    Let us look at one natural example of a non-strong theory:
    
    \begin{remark}
        Let $k$ be an imperfect field of characteristic $p$, considered as a structure in the language of fields. Then $\bdn(k) \geq \aleph_0$.
    \end{remark}
    \begin{proof}   
        Let $e_0,e_1 \in k$ be two linearly independent elements over $k^p$. Then, we have the definable injective map 
        \[\begin{array}{cccc}
        f_2 : & k\times k & \to & k \\
          & (a,b) & \mapsto & a^pe_0+b^pe_1 \\
        \end{array}\]
         By induction, we define for $n\geq 2$:
        \[\begin{array}{cccc}
        f_{n+1} : & k^{(n+1)} & \to & k \\
          & (a_0,\ldots, a_{n-1}) & \mapsto & f_n(a_0,\ldots, a_{n-3},f_2(a_{n-2},a_{n-1})). \\
        \end{array}\]
         
         Then, consider the formula for $n\geq 0$: 
        \[\phi_n(x,y_n) \ \equiv \ \exists y_0,\ldots, y_{n-1}, y_{n+1} \ x= f_{n+2}(y_0,\ldots, y_{n+1}),\]
        and pairwise distinct parameters $b_{n,j} \in k$, for $j<\omega$.
        Then \[\{\phi_n(x,y_n), (b_{n,j})_{j<\omega}\}_{n<\omega}\]
        is an inp-pattern of depth $\aleph_0$.
    \end{proof}
  Similarly, one can show that if a model $\mathcal{M}$ is bi-interpretable on a unary set with the direct product $\mathcal{M}\times \mathcal{M}$, then $\bdn{\mathcal{M}} \geq \aleph_0$ (and it is never of the form $\lambda_-$ where $\lambda$ is a cardinal of cofinality $\aleph_0$).

\subsection{On model theory of algebraic structures}\label{SectionModelTheoryAlgebraicStructures}
\subsubsection{Valued fields}

We gather here some facts on valued fields. After some general statement on indiscernible sequences, we will introduce the $\RV$-sort. Then, we will list the theories of valued fields that we will consider. We will assume the theory of Kaplansky known. We will need a few lemmas, such as a kind of transitivity of pseudo-limits, and a case study of indiscernible sequences. A valued field will be typically denoted by $\mathcal{K}=(K,\Gamma,k,\val)$ where $K$ is the field (main sort), $\Gamma$ the value group and $k$ the residue field. The valuation is denoted by $\val$, the maximal ideal $\mathfrak{m}$ and the valuation ring $\mathcal{O}$. We recall the two traditional languages of valued fields.\\

\textbf{Notation and languages}\\
    We will work in different (many-sorted) languages. Let us define two of them:
    \begin{itemize}
        \item $\mathrm{L}_{\text{div}}=\{K,0,1,+,\cdot,\mid\} $, where $\mid$ is a binary relation symbol, interpreted by the division: 
        \[\text{for }a,b \in K, a \mid b \text{ if and only if } \val(a) \leq \val(b).\]
        \item $\mathrm{L}_{\Gamma,k}= \{K,0,1,+,\cdot\} \cup\{k,0,1,+,\cdot\} \cup \{\Gamma,0,\infty,+,<\} \cup \{\val: K \rightarrow \Gamma, \Res: K^2 \rightarrow k\}$. \nomenclature[L]{$\mathrm{L}_{\Gamma,k}$}{}
    \end{itemize}
    where $\Res: K^2 \rightarrow k$ is the two-place residue map, interpreted as follows:
    \[ \Res(a,b) = 
        \begin{cases}
            \res(a/b) \text{ if } \val(a)\geq \val(b) \neq \infty,\\
            0 \text{ otherwise}.
        \end{cases}\]
    In the next paragraphs, we will also introduce the many-sorted languages $\mathrm{L}_{\RV}$ and $\mathrm{L}_{\RV_{<\omega}}$ which involve the leading term structures $\RV$ and $\RV_{<\omega}$.
        
    By bi-interpretability, a theory of valued fields can be expressed indifferently in either of these languages. Let $\mathcal{K}$ be a valued field. If the context is clear, we will often abusively denote by $K,\Gamma,k,\RV,...$ the sorts in $\mathcal{K}$. In general, the sorts of a valued field $\mathcal{L}$ will be denoted by $L,\Gamma_L,k_L,\RV_L...$ and of a valued field $\mathcal{K}'$ by $K',\Gamma',k',\RV',\dots$ etc.\\

\paragraph{\textbf{Pseudo-Cauchy sequences}}

    We will discuss here some simple facts about mutually indiscernible arrays in a valued field $\mathcal{K}$. We will denote by $\bar{\mathbb{Z}}$ the set of integers with extreme elements $\lbrace -\infty, \infty \rbrace$. We will assume that the reader is familiar with pseudo-Cauchy sequences. We recall however the basic definition:
    
    \begin{definition}
        Let $(I,<)$ be a totally ordered index set without greatest element. A sequence $(a_i)_{i \in I}$ of elements of $K$ is \textit{pseudo-Cauchy} \index{Pseudo-Cauchy sequences} if there is $i\in I$ such that for all indices $i<i_1 < i_2 <i_3$, $\val(a_{i_2}-a_{i_1}) < \val(a_{i_3}-a_{i_2})$. We say that $a\in K$ is a pseudo limit of the \textit{pseudo-Cauchy} sequence $(a_i)_{i\in I}$ and we write $(a_i)_{i\in I}{\Rightarrow} a$ if there is $i\in I$ such that for all indices $i<i_1 < i_2$, we have $\val(a- a_{i_1})= \val(a_{i_2}-a_{i_1})$.
     \end{definition}
    
   The next two lemmas give some useful properties of indiscernible pseudo-Cauchy sequences. 
 
    	\begin{lemma}\label{lemmaind1}
    	
	    \begin{enumerate}
	        \item Assume $(a_i)_{i <\omega}$ is an indiscernible sequence and $a$ is a pseudo limit of $(a_i)_{i<\omega}$. Then for any $i$, $\val(a_i-a)=\val(a_i-a_{i+1})$ depends only on $i$ and not on the chosen limit $a$ (for general pseudo-Cauchy sequence, this holds only for $i$ big enough).
	        \item For three mutually indiscernible sequences $(a_i)_{i<\omega}$, $(b_i)_{i<\omega}$ and $(c_i)_{i<\omega}$, if $(a_i)_{i<\omega}{\Rightarrow}b_0$ and $(b_i)_{i<\omega}{\Rightarrow}c_0$, then we have $(a_i)_{i<\omega}{\Rightarrow}c_0$. 
	        
	\begin{center}
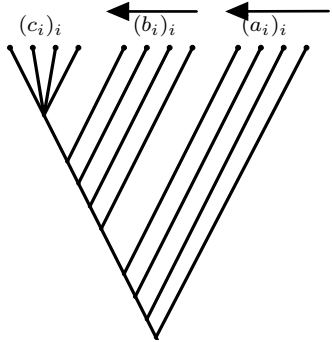
\begin{minipage}{0.70\linewidth}

	\begin{tikzpicture}[line cap=round,line join=round,>=triangle 45,x=1.0cm,y=1.0cm,scale=0.6]
\clip(-2.6813490653752179,-1.444744554454015) rectangle (16.301104990983102,6.5306249304654555);
\draw [line width=1.2pt] (6.,-3.)-- (2.,5.);
\draw [line width=1.2pt] (2.745113232048317,3.509773535903366)-- (3.5,5.);
\draw [line width=1.2pt] (3.2572757384796804,2.485448523040639)-- (4.5,5.);
\draw [line width=1.2pt] (3.5,2.)-- (5.,5.);
\draw [line width=1.2pt] (3.7524098570521747,1.4951802858956507)-- (5.5,5.);
\draw [line width=1.2pt] (4.,1.)-- (6.,5.);
\draw [line width=1.2pt] (7.,5.)-- (4.5,0.);
\draw [line width=1.2pt] (7.5,5.)-- (4.74424472152211,-0.5033601209440769);
\draw [line width=1.2pt] (8.,5.)-- (5.,-1.);
\draw [line width=1.2pt] (8.5,5.)-- (5.2,-1.4);
\draw [line width=1.2pt] (2.5,5.)-- (2.745113232048317,3.509773535903366);
\draw [line width=1.2pt] (3.,5.)-- (2.745113232048317,3.509773535903366);
\draw [->,line width=1pt] (8.983574255409578,5.808506646116047) -- (6.715348226133187,5.801655878143558);
\draw [->,line width=1pt] (6.089660595481235,5.801655878143558) -- (4.1045837982323555,5.801655878143558);
\begin{scriptsize}
\draw [fill=black] (6.,-3.) circle (1.0pt);
\draw [fill=black] (2.,5.) circle (1.5pt);
\draw [fill=black] (2.5,5.) circle (1.5pt);
\draw[color=black] (2.671546559378185,5.47561091638033) node {$(c_{i})_i$};
\draw [fill=black] (3.,5.) circle (1.5pt);
\draw[color=black] (5.1789879436845925,5.47561091638033) node {$(b_i)_i$};
\draw [fill=black] (2.745113232048317,3.509773535903366) circle (1.0pt);
\draw [fill=black] (3.5,5.) circle (1.5pt);
\draw [fill=black] (3.2572757384796804,2.485448523040639) circle (1.0pt);
\draw [fill=black] (4.5,5.) circle (1.5pt);
\draw [fill=black] (3.5,2.) circle (1.0pt);
\draw [fill=black] (5.,5.) circle (1.5pt);
\draw [fill=black] (3.7524098570521747,1.4951802858956507) circle (1.0pt);
\draw [fill=black] (5.5,5.) circle (1.5pt);
\draw [fill=black] (4.,1.) circle (1.0pt);
\draw [fill=black] (6.,5.) circle (1.5pt);
\draw [fill=black] (7.,5.) circle (1.5pt);
\draw [fill=black] (4.5,0.) circle (1.0pt);
\draw [fill=black] (7.5,5.) circle (1.5pt);
\draw [fill=black] (4.74424472152211,-0.5033601209440769) circle (1.0pt);
\draw [fill=black] (8.,5.) circle (1.5pt);
\draw[color=black] (7.57896881089334,5.47561091638033) node {$(a_{i})_i$};
\draw [fill=black] (5.,-1.) circle (1.0pt);
\draw [fill=black] (8.5,5.) circle (1.5pt);
\draw [fill=black] (5.2,-1.4) circle (1.0pt);
\end{scriptsize}
\end{tikzpicture}

\end{minipage}
\end{center}
	        \item If $(a_i)_{i \in \bar{\mathbb{Z}}}$ is an indiscernible sequence in $K$, then $(a_{i})_{i\in \omega}{\Rightarrow}a_{\infty}$ or $(a_{-i})_{i\in \omega}{\Rightarrow}a_{-\infty}$ or for $i\neq j$, $\val(a_{i}-a_{j})$ is constant (in this last case, $(a_i)_{i \in \bar{\mathbb{Z}}}$ will be called \textit{a fan}).
	    \end{enumerate}
	\end{lemma}
	\begin{proof}
	    \begin{enumerate}
	        \item By definition of a pseudo-Cauchy sequence, $(\val(a_i-a_{i+1}))_i$ is eventually strictly increasing. By indiscernibility, it is strictly increasing. Let $i_0$ be such that $\val(a-a_i)=\val(a_{i+1}-a_i)$ for all $i>i_0$. Then $\val(a-a_{i_0}) = \min (\val(a-a_{i_0+1}), \val(a_{i_0+1}-a_{i_0}) ) = \min   (\val(a_{i_0+2}-a_{i_0+1}), \val(a_{i_0+1}-a_{i_0}) ) = \val (a_{i_0+1}-a_{i_0})$. It holds also for $i=i_0$ and we can reiterate.
	        \item Notice that, by mutual indiscernibility and (1), $\val(a_i-b_0)=\val(a_i-a_{i+1})=\val(a_i-b_j)$ for any $i,j < \omega$, \textit{i.e.}  $(a_i)_{i<\omega}{\Rightarrow}b_j$ for any $j$. Similarly, $(b_i)_{i<\omega}{\Rightarrow}c_j$ for any $j$. We have $\val(b_0-b_1) \geq \val(b_0-a_i)= \val(a_i-b_1)$. If $\val(b_0-b_1)=\val(b_0-a_i)$, we have by mutual indiscernibility that $(\val(b_0-a_i))_{i<\omega}$ is constant, which is a contradiction with $(a_i)_{i<\omega}{\Rightarrow}b_0$. Then, we have $\val(b_0-c_0)=\val(b_0-b_1) > \val(b_0-a_i)$.  As $\val(a_i-c_0) \geq \min(\val(a_{i}-b_0), \val(b_0-c_0))$, we deduce that $\val(a_i-c_0)=\val(a_i-b_0)$ for all $i$, \textit{i.e.} $(a_i){\Rightarrow} c_0$.
	        
	        \item It is immediate by indiscernibility (consider for example $\val(a_0-a_1)$ and $\val(a_1-a_2)$).
	    \end{enumerate}
	\end{proof}
	\begin{lemma}\label{lemmaind2}
	Let $(a_{j})_{j \in {\mathbb{Z}}}$ and $(b_{l})_{l \in {\mathbb{Z}}}$ two mutually indiscernible sequences in $K$ such that ${(\val(a_j-b_l))_{j,l}}$ is not constant.   
	At least one of the following occurs:
	\begin{enumerate}
				\item $(a_j)_{j < \omega}{\Rightarrow}b_0$,
				\item $(b_{l})_{l < \omega}{\Rightarrow}a_{0}$,
				\item $(a_{-j})_{j < \omega}{\Rightarrow}b_{0}$,
				\item $(b_{-l})_{l < \omega}{\Rightarrow}a_{0}$.
	\end{enumerate}
		
		Note that if for example $(b_{l})_{l < \omega}{\Rightarrow}a_{0}$, then by mutual indiscernibility, ${(b_{l})_{l < \omega}{\Rightarrow}a_{j}}$ for every $j \in \mathbb{Z}$.
	\end{lemma}

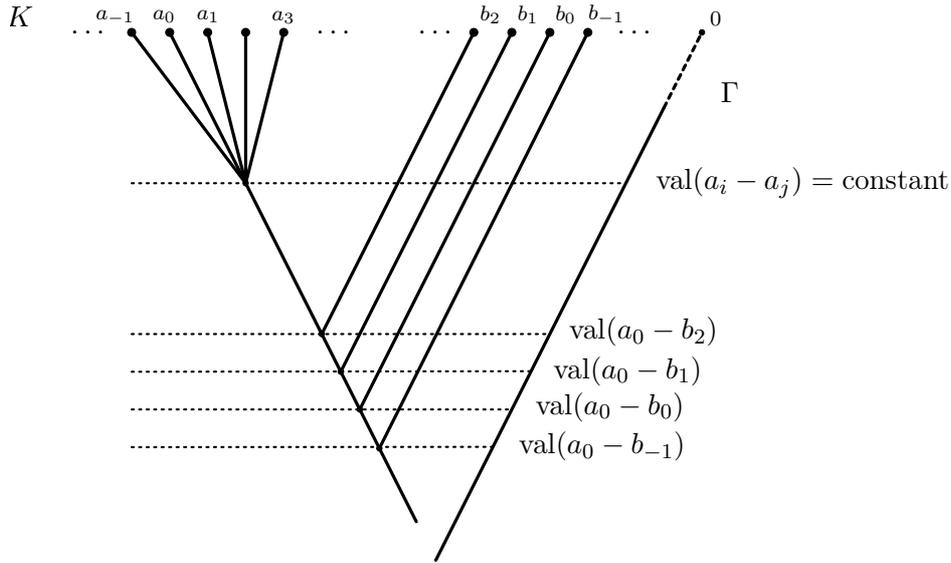
	\begin{proof}
		Since $\val(a_{j}-b_{l})$ is not constant, using the mutual indiscernibility, one of the following occurs:
		\begin{multicols}{2}		
		\begin{enumerate}
		
			\item ${\val(a_{0}-b_{0}) < \val(a_{1}-b_{0})}$,
			\item ${\val(a_{0}-b_{0}) < \val(a_{0}-b_{1})}$,
			\item ${\val(a_{0}-b_{0}) < \val(a_{-1}-b_{0})}$,
			\item ${\val(a_{0}-b_{0}) < \val(a_{0}-b_{-1})}$.
		\end{enumerate}
		\end{multicols}	
			Indeed, if 1. and 3. do not hold, then the sequence ${(\val(b_{0}-a_{j}))_{j \in \mathbb{Z}}}$ is constant. If 2. and 4. do not hold,  then the sequence ${(\val(b_{l}-a_{0}))_{l \in \mathbb{Z}}}$ is constant. This cannot be true for both sequences as it would contradict the assumption. We conclude by indiscernibility. 

\begin{center}

		\begin{tikzpicture}[line cap=round,line join=round,>=triangle 45,x=1.0cm,y=1.0cm]
\clip(-7.087801400314143,3.62355157566648) rectangle (7.02149290705931,10.624429216441536);
\draw[color=black] (2.3691804349710825,9.212012843205708) node {$\Gamma$};
\draw[color=black] (-6.9691804349710825,10.212012843205708) node {$K$};
\draw [line width=1.2pt] (-4.,8.)-- (-4.5,10.);
\draw [line width=1.2pt] (-4.,8.)-- (-4.,10.);
\draw [line width=1.2pt] (-4.,8.)-- (-3.5,10.);
\draw (-3.2064331491892544,10.214496117852873) node[anchor=north west] {$\cdots$};
\draw [line width=1.2pt] (-3.,6.)-- (-1.,10.);
\draw [line width=1.2pt] (-2.75,5.5)-- (-0.5,10.);
\draw [line width=1.2pt] (0.,10.)-- (-2.5,5.);
\draw (-1.866745355775492,10.214496117852873) node[anchor=north west] {$\cdots$};
\draw [line width=1.2pt] (0.5,10.)-- (-2.242086502327691,4.484173004655382);
\draw [line width=1.2pt] (-5.,10.)-- (-1.755846959657945,3.51169391931589);
\draw [line width=1.2pt] (-4.,8.)-- (-5.5,10.);
\draw (-6.434589399163642,10.214496117852873) node[anchor=north west] {$\cdots$};
\draw (0.7527640338747059,10.214496117852873) node[anchor=north west] {$\cdots$};
\draw (1.2614108941988953,8.34239195120329) node[anchor=north west] {$\val(a_i-a_j) = \text{constant} $};
\draw (-0.5394739355975589,4.836852778698755) node[anchor=north west] {$\val(a_0-b_{-1})$};
\draw (0.12039226157976783,6.362793359671317) node[anchor=north west] {$\val(a_0-b_2)$};
\draw [line width=1.2pt] (-1.5,3.)-- (1.5,9.);
\draw [line width=1.2pt,dash pattern=on 2pt off 2pt] (1.5,9.)-- (2.,10.);
\draw [line width=0.8pt,dotted] (-5.5,8.)-- (1.,8.);
\draw [line width=0.8pt,dotted] (-5.5,6.)-- (0.,6.);
\draw [line width=0.8pt,dotted] (-5.5,5.5)-- (-0.24655329867187703,5.506893402656246);
\draw [line width=0.8pt,dotted] (-5.5,5.)-- (-0.5,5.);
\draw [line width=0.8pt,dotted] (-5.5,4.5)-- (-0.748053812275157,4.503892375449686);
\draw (-0.08581592503814676,5.840399286905935) node[anchor=north west] {$\val(a_0-b_1)$};
\draw (-0.31951853653844997,5.359246851464136) node[anchor=north west] {$\val(a_0-b_0)$};
\begin{scriptsize}
\draw [fill=black] (-5.,10.) circle (1.5pt);
\draw[color=black] (-5.0691804349710825,10.212012843205708) node {$a_{0}$};
\draw [fill=black] (-4.,8.) circle (1.0pt);
\draw [fill=black] (-4.5,10.) circle (1.5pt);
\draw[color=black] (-4.491797512440922,10.212012843205708) node {$a_{1}$};
\draw [fill=black] (-4.,10.) circle (1.5pt);
\draw [fill=black] (-3.5,10.) circle (1.5pt);
\draw[color=black] (-3.501998216674932,10.198265630764514) node {$a_{3}$};
\draw [fill=black] (-3.,6.) circle (1.0pt);
\draw [fill=black] (-1.,10.) circle (1.5pt);
\draw[color=black] (-0.7800501533184593,10.239507268088095) node {$b_2$};
\draw [fill=black] (-0.5,10.) circle (1.5pt);
\draw[color=black] (-0.28515050543546416,10.239507268088095) node {$b_1$};
\draw [fill=black] (-2.75,5.5) circle (1.0pt);
\draw [fill=black] (0.,10.) circle (1.5pt);
\draw[color=black] (0.20974914244753085,10.239507268088095) node {$b_0$};
\draw [fill=black] (-2.5,5.) circle (1.0pt);
\draw [fill=black] (0.5,10.) circle (1.5pt);
\draw[color=black] (0.7458904276541087,10.239507268088095) node {$b_{-1}$};
\draw [fill=black] (-2.242086502327691,4.484173004655382) circle (1.0pt);
\draw [fill=black] (-5.5,10.) circle (1.5pt);
\draw[color=black] (-5.7152994197072156,10.198265630764514) node {$a_{-1}$};
\draw [fill=black] (2.,10.) circle (1.0pt);
\draw[color=black] (2.182474127758914,10.177644812102722) node {$0$};
\end{scriptsize}
\end{tikzpicture}

\end{center}
	
	\end{proof}

\paragraph{\textbf{The leading term structure}}\label{SusubsectionRVSort}
We will now define the \emph{$\RV$-sort} (or \emph{$\RV$-sorts}, as we may need to consider more than one sort) --an intermediate structure between the valued field and its value group and residue field. We also introduce corresponding languages $\mathrm{L}_{\RV}$ and $\mathrm{L}_{\RV_{<\omega}}$. This paragraph is largely inspired by \cite{Fle11}, which one can use as a reference. Let $K$ be a Henselian valued field of characteristic $(0,p)$ with $p \geq 0$, of value group $\Gamma$ and residue field $k$. If $\delta \in \Gamma_{\geq 0}$, we denote by $\mathfrak{m}_\delta$ the ideal of the valuation ring $\mathcal{O}$ defined by $\lbrace x\in \mathcal{O} \ \vert \ v(x)>\delta \rbrace$.
The \textit{leading term structure} of order $\delta$ \index{Leading term structure $\RV_\delta$ of order $\delta$} \nomenclature[]{$\RV_\delta$}{Leading term structure of order $\delta$} is the quotient group 
$$ \RV_\delta^\star:= K^\star/(1+\mathfrak{m}_\delta).$$
The quotient map is denoted by $\rv_\delta: K^\star \rightarrow \RV_\delta^\star$. The valuation $\val:K^\star \rightarrow \Gamma$ induces a group homomorphism $\val_{\rv_\delta}: \RV_\delta^\star \rightarrow \Gamma$. Since $\mathfrak{m}=\mathfrak{m}_0$ and ${k^\star:= (\mathcal{O}/\mathfrak{m}_0)^\star \simeq \mathcal{O}^\times/{(1+\mathfrak{m}_0)}}$, we have the following short exact sequence: 
$$1 \rightarrow k^\star \overset{\iota}{\rightarrow} \RV_{0}^\star \overset{\val_{\rv_0}}{\rightarrow} \Gamma \rightarrow 0.$$
In general, we denote by $\mathcal{O}_\delta$ \nomenclature[]{$\mathcal{O}_\delta$}{Residue ring of order $\delta$} the ring $\mathcal{O}/m_\delta$, called the \textit{residue ring of order} $\delta$ \index{Residue ring $\mathcal{O}_\delta$ of order $\delta$}. One has $\mathcal{O}_\delta^\times \simeq \mathcal{O}^\times/{(1+\mathfrak{m}_\delta)}$ and the following exact sequence:
\[1 \rightarrow \mathcal{O}_\delta^\times \overset{\iota_{\delta}}{\rightarrow} \RV_{\delta}^\star \overset{\val_{\rv_\delta}}{\rightarrow} \Gamma \rightarrow 0.\]

Furthermore, as $\mathfrak{m}_\gamma \subseteq \mathfrak{m}_\delta$ for any $\delta \leq \gamma$ in $\Gamma_{\geq 0}$, we have a projection map
$ \RV_\gamma^\star \rightarrow \RV_\delta^\star$ denoted by $ \rv_{\gamma \rightarrow \delta}$ or simply by $\rv_\delta$.
We add a new constant $\mathbf{0}$ to the sort $\RV_\delta^\star$ and we write $\RV_\delta:=\RV_\delta^\star\cup \lbrace \mathbf{0} \rbrace$. We set the following properties:
	\begin{itemize}
 		\item for all $\mathbf{x}\in \RV_\delta$, $\mathbf{0}\cdot \mathbf{x} = \mathbf{x} \cdot \mathbf{0} = \mathbf{0}$.
 		\item $val_{\rv_\delta}(\mathbf{0})= \infty, \quad \rv_\delta (0)= \mathbf{0}$.
	\end{itemize}

\begin{proposition}
For any $a,b \in K$ and $\delta \in \Gamma_{\geq 0}$, $\rv_\delta(a)=\rv_\delta(b)$ if and only if $\val(a-b)> \val(b)+\delta$ or $a=b=0$.
\end{proposition}
\begin{proof}
	This follows easily from the definition: assume $\rv_\delta(a)=\rv_\delta(b)$ and $a\neq 0$. Then $a=b(1+\mu)$ for some $\mu \in m_\delta$ and $\val(a-b)=\val(b)+\val(\mu) > \val(b)+\delta$. Conversely, if $\val(a-b)> \val(b)+\delta$, one can write $a=b(1+\frac{(a-b)}{b})$.
\end{proof}
As a group quotient, the sort $\RV_{\delta}$ is endowed with a multiplication. As we will see, it also inherits from the field some kind of addition. 
\begin{notation}
Let $0 \leq \delta_1,\delta_2,\delta_3$ be three elements of $\Gamma$ and $\mathbf{x} \in \RV_{\delta_1}$, $\mathbf{y} \in \RV_{\delta_2}$, $\mathbf{z} \in \RV_{\delta_3}$ three variables.
Then we define the following formulas:
\begin{align*}
 \oplus_{\delta_1,\delta_2,\delta_3} (\mathbf{x},\mathbf{y},\mathbf{z}) & \equiv & \exists a,b \in K \ \nomenclature[]{$\oplus_{\delta_1,\delta_2,\delta_3}(\mathbf{x},\mathbf{y},\mathbf{z})$}{} \rv_{\delta_1}(a)=\mathbf{x} \wedge \rv_{\delta_2}(b)=\mathbf{y} \wedge \rv_{\delta_3}(a+b)=\mathbf{z}
 \end{align*}
\end{notation}

In our study of valued fields, we will consider the structures $\RV$ and $\RV_{<\omega}$, that we define now:

\begin{definition}
    The $\RV$-sort of a valued field $\mathcal{K}$ is the first order leading term structure \nomenclature[]{$\RV$}{$\RV$-sort or first order leading term structure}\index{$\RV$-sort}
    \[\RV=(\RV_0, \cdot, \oplus_{0,0,0}, \mathbf{0}, \mathbf{1}) \]
    endowed with its natural structure of abelian group and the ternary predicate described above. Following the usual convention, we drop the index $0$, and write $\RV, \oplus$ and $\rv$  instead of $\RV_0, \oplus_{0,0,0}$ and $\rv_0$.
\end{definition}
    
\begin{fact} [{ Flenner, \cite[Proposition 2.8]{Fle11}}]
The three-sorted structure $\{(\RV,\mathbf{1},\cdot,\mathbf{0},\mathbf{1}), (k,0,1,+,\cdot), (\Gamma,0,+,<), \iota,\val_{\rv} \}$ and the one-sorted structure $\{\RV, \mathbf{0}, \cdot, \oplus\}$ are bi-interpretable on unary sets.
\end{fact} 
This will mean, in the context of this paper, that these two points of view are equivalent, and we will swap between one to the other indifferently (see Fact \ref{FactBdnInterpretUnarySet}).

We defined the \textit{leading term language} $\mathrm{L}_{\RV}$\nomenclature[L]{$\mathrm{L}_{\RV}$}{} as the multisorted language with 
	
	 \begin{itemize}
	    \item a sort for $K$ and $\RV$.
	    \item the ring language for $K$,
	    \item the (multiplicative) group language as well as the symbol $0$ for $\RV$.
	    \item the ternary relation symbol $\oplus$ and the function symbols $\rv$.
	 \end{itemize}
The structure $\mathcal{K}=\left(K,\RV,\rv\right)$ becomes a structure in this language where all symbols are interpreted as before. This language is also bi-interpretable (without parameters) with the usual languages of valued fields, e.g. with $\mathrm{L}_{\text{div}}$  (see \cite[Proposition 2.8]{Fle11}).

Also, notice that the symbol $\oplus$ \nomenclature[]{$a\oplus b$}{} suggests a binary operation. Occasionally, we will indeed write  $\rv(a)\oplus \rv(b)$ for $a,b\in K$ to denote the following element:
	\[ \rv(a) \oplus \rv(b) := \begin{cases} \rv(a+b) &\ \text{if } \val(a+b)= \min(\val(a),\val(b)), \\
	\mathbf{0} &\text{ otherwise}.
	\end{cases}\]
	It is not hard to see that this is independent of the choice of representatives of $\rv(a)$ and $\rv(b)$. We will write $\bigoplus_{i\in I} \mathbf{a}_i $
    for $I$ a set of indices and $\mathbf{a}_i \in \RV$, when such a sum does not depend on any choices of parentheses.  
Notice that the law $\oplus$ is not an addition. If $a,b\in K$ with $\val(a)<\val(b)$, we have that $\rv(a) \oplus \rv(b) = \rv(a) = \rv(a+b)$. Also, notice that it is  in general not true that $\rv(a+b)= \rv(a) \oplus \rv(b)$ (choose $a,b\in K$ such that $\rv(a) = - \rv(b)$ and $a\neq -b$). When we have indeed that  $\rv(a+b)= \rv(a) \oplus \rv(b)$, we say that the sum $\rv(a) \oplus \rv(b)$ is well-defined.

In the specific context of mixed characteristic Henselian valued fields, we might have to consider a larger structure:
\begin{definition}
Assume  that $\mathcal{K}$ is a valued field of characteristic $0$ and residue characteristic $p\geq 0$. We reserve now the notation $\delta_n$ for $\delta_n= \val(p^n)$. 
We write $\RV_{<\omega}$ \index{$\RV_{<\omega}$-sort} \nomenclature[]{$\RV_{<\omega}$}{Leading term structure of finite order} for the union of sorts leading term structure of finite order $\left\{(\RV_{\delta_n})_{n<\omega},(\oplus_{\delta_l,\delta_m, \delta_n})_{n<l,m}, (\rv_{\delta_n\rightarrow \delta_m})_{m<n<\omega}  \right\}$ endowed with ternary predicates $\oplus_{\delta_l,\delta_m, \delta_n}$ and a projective system of maps $(\rv_{\delta_n\rightarrow \delta_m})_{m<n<\omega}$. We also write $\val_{\rv_{<\omega}}: \RV_{<\omega} \rightarrow \Gamma\cup \lbrace \infty \rbrace $ for $\bigcup_{n<\omega}(\val_{\rv_{\delta_n}}: \RV_{\delta_n}\rightarrow \Gamma\cup \lbrace \infty \rbrace)$, etc.
\end{definition}

\begin{remark}\label{RemarkIdentifyRVRVOmega}
In equicharacteristic $0$, we have that $\delta_n:=\val(p^n)=0$ for all $n<\omega$.  This leads to identifying $\bigcup_{n<\omega}\RV_{\delta_n}$ with $\RV=\RV_0$.
\end{remark}

In Section \ref{SectionUnrimified}, we will have to use another language to describe the induced structure on $\RV_{<\omega}$:
    \begin{fact} [{ Flenner, \cite[Proposition 2.8]{Fle11}}] \label{FactRVBiInterpretability}
        The structure
        \[\left\{(\RV_{\delta_n})_{n<\omega},(\oplus_{\delta_l,\delta_m, \delta_n})_{n<l,m}, (\rv_{\delta_n\rightarrow \delta_m})_{m<n<\omega}  \right\}\]
        and the structure 
        \begin{align*}
            \left\{(\RV_{\delta_n})_{n<\omega},(\mathcal{O}_{\delta_n},\cdot,+,0,1)_{n<\omega}, (\Gamma,+,0,<), (\val_{\rv_{\delta_n}})_{n<\omega},\right.\\ \left. (\mathcal{O}_{\delta_n}^\times\rightarrow \RV_{\delta_n}^\times)_{n<\omega},(\rv_{\delta_n\rightarrow \delta_m})_{m<n<\omega}\right\} 
        \end{align*}
        are bi-interpretable on unary sets. 
    \end{fact}
    As before, this is only to say that one can recover the valuation using the symbols $\oplus$ (see  \cite[Proposition 2.8]{Fle11}). 
     And again, this fact means that we will be able to swap between one language to the other indifferently (see Fact \ref{FactBdnInterpretUnarySet}).

We defined the language $\mathrm{L}_{\RV_{<\omega}}$ \nomenclature[L]{$\mathrm{L}_{\RV_{<\omega}}$}{} as the multisorted language with 
	 
	 \begin{itemize}
	    \item sorts for $K$ and $\RV_{\delta_n}$ for $n<\omega$.
	    \item the ring language for $K$,
	    \item for all $n<\omega$, the (multiplicative) group language as well as the symbol $0$ for $\RV_{\delta_n}$.
	    \item relation symbols $\oplus_{\delta_l,\delta_m,\delta_n}$ for $n \leq l, m$ integers, function symbols ${\rv_{\delta_n}: \mathcal{K} \rightarrow \RV_{\delta_n}}$ and \\ ${\rv_{\delta_n\rightarrow \delta_m}: \RV_{\delta_n} \rightarrow \RV_{\delta_m}}$ for $n>m$.
	 \end{itemize}
The structure $\mathcal{K}=\left(K,(\RV_{\delta_n})_{n<\omega},(\oplus_{\delta_l,\delta_m, \delta_n})_{n<l,m},(\rv_{\delta_n})_{n<\omega}, (\rv_{\delta_n\rightarrow \delta_m})_{m<n<\omega} \right)$ becomes a structure in this language where all symbols are interpreted as before. This language is also bi-interpretable (without parameters) with the usual languages of valued fields, e.g. with $\mathrm{L}_{\text{div}}$  (see \cite[Proposition 2.8]{Fle11}).

Let $\mathcal{K}$ be any valued field. Let us state few lemmas.

\begin{notation}
Let $\delta_1,\delta_2, \delta_3 \in \Gamma$ be three values. We write:

\begin{align*}
\WD_{\delta_1,\delta_2,\delta_3}(\mathbf{x},\mathbf{y}) \equiv  \exists ! \mathbf{z}\in \RV_{\delta_3} \ \oplus_{\delta_1,\delta_2,\delta_3} (\mathbf{x},\mathbf{y},\mathbf{z}).\end{align*}\\
\end{notation} \nomenclature[]{$\WD_{\delta_1,\delta_2,\delta_3}(\mathbf{x},\mathbf{y})$, $\WD_{\delta_3}(\mathbf{x},\mathbf{y})$}{}

If the context is clear and in order to simplify notations, we will write:
\begin{itemize}

\item  $\WD_{\delta_3}$ instead of $\WD_{\delta_1,\delta_2,\delta_3}$, 
\item for any formula $\phi(\mathbf{z})$ with $\mathbf{z}\in \RV_{\delta_3}$, $\mathbf{x}\in \RV_{\delta_1}$ and $\mathbf{y} \in \RV_{\delta_2}$:
 \[{\phi(\rv_{\delta_3}(\mathbf{x})+\rv_{\delta_3}(\mathbf{y}))}\] or \[{\phi(\rv_{\delta_3}(\mathbf{x})+\rv_{\delta_3}(\mathbf{y})) \wedge \WD_{\delta_3}(\mathbf{x},\mathbf{y})}\] instead of
\[\exists \mathbf{z} \in \RV_{\delta_3} \ \oplus_{\delta_1,\delta_2,\delta_3} (\mathbf{x},\mathbf{y},\mathbf{z}) \wedge \phi(\mathbf{z}) \wedge \WD_{\delta_1,\delta_2,\delta_3}(\mathbf{x},\mathbf{y}).\]
\end{itemize}

\begin{example}
Take $K= \mathbb{R}((t))$ the field of power series over the reals endowed with the $t$-adic valuation. Consider $x=t^2+t^3+t^4+t^5, \ x'=t^2+t^3+t^4+2t^5 \in K$, $y=-t^2-t^3+t^4-t^5 \in K$ and $z=2t^4,\ z'= 2t^4+t^5 \in K$.

Then, we have $\rv_2(x)=\rv_2(x')$ since $\val(x-x')=5 > \val(x)+2$ but $\rv_1(z) \neq \rv_1(z')$ since $\val(z-z')=5 \ngtr \val(z)+1 $.
We have:
\[ \models \oplus_{2,2,1} (\rv_2(x),\rv_2(y),\rv_1(z)), \quad \models  \oplus_{2,2,1} (\rv_2(x),\rv_2(y),\rv_1(z')),\]
Hence, the sum is not well-defined in $\RV_1$:
\[ \models \neg \WD_{2,2,1}(\rv_2(x),\rv_2(y)) \]
We need to pass to $\RV_3$ in order to get a well-defined sum in $\RV_1$: 
\[ \models  \oplus_{3,3,1} (\rv_3(x),\rv_3(y),\rv_1(z)), \quad \neg\models  \oplus_{3,3,1} (\rv_3(x),\rv_3(y),\rv_1(z')),\] 
\[  \quad \models \WD_{3,3,1}(\rv_3(x),\rv_3(y)). \]
\end{example}

More generally, we have the following proposition:

\begin{proposition}\label{WDdelta}
Let $0 \leq \gamma \leq \delta$ be two elements of $\Gamma_{ \geq 0}$ and $\epsilon = \delta-\gamma \geq 0$. Then for every $a,b \in K^\star$:
$$\WD_{\gamma}(\rv_{\delta}(a),\rv_{\delta}(b)) \quad \text{if and only if} \quad \val(a+b) \leq \min \lbrace \val(a),\val(b) \rbrace + \epsilon.$$
\end{proposition}

\begin{center}
\begin{tikzpicture}[line cap=round,line join=round,>=triangle 45,x=1cm,y=1cm,scale=0.75]
\clip(-3.2383443958462172,1.4009780613773635) rectangle (8.226479699154212,6.683890024741572);
\draw [line width=2pt] (0,0)-- (4,6);
\draw [line width=2pt] (2,6)-- (3,4.5);
\draw [line width=2pt] (-1,6)-- (0,4.5);
\draw [line width=2pt] (0,4.5)-- (1,6);
\draw [line width=2pt] (0,4.5)-- (1.5200297019506115,2.2800445529259172);
\draw [line width=1.2pt,dotted] (-1.5,4.5)-- (4,4.5);
\draw [line width=1.2pt,dotted] (-1.5001915218714257,2.2472413709286223)-- (4.000148141043704,2.243946937151217);
\draw (3.9744122355325544,4.90158795503768) node[anchor=north west] {$\val(a+b)$};
\draw (4.113605784559162,2.7744911706119573) node[anchor=north west] {$\val(a)=\val(b)$};
\draw [line width=1.2pt,dotted] (-1.499800534220066,5.506848846585093)-- (4.000199465779933,5.506848846585093);
\draw (3.999720153537392,5.992057167202163) node[anchor=north west] {$\val(a)+\delta$};
\draw (4.214837456578514,3.8500776858175616) node[anchor=north west] {$\updownarrow \epsilon$};
\draw (4.25279933358577,5.393860684112664) node[anchor=north west] {$\updownarrow \gamma$};
\begin{scriptsize}
\draw [fill=black] (4,6) circle (2.5pt);
\draw[color=black] (4.202183497576095,6.273310834780777) node {$0$};
\draw [fill=black] (2,6) circle (2.5pt);
\draw[color=black] (2.30408964721326,6.273310834780777) node {$a+b$};
\draw [fill=black] (3,4.5) circle (2.5pt);
\draw [fill=black] (-1,6) circle (2.5pt);
\draw[color=black] (-0.8514382263842078,6.273310834780777) node {$a$};
\draw [fill=black] (0,4.5) circle (2.5pt);
\draw [fill=black] (1,6) circle (2.5pt);
\draw[color=black] (0.9038110500124913,6.273310834780777) node {$-b$};
\draw [fill=black] (1.5200297019506115,2.2800445529259172) circle (2.5pt);
\end{scriptsize}
\end{tikzpicture}
\end{center}

\begin{proof}
Assume $\val(a+b) \leq \min \lbrace \val(a),\val(b) \rbrace + \epsilon$.
Let $a' \in K$ such that $\rv_{\delta}(a')= \rv_{\delta}(a)$. This is equivalent to $\val(a-a')> \val(a)+ \delta$, thus we have:  
$$\val\left(a+b - (a'+b)\right) = \val(a-a') > \val(a)+\delta = \val(a)+\epsilon + \gamma \geq \val(a+b) +\gamma. $$
Hence, $\rv_{\gamma}(a'+b)= \rv_{\gamma}(a+b)$. We have proved the implication from right to left. \\

Conversely, assume that $\val(a+b) > \min \lbrace \val(a),\val(b) \rbrace + \epsilon$ and $\min \lbrace \val(a),\val(b) \rbrace = \val(a)$. Let $\eta= \val(a+b) + \gamma$ and take any $c \in K$ of valuation $\eta$. Then $\rv_{\delta}(a)=\rv_{\delta}(a+c)$ since ${\val(a+c-a) = \eta > \val(a) +\delta} $ and ${\rv_{\gamma}(a+b)\neq \rv_{\gamma}(a+c+b)}$ since $\val\left(a+c+b-(a+b)\right) = \eta = \val(a+b)+\gamma$.

\end{proof}

\begin{remark}\label{remark}
To prove $\val(a+b) \leq \min \lbrace \val(a),\val(b) \rbrace + \epsilon$ with $\epsilon \geq 0$, it is enough to show that $\val(a+b) \leq \val(a) + \epsilon$ (or $\val(a+b) \leq \val(b) + \epsilon$). Indeed, if $\val(a)=\val(b)$ then this is clear. If $\val(a)<\val(b)$ or $\val(b)<\val(a)$, this is also clear since we have $\val(a+b)=\val(a) \leq \val(a)+ \epsilon$ in the first case and $\val(a+b)=\val(b)<\val(a)+\epsilon$ in the second.
\end{remark}

 The following lemma is immediate:
\begin{lemma} \label{WDcase}
Let $a,b$ and $c=a-b$ be elements of $K$ and let $\gamma \in \Gamma_{\geq 0}$. At least one of the following holds: 

    \begin{eqnarray}
        \models \WD_\gamma(\rv_\gamma(a),\rv_\gamma(b-a)) \\
	    \models \WD_\gamma(\rv_\gamma(b),\rv_\gamma(a-b))
	\end{eqnarray}
\end{lemma}

\begin{proof}

Notice that exactly one of the following occurs:
	\begin{multicols}{2}
	\begin{enumerate}
		\item $\val{a} = \val{c} < \val{b},$
		\item $\val{b} = \val{c} < \val{a},$
		\item $\val{a} = \val{b} < \val{c},$
		\item $\val{a} = \val{b} = \val{c}.$
	\end{enumerate}	
	\end{multicols}

\begin{center}
\begin{minipage}{0.80\linewidth}
	
\begin{tikzpicture}[line cap=round,line join=round,>=triangle 45,x=1.0cm,y=1.0cm,scale=0.55]
\clip(-7.0610774325468855,1.617819523136445) rectangle (16.682404837938282,13.468622843325615);
\draw [line width=1.2pt] (2.,7.)-- (2.504587309046625,6.00917461809325);
\draw [line width=1.2pt] (1.4971301491130982,3.9942602982261963)-- (0.,7.);
\draw [line width=1.2pt] (0.48589692377660665,6.024480486462519)-- (1.,7.);
\draw [line width=1.2pt] (1.,3.)-- (2.504587309046625,6.00917461809325);
\draw [line width=1.2pt,dash pattern=on 5pt off 5pt] (2.504587309046625,6.00917461809325)-- (3.,7.);
\draw [line width=1.2pt] (5.,7.)-- (5.5,4.);
\draw [line width=1.2pt] (6.,7.)-- (5.5,4.);
\draw [line width=1.2pt] (4.,7.)-- (5.5,4.);
\draw [line width=1.2pt] (5.,3.)-- (6.5,6.);
\draw [line width=1.2pt,dash pattern=on 5pt off 5pt] (6.5,6.)-- (7.,7.);
\draw [line width=1.2pt] (1.993551188801738,11.989497869058052)-- (2.4951669348492516,10.992729361153081);
\draw [line width=1.2pt] (1.4877097749157246,8.977815041286027)-- (0.,12.);
\draw [line width=1.2pt] (0.49564485869469677,10.993129935644955)-- (0.9935511888017377,11.989497869058052);
\draw [line width=1.2pt] (0.9935511888017377,7.989497869058052)-- (2.4951669348492516,10.992729361153081);
\draw [line width=1.2pt,dash pattern=on 5pt off 5pt] (2.4951669348492516,10.992729361153081)-- (2.993551188801737,11.989497869058052);
\draw [line width=1.2pt] (6.,12.)-- (6.5045873090466255,11.009174618093251);
\draw [line width=1.2pt] (5.497130149113098,8.994260298226196)-- (4.,12.);
\draw [line width=1.2pt] (4.485896923776606,11.02448048646252)-- (5.,12.);
\draw [line width=1.2pt] (5.,8.)-- (6.5045873090466255,11.009174618093251);
\draw [line width=1.2pt,dash pattern=on 5pt off 5pt] (6.5045873090466255,11.009174618093251)-- (7.,12.);
\draw (2.3609393414551656,10.579204365965005) node[anchor=north west] {$\Gamma$};
\draw (6.48568890696273,10.642017811125019) node[anchor=north west] {$\Gamma$};
\draw (2.3190637113484898,5.679755643483969) node[anchor=north west] {$\Gamma$};
\draw (6.339124201589365,5.679755643483969) node[anchor=north west] {$\Gamma$};
\draw (-0.6375413208521254,11.982037974538635) node[anchor=north west] {1.};
\draw (3.338146059708777,12.02391360464531) node[anchor=north west] {3.};
\draw (-0.7237900606387739,7.019775806897586) node[anchor=north west] {2.};
\draw (3.359083874762115,6.956962361737573) node[anchor=north west] {4.};
\begin{scriptsize}
\draw [fill=black] (1.,3.) circle (1.0pt);
\draw [fill=black] (3.,7.) circle (1.5pt);
\draw[color=black] (3.135638498428668,7.365249755277659) node {$0$};
\draw [fill=black] (2.,7.) circle (1.5pt);
\draw[color=black] (2.0259343006017594,7.3024363101176455) node {$a$};
\draw [fill=black] (2.504587309046625,6.00917461809325) circle (1.0pt);
\draw [fill=black] (1.4971301491130982,3.9942602982261963) circle (1.0pt);
\draw [fill=black] (0.,7.) circle (1.5pt);
\draw[color=black] (0.015904055481321816,7.239622864957632) node {$a-b$};
\draw [fill=black] (0.48589692377660665,6.024480486462519) circle (1.0pt);
\draw [fill=black] (1.,7.) circle (1.5pt);
\draw[color=black] (1.0837326232015543,7.281498495064308) node {$b$};
\draw [fill=black] (5.,3.) circle (1.0pt);
\draw [fill=black] (7.,7.) circle (1.5pt);
\draw[color=black] (7.28132587898957,7.3024363101176455) node {$0$};
\draw [fill=black] (5.,7.) circle (1.5pt);
\draw[color=black] (5.124730928495767,7.344311940224321) node {$b$};
\draw [fill=black] (5.5,4.) circle (1.0pt);
\draw [fill=black] (6.,7.) circle (1.5pt);
\draw[color=black] (6.234435126322675,7.365249755277659) node {$a-b$};
\draw [fill=black] (4.,7.) circle (1.5pt);
\draw[color=black] (4.09877799088221,7.281498495064308) node {$a$};
\draw [fill=black] (0.9935511888017377,7.989497869058052) circle (1.0pt);
\draw [fill=black] (2.993551188801737,11.989497869058052) circle (1.5pt);
\draw[color=black] (3.1775141285353437,12.285636292812033) node {$0$};
\draw [fill=black] (1.993551188801738,11.989497869058052) circle (1.5pt);
\draw[color=black] (2.151561190921787,12.285636292812033) node {$b$};
\draw [fill=black] (2.4951669348492516,10.992729361153081) circle (1.0pt);
\draw [fill=black] (1.4877097749157246,8.977815041286027) circle (1.0pt);
\draw [fill=black] (0.,12.) circle (1.5pt);
\draw[color=black] (0.050777968558799759,12.250885032598682) node {$a$};
\draw [fill=black] (0.49564485869469677,10.993129935644955) circle (1.0pt);
\draw [fill=black] (0.9935511888017377,11.989497869058052) circle (1.5pt);
\draw[color=black] (1.1256082533082301,12.306574107865371) node {$a-b$};
\draw [fill=black] (5.,8.) circle (1.0pt);
\draw [fill=black] (7.,12.) circle (1.5pt);
\draw[color=black] (7.28132587898957,12.306574107865371) node {$0$};
\draw [fill=black] (6.,12.) circle (1.5pt);
\draw[color=black] (6.08787042094931,12.348449737972047) node {$a-b$};
\draw [fill=black] (6.5045873090466255,11.009174618093251) circle (1.0pt);
\draw [fill=black] (5.497130149113098,8.994260298226196) circle (1.0pt);
\draw [fill=black] (4.,12.) circle (1.5pt);
\draw[color=black] (4.182529251095563,12.285636292812033) node {$a$};
\draw [fill=black] (4.485896923776606,11.02448048646252) circle (1.0pt);
\draw [fill=black] (5.,12.) circle (1.5pt);
\draw[color=black] (5.1456687435491055,12.369387553025383) node {$b$};
\end{scriptsize}
\end{tikzpicture}

\end{minipage}
\end{center}

    Let $\mathbf{a}=\rv_\gamma(a),\mathbf{b}=\rv_\gamma(b)$  and $\mathbf{c}=\rv_\gamma(c)$. In cases 2,3 and 4, the difference between $\mathbf{a}$ and $\mathbf{c}$ is well-defined. 
	$$\models \WD_\gamma(\mathbf{a},-\mathbf{c}).$$
	
	In cases 1,3 and 4, the sum of $\mathbf{b}$ and $\mathbf{c}$ is well-defined:
	$$\models \WD_\gamma(\mathbf{b},\mathbf{c}).$$
\end{proof}

\paragraph{\textbf{Benign theory of Henselian valued fields}}\label{Preliminaries Benign theory of Henselian valued fields}
    Later in this text, we prove transfer principles for some rather nice Henselian valued fields, that we called here `benign' (see Definition \ref{DefBenign} below). The goal of this paragraph is to discuss some essential properties that these benign Henselian valued fields share and that we will use for proving Theorem \ref{ThmBdnHenValFieCha00}. 

    Let $T$ be a (possibly incomplete) theory of Henselian valued fields. We need first to recall the definition of \emph{an angular component} \index{Angular component} (or $\ac$\emph{-map}) \nomenclature[]{$\ac$}{Angular component}. It is a group homomorphism usually denoted by $\ac: (K^\star,\cdot) \rightarrow (k^\star,\cdot)$ such that $\restriction{ac}{\mathcal{O}^\times}=\restriction{\res}{\mathcal{O}^\times}$. We also set $\ac(0)=0$. We have the following diagram:

$$\xymatrix{
1 \ar@{->}[r] & O^\times \ar@{->}[r]\ar@{->}[d]_{\res}& K^\star \ar@{->}[r]_\val\ar@{->}[d]_{\rv}\ar@/_1.0pc/@{->}[dl]^\ac & \Gamma\ar@{->}[r] \ar@{=}[d] &  0 \\
1 \ar@{->}[r]& O^\times /1+\mathfrak{m} \simeq k^\star\ar@{->}[r]  & \RV^{\star} \ar@/_1.5pc/@{-->}[l]^{\ac_{\rv}}\ar@{->}[r]^{\val_{\rv}} & \Gamma \ar@{->}[r]  & 0 }$$ 

     One remarks that an angular component gives a section $ac_{\rv}:\RV \rightarrow k^\star$ and, as a consequence, the sort $\RV$ becomes isomorphic as a group to the direct product $\Gamma \times k^\star$.
    Such a map always exists in an $\aleph_1$-saturated valued field $\mathcal{K}$: as $\mathcal{O}^{\times}$ is a pure subgroup of $K^\star$, there is a section $s:\Gamma \rightarrow K^\star$  of the valuation (see Fact \ref{FactSectionPureSubgroupAleph1Saturated}). Then, the function $\ac:a \mapsto \res(a/s(v(a)))$ is an $\ac$-map. Any theory $T$ of Henselian valued fields in a language $\mathrm{L}_{\Gamma,k}$ admits a natural expansion -- denoted by $T_{\ac}$ -- in the language $\mathrm{L}_{\Gamma,k,\ac}=\mathrm{L}_{\Gamma,k} \cup \{\ac: K \rightarrow k\}$ \nomenclature[L]{$\mathrm{L}_{\Gamma,k\ac}=\mathrm{L}_{\Gamma,k} \cup \{\ac: K \rightarrow k\}$}{} by adding the axiom saying that $\ac$ is an angular component. 
    
    
    
    
    An important model theoretic property is relative quantifier elimination: 
    \begin{equation}
    \tag{$(\text{EQ})_{\Gamma,k,ac}$} \begin{array}{l}
        \text{ $T_{\ac}$ has quantifier elimination (resplendently)}  \\ 
          \text{relatively  to $\Gamma$ and $k$ in the language $\mathrm{L}_{\Gamma,k,\ac}$. }
          \end{array}
    \end{equation}
        \begin{equation}
    \tag{$(\text{EQ})_{\RV}$} \begin{array}{l}
        \text{ $T$ has quantifier elimination (resplendently)}  \\ 
          \text{relatively  to $\RV$ in the language $\mathrm{L}_{\RV}$. }
          \end{array}
    \end{equation}
     \nomenclature[P]{$(\text{EQ})_{\Gamma,k,ac}$,$(\text{EQ})_{\RV}$,$(\text{SE})_{\Gamma,k}$,$(\text{SE})_{\RV}$}{}
    
    Notice that according to the terminology in Paragraph \ref{SubsectionRelativeQuantifierElimination}, $\{\Gamma\},\{k\}$ and $\{\RV\}$ are closed sets of sorts. Then resplendency automatically follows from relative quantifier elimination (Fact \ref{FactClosedSortResplQE}).

    \begin{observation}
         $(\text{EQ})_{\RV}$ implies $(\text{EQ})_{\Gamma,k,ac}$.
        
    \end{observation}
    We include a proof of this observation for completeness.
    \begin{proof}
        We sketch a proof using the usual back-and-forth criterion. 
        We assume $(\text{EQ})_{\RV}$. Consider two models $\mathcal{M}=\{K_M,\Gamma_M,k_M\}$ and $\mathcal{N}=\{K_N,\Gamma_N,k_N\}$ of $T$ in the language $\mathrm{L}_{\Gamma,k,\ac}$, and a partial automorphism $f=(f_K,f_\Gamma,f_k):A=(K_A,\Gamma_A,k_A) \rightarrow B=(K_B,\Gamma_B,k_B)$ between a substructure $A\subseteq \mathcal{M}$ and a substructure $B\subseteq \mathcal{N}$. Moreover, we assume $f_k$ and $f_\Gamma$ to be elementary as morphisms respectively of fields and of ordered abelian groups. We want to extend $f$ to an elementary embedding of $M$  into $N$. By elementarity, we may extend $f_\Gamma$ (resp.\ $f_k$) to an elementary embedding of ordered abelian groups $\tilde{f}_\Gamma: \Gamma_M \rightarrow \Gamma_N$ (resp.\ to an elementary embedding of fields $\tilde{f}_k: k_M \rightarrow k_N$). Then, by studying quantifier-free formulas, one sees that $\tilde{f}=f \cup  \tilde{f}_\Gamma \cup \tilde{f}_k$ is a partial isomorphism of substructures. Without loss, assume that $\Gamma_A=\Gamma_M$ and $k_A=k_M$ and reset the notation.
        As the $\ac$-map induces a splitting of the exact sequence
        \[ 1 \rightarrow k^\star \overset{\iota}{\rightarrow} \RV^\star \overset{\val_{\rv}}{\rightarrow} \Gamma \rightarrow 0, \]
        we have the bijections  $\RV_M^\star \simeq  k_M^\star \times \Gamma_M$ and $\RV_N^\star \simeq  k_N^\star \times \Gamma_N$. Hence, the partial isomorphism $f$ induces an elementary embedding of $\RV$-structure $f_{\RV}: (\RV_M,\oplus, \cdot,\textbf{1},\textbf{0}) \rightarrow (\RV_N, \oplus , \cdot,\textbf{1},\textbf{0})$, and $f_K \cup f_{\RV}$ is a partial isomorphism of substructures in the language $\mathrm{L}_{\RV}$. By relative quantifier elimination down to $\RV$, $f_K \cup f_{\RV}$ extends to an elementary embedding $\tilde{f}=(\tilde{f}_K,f_{\RV})$ of $\{M,\RV_M\}$ into $\{N,\RV_N\}$.
        One sees that $\tilde{f}_K\cup f_\Gamma \cup f_k: \mathcal{M} \rightarrow \mathcal{N}$ is an embedding extending the original partial isomorphism $f$. By back-and-forth, $T$ satisfies $(\text{EQ})_{\Gamma,k,\ac}$.
    \end{proof}

    More specifically, we will have to study 1-dimensional definable sets $D \subset K$. Flenner showed in \cite{Fle11} that in Henselian valued fields of characteristic $0$, definable sets can be written with field-sorted linear terms (see Fact \ref{factfle}). This property will be of essential use. Let us give it also an abbreviation:
    
    \begin{definition}
       Let $T$ be the theory of a Henselian valued field $K$ in the language $\mathrm{L}_{\RV}$.
       We denote by $(\text{Lin})_{\RV}$ \nomenclature[P]{$(\text{Lin})_{\RV}$}{} the following property: any formula $\phi(x)$ with parameters in $K$ and with $\vert x \vert =1$ is equivalent to a formula of the form
        \begin{align}\label{FormulaRVLinearTerms}\phi_{\RV}(\rv(x-a_1),\ldots,\rv(x-a_r),\alpha)\end{align}
        where $r\in \mathbb{N}$ and $\phi_{\RV}$ is an $\RV$-formula with a tuple of parameters $\alpha \in \RV(K)$ and $a_1,\ldots,a_r \in K^1$.
    \end{definition}
    Notice that it is an improvement of a relative quantifier elimination down to $\RV$ for unary-definable sets: the term inside $\rv$ is linear in $x$ where $(\text{EQ})_{\RV}$ gives only a polynomial in $x$.
 
    Its algebraic counterpart seems to be the following:
    \begin{definition}
        A valued field is called \emph{algebraically maximal} if it admits no immediate algebraic extension. \index{Valued field! algebraically maximal}
    \end{definition}
    In particular, Henselian valued fields of equicharacteristic $0$ are algebraically maximal (by the fundamental equality) as well as algebraically closed valued fields.
    Delon proved that this is actually a first order property. For details, we refer to Delon's thesis \cite{Del} and a recent work of Halevi and Hasson in \cite{HH19}.
    
    We show now that algebraically maximal valued fields with quantifier elimination relative to $\RV$ enjoy the property $(\text{Lin})_{\RV}$. This fact was suggested by Yatir Halevi.
    Notice that a similar statement has been proved by Peter Sinclair for valued fields in the Denef-Pas language $\mathrm{L}_{\Gamma,k,ac}$ (see \cite[Theorem 2.1.1.]{Sin18}). We thank both of them for their enlightenment.

Let $\mathcal{K}=\{(K,\cdot,+,0,1), (\RV,\textbf{0},\textbf{1}, \cdot,\oplus), \rv:K^\star\rightarrow \RV^\star \}$ be a Henselian valued field viewed as a structure in the language $\mathrm{L}_{\RV}$. Let $\mathbb{K} \succ \mathcal{K}$ be a monster model. Let $x\in \mathbb{K} \setminus K$. We denote by $I_K(x)$ the set of values ${\{ \val(x-a) \vert \ a \in K\}}$. 
    \begin{fact}[Delon]
            $I_K(x)$ has no maximum if and only if the extension $K(x)/K$ is immediate.
    \end{fact}

    \begin{lemma}
        Assume that $\mathcal{K}$ is algebraically maximal. Let $x\in \mathbb{K} \setminus K$.
        \begin{itemize}
            \item If $I_K(x)$ has no maximum, let $(\gamma_i=\val(x-a_i))_{i\in I}$ be a co-final sequence of values in $I_K(x)$. Then 
            the quantifier-free type $\qftp(x/K)$ is implied by the type $\{\val(x-a_i)=\gamma_i\}_{i\in I}$. 
            \item If $I_K(x)$ has a maximum, then there is $a\in K$ such that $\tp(\rv(x-a)/\RV(K)$ determined $\qftp(x/K)$. Moreover, $\RV(K(x))$ is generated by $\RV(K)$ and $\rv(x-a)$.
        \end{itemize}
    \end{lemma}
    
    \begin{proof}
        Assume that $I_K(x)$ has no maximum, and let $(\gamma_i=\val(x-a_i))_{i\in I}$ be a co-final sequence of values in $I_K(x)$. Then, the sequence $(a_i)_{i\in I}$ is a pseudo-Cauchy sequence in $K$, with no pseudo-limit in $K$ and which pseudo-converges to $x$. The extension is immediate and $\mathcal{K}$ is algebraically maximal. By \cite{Kap42}, the pseudo-Cauchy sequence $(a_i)_{i\in I}$ is of transcendental type. Then, if $x'\in \mathbb{K}$ is another pseudo-limit of $(a_i)_{i\in I}$, the two extensions $K(x)$ and $K(x')$ are isomorphic over $K$. In other words, the quantifier-free type of $x$ over $K$ is uniquely determined by $\{\val(x-a_i)=\gamma_i\}_{i\in I}$.  \\
        Assume that $I_K(x)$ has a maximum $\gamma=\val(x-a)$. Then, we distinguish two cases.\\
             \textbf{Case 1:} We have that $\val(x-a) \in \Gamma_K$.\\ 
             Then, we have that $\rv(x-a)\notin \RV(K)$, as otherwise, it would exist $b\in K$ such that $\val(x-a-b)>\val(x-a)$, contradicting  the maximality of $\val(x-a)$. Let $c \in K$ such that $\val(x-a)=\val(c)$. Then, since $k_K$ is relatively algebraically closed in $k_{\mathbb{K}}$ (as $\mathbb{K}$ is an elementary extension of $K$) and since $\rv(\frac{x-a}{c})$ is in $k_{\mathbb{K}}\setminus k_K$, we have that $\res(\frac{x-a}{c})$ is transcendental over $k_K$. Without loss of generality, assume that $\frac{x-a}{c}=x$. The extension $K(x)$ is the Gauss extension, thus it is unique up to $K$-isomorphism (see e.g. \cite{EP05}). In particular, the quantifier-free type of $x$ over $K$ is uniquely determined by $\tp(\rv(x)/\RV(K))$. We show now that $\RV(K(x))$ is generated by $\RV(K)$ and $\rv(x)$ in the following sense:  
             Consider $P(X):=\sum_{i<n}a_iX^i$ a non-trivial polynomial in $K$, and assume that $a_{i_0},\ldots,a_{i_{k-1}}$ are the coefficient of minimum value. Since $\rv(x^n) \in k^\star_{\mathbb{K}}$ for all $n$, and since $\rv(x)$ is transcendental over $k_K$, we have:
             \[\rv\left(\frac{P(x)}{a_{i_0}}\right)= \sum_{j<k} \frac{\rv(a_{i_j})}{\rv(a_{i_0})}\rv(x^{i_j}) \in k^\star_{\mathbb{K}},\]
             and so
             \[\rv(P(x)) = \bigoplus_{j<k} \rv(a_{i_j})\rv(x)^{i_j}=\bigoplus_{i<n} \rv(a_{i})\rv(x)^{i}.\]
            \textbf{Case 2:} We have that $\val(x-a) \notin \Gamma_K$. Then for all $n \in \mathbb{N}^\star$, $n\cdot \val(x-a) \notin \Gamma_K$, as $\mathbb{K}$ is an elementary extension of $K$. Then, for any polynomial $P(x) \in K(x)$, $\val(P(x-a))$ can be expressed in terms of $\val(x-a)$. The isomorphism type of $K(x)$ over $K$ is uniquely determined by $\rv(x-a)=\alpha$, or in other words the quantifier-free type of $x$ over $K$ is uniquely determined by $\tp(\rv(x-a)/\RV(K))$. Without loss of generality, we may assume that $x=x-a$. One sees as well that $\RV(K(x))$ is generated by $\RV(K)$ and $\rv(x)$: consider $P(X):=\sum_{i<n}a_iX^i$ a non-trivial polynomial in $K$. As $n\cdot\val(x)$ is not in $\Gamma_K$ for all $n$, we have that:
            \[\rv(P(x)) = \bigoplus_{i<n} \rv(a_{i})\rv(x^{i}).\]
        
    \end{proof}

    \begin{theorem}\label{TheoremLinearTerms}
        Assume that $\mathcal{K}$ is algebraically maximal and admits quantifier elimination relative to $\RV$. Then it also satisfies the property $\text{(Lin)}_{\RV}$.
    \end{theorem}

    \begin{proof}
    By compactness, it is enough to show that any (complete) $1$-type $p(x)=\tp(b/K)$ over $K$ is determined by formulas of the form \ref{FormulaRVLinearTerms}. 
    \begin{itemize}
        \item If $p(x)=\tp(b/K)$ is a realised type, i.e. $b\in K$, then the type is determined by $\{\rv(x-b)=\mathbf{0}\}$.
        \item If $K(b)/K$ is immediate, then by the previous lemma, $\qftp(b/K)$ is determined by the type $\{\val(x-a_i)=\gamma_i\}_{i\in I}$, where $\gamma_i$ and $a_i$ are given by the previous lemma. This can be written in the language $\mathrm{L}_{\RV}$: for $i\in I$ choose $c_i\in K $ of value $\gamma_i$. Then \[\val(x-a_i)=\gamma_i \ \Leftrightarrow \ \rv(x-a_i) \oplus \rv(c_i) \neq \rv(x-a_i)\  \wedge \rv(x-a_i) \oplus \rv(c_i) \neq \rv(c_i). \]
        As $\RV(K(b))=\RV(K)$, and by quantifier elimination relative to the $\RV$-sort, $\{\val(x-a_i)=\gamma_i\}_{i\in I} \cup \tp(\emptyset /\RV(K))$ determines $p(x)=\tp(b/K)$.
        \item If $K(b)/K$ is non-immediate, then by the previous lemma, there is an $a\in K$ such that $\qftp(b/K)$ is determined by $q(\rv(x-a))$ where $q=\tp(\rv(b-a)/\RV)$. As $\RV(K(x))$ is generated by $\RV(K)$ and $\rv(b-a)$, and by quantifier elimination relative to $\RV$, we got that $q(\rv(x-a))$ determines $p(x)=\tp(b/K)$.
    \end{itemize}
    
    \end{proof}
    
    \begin{definition}
        A valued field of equicharacteristic $p>0$ is said Kaplansky if the value group is $p$-divisible, the residue field is perfect and does not admit any finite separable extensions of degree divisible by $p$.
    \end{definition}
        \begin{definition}\label{DefBenign} \index{Valued field! benign Henselian}
        Any $\{\Gamma\}$-$\{k\}$-enrichment of one of the following theories of Henselian valued fields is called \emph{benign}:

    \begin{enumerate}
        \item Henselian valued fields of characteristic $(0,0)$,
        \item algebraically closed valued fields,
        \item algebraically maximal Kaplansky Henselian valued fields.
    \end{enumerate}
        A model of a benign theory will be called a \emph{benign Henselian valued field}.
    \end{definition}

    As promised, we have:
      \begin{fact}\label{FactBenignTheories}
    Benign theories satisfy $(\text{EQ})_{\RV}$ and $(\text{Lin})_{\RV}$. 
    \end{fact}
    By the discussion above, they also satisfy $(\text{EQ})_{\Gamma,k,ac}$.

    \begin{proof} We just give examples of references for $(\text{EQ})_{\RV}$ and $(\text{EQ})_{\Gamma,k,ac}$. Notice that we might not refer to original proofs.
    The fact that Henselian valued fields of characteristic $(0,0)$ has property $(\text{EQ})_{\Gamma,k,ac}$ is the classical theorem of Pas. The proof that it has $(\text{EQ})_{\RV}$ is in \cite{Fle11}. Algebraically closed valued fields (in any characteristic) eliminate quantifiers by the theorem of Robinson. One deduces the property $(\text{EQ})_{\Gamma,k,ac}$ from it. One can find a proof that algebraically closed valued fields of any characteristic have $(\text{EQ})_{\RV}$ in \cite{HKR18}. Algebraically maximal Kaplansky valued fields have $(\text{EQ})_{\Gamma,k,ac}$ and $(\text{EQ})_{\RV}$ by \cite{HH19}. 
    As all these fields are algebraically maximal, they satisfy the condition of Theorem \ref{TheoremLinearTerms}, and thus enjoy the property $(\text{Lin})_{\RV}$.

    Finally, all these properties hold for any $\{\Gamma\}$-$\{k\}$-enrichment, as it is a particular case of $\{\RV\}$-enrichment, and as the sorts $\Gamma, k$ and $\RV$ are closed (Fact \ref{FactClosedSortResplQE}).
    
    \end{proof}

We will complete our study with some transfer principle for unramified mixed characteristic Henselian valued fields with perfect residue field. As it requires further techniques, it needs to be treated independently. We first need to introduce the Witt vector construction.\newline
    
\paragraph{\textbf{Witt vectors}}

    
     In the theory of unramified mixed characteristic Henselian valued fields, we will understand $\RV_{\delta_n}$-structures thanks to the well known Witt vector construction. We start by briefly recalling the definition. For more details, one can see \cite{Ser80}. Let $k$ be a field of characteristic $p$.
    
    \begin{definition}[Witt vectors] Fix $X_0,X_1,X_2, \dots$ and $Y_0,Y_1,Y_2, \ldots$ some indeterminates. We consider the polynomials $W_0,W_1,W_2,\ldots$ in $\mathbb{Z}[X_0, X_1, X_2, \ldots]$, called \textit{Witt polynomials} and defined by: 
    \[W_i= \sum_{i=0}^n p^iX_i^{p^{n-i}}.\]
The ring of \textit{Witt vectors over }$k$\index{Ring of Witt vectors}, denoted by $W(k)$\nomenclature[]{$W(k)$}{}, is a ring of base set $k^\omega$. The sum of $\bar{x}=(x_0,x_1,\ldots)$ and $\bar{y}=(y_0,y_1,\ldots)$ is given by: 
    \[ \bar{x}+\bar{y}= (S_n(x_0,\ldots,x_{n-1},y_0,\ldots,y_{n-1}))_n \]
    where $S_n(X_0,\ldots,X_{n-1})$ is the unique polynomial in $\mathbb{Z}[X_0, \ldots , X_{n-1}, Y_{0}, \ldots, Y_{n-1}]$  such that $$W_n(X_0,\ldots,X_{n-1})+W_n(Y_0,\ldots,Y_{n-1})=W_n(S_0,\ldots,S_{n-1}).$$
    The product $\bar{x} \cdot \bar{y}$ is defined similarly. These operations make $W(k)$ into a commutative ring. 
    
    The \textit{residue map} $\pi$ is simply the projection to the first coordinate. The natural section of the residue map, the so called \textit{Teichmüller lift}, is defined as follows: 
        \[\begin{array}{lrcl}
\tau: & k & \longrightarrow & W(k) \\
    & a & \longmapsto & [a]:=(a,0,0,\ldots) \end{array}.\] 
        
        Finally, all the above definitions make sense if we restrict the base-set to $k^n$. One gets then the \textit{truncated ring of Witt vectors of length} $n$ denoted by $W_n(k)$, as well as its Teichmüller map $\tau_n:  k \rightarrow W_n(k)$.
    \end{definition}

    \begin{observation}
        $(W_n(k),+,\cdot,\pi)$ is interpretable in the field $(k,+,\cdot,1,0)$, with base set $k^n$. It is clear that $\bdn(W_n(k)):=\kappa_{inp}^1(W_n(k)) \leq \kappa_{inp}^n(k)$ (we will show that they are in fact equal).
    \end{observation}
    
    Recall that a $p$-ring is a complete local ring $A$ of maximal ideal $pA$ and perfect residue field $A/pA$. It is strict if $p^n\neq 0$ for every $n \in \mathbb{N}$. Here are some basic facts about Witt vectors:
    \begin{fact}[{see e.g. \cite[Chap. 6]{vdD14}}] \label{FactWittVectors}
    Recall that the field $k$ is perfect.
        \begin{enumerate}
            \item  The ring of Witt vectors $W(k)$ is a strict local $p$-ring of residue field $k$, unique up to isomorphism with these properties.
            \item The Teichmüller map is given by the following:  let $a\in k$, then $\tau(a)$ is the limit (for the topology given by the maximal ideal $pW(k)$ ) of any sequence $(a_n^{p^n})_{n < \omega}$ such that $\pi(a_n)^{p^n}=a$ for all $n$.
            \item In particular $\tau_n$ is definable in the structure $(W_n(k),+,\cdot,\pi)$. Indeed  $\tau_n(a)$ is the (unique) element $a_{n-1}^{p^{n-1}}\in W_n$ such that $\pi(a_{n-1})^{p^{n-1}}=a$.
            \item for $x=(x_{n})_{n<\omega} \in W(k)$, one has $x=\sum_{n<\omega}[x_n^{p^{-n}}]p^{n}.$
            \item In particular, the map $\chi_{i}: W(k) \rightarrow k, \ x=(x_0,x_1,\ldots)\mapsto x_i$ is definable in the structure $(W(k),+,\cdot,\pi)$. One has indeed 
            \[x_i= \pi\left(\frac{x-p^{i-1}[x_{i-1}^{p^{1-i}}]- \cdots - p[x_1^{p^{-1}}]-[x_0]}{p^{i}}\right)^{p^i}.\]
            Similarly, for $0\leq i \leq n-1$ the map $\chi_{n,i}: W_n(k) \rightarrow k, \ x=(x_0,x_1,\ldots, x_{n-1})\mapsto x_i$ is definable in $(W_n(k),+,\cdot,\pi)$. 
        \end{enumerate}

    \end{fact}
    
    We deduce the following:
     \begin{corollary}\label{CorollaryWnBiInterpretableWithk}
        \begin{itemize}
            \item The structure $(W_n(k),+,\cdot,\pi:W_n(k) \rightarrow k)$ is bi-interpretable on unary sets with the structure $(k^{n},k,+,\cdot, p_i, i<n)$, where $p_i:k^{n}\rightarrow k, \ (x_0,\ldots,x_{n-1})\mapsto x_i$ is the projection map. In other words, there is a bijection $W_n(k)\simeq k^{n}$ which leads to identify definable sets. 
            \item Similarly, the structure $(W(k),+,\cdot,\pi:W(k) \rightarrow k)$ is bi-interpretable on unary sets (with no parameters) with the structure $(k^{\omega},k,+,\cdot, p_i, i<\omega)$ where $p_i:k^{\omega}\rightarrow k, \ (x_0,x_1,\ldots)\mapsto x_i$.
            \end{itemize}
    \end{corollary}
    This corollary, coupled with Fact \ref{FactBdnInterpretUnarySet}, will be one of the main argument to treat our reduction principle in the context of unramified mixed characteristic valued fields with perfect residue field.

        

\paragraph{\textbf{Unramified mixed characteristic Henselian valued fields}}\label{Preliminaries Unramified mixed characteristic Henselian valued fields}

    We give a short overview on unramified valued fields, by presenting the similarities with benign valued fields.
    The (partial) theory of Henselian valued fields of characteristic $0$ does not satisfy either $(\text{EQ})_{\Gamma,k,ac}$ or $(\text{EQ})_{\RV}$. We indeed need to get `information' modulo $\mathfrak{m}_{\delta_n}$ in a quantifier-free way. 
    Let us recall the leading term language of finite order: 
    \[\mathrm{L}_{\RV_{<\omega}}=\left\{ K,(\RV_{\delta_n})_{n<\omega},(\oplus_{\delta_l,\delta_m, \delta_n})_{n<l,m},(\rv_{\delta_n})_{n<\omega}, (\rv_{\delta_n\rightarrow \delta_m})_{m<n<\omega} \right\},\]
      where $\oplus_n$ are ternary relation symbols and $\delta_n=\val(p^n)$. Let us just define all the analogous properties:
     \begin{align} \tag*{$(\text{AKE})_{\RV_{<\omega}}$}  \mathcal{K},\mathcal{K}' \models T,  \mathcal{K}\subseteq \mathcal{K}', \text{ we have } \mathcal{K} \preceq\mathcal{K}' \Leftrightarrow \RV_{<\omega} \preceq \RV_{<\omega}'.
    \end{align} \nomenclature[P]{$(\text{AKE})_{\RV_{<\omega}}$}{}

    Let us cite two main results in \cite{Fle11}. First we have:
    \begin{fact}\cite[Proposition 4.3]{Fle11}\label{FactRelQERV}
        Let $T$ be the theory of characteristic 0 Henselian valued fields in the language $\mathrm{L}_{\RV_{<\omega}}$. Then $T$ eliminates field-sorted quantifiers. 
                    \begin{center}
                \begin{tikzpicture}
                    \node{$\mathcal{K}$ }
                        child { node {$\RV_{<\omega}$}
                        };
                \end{tikzpicture}
            \end{center}
    \end{fact}

    This result was already proved in \cite{Bas91}. Again, an important consequence is that the multisorted substructure  $\left((\RV_{\delta_n})_{n<\omega},(\oplus_{\delta_l,\delta_m, \delta_n})_{n<l,m}, (\rv_{\delta_n\rightarrow \delta_m})_{m<n<\omega} \right)$ is stably embedded and pure. Secondly, we have its one-dimensional improved version: 
    
    \begin{fact}\label{factfle}\cite[Proposition 5.1]{Fle11}
	Let $T$ be the theory of Henselian valued fields $K$ of characteristic $0$ in the language $\mathrm{L}_{\RV_{<\omega}}$.
       It has the following property denoted by $(\text{Lin})_{\RV_{<\omega}}$\nomenclature[P]{$(\text{Lin})_{\RV_{<\omega}}$}{}: any formula $\phi(x)$ with parameters in $K$ and with $\vert x \vert =1$ is given by a formula of the form
        \begin{align}\label{FormulaRV<omegaLinearTerms}\phi_{\RV_{\delta_n}}(\rv(x-a_1),\ldots,\rv(x-a_r),\alpha)\end{align}
        where $\phi_{\RV_{\delta_n}}(\mathbf{x}_1,\ldots, \mathbf{x}_r, \mathbf{y})$ is an $\RV_{\delta_n}$-formula, with a tuple of parameters $\alpha \in \RV_{\delta_n}(K)$ and $a_1,\ldots,a_r \in K$ and $r\in \mathbb{N}$.
    \end{fact}

    Again, notice that the improvement comes from the fact that the term inside $\rv_{\delta_n}$ is linear in $x$ where Fact \ref{FactRelQERV} gives only a polynomial in $x$. These theorems also include the case of equicharacteristic $0$, and it gives the same result as cited above. Indeed, in equicharacteristic $0$, we may identify $\bigcup_{n<\omega}\RV_{\delta_n}$ with $\RV=\RV_0$ (Remark \ref{RemarkIdentifyRVRVOmega}).
    We continue with a remark on enrichment.
    \begin{remark}
        Fact \ref{FactRelQERV} above holds in any $\RV_{<\omega}$-enrichment of $\mathrm{L}_{\RV}$. Indeed, first note that the $\RV_{<\omega}$-sort is closed in the language $\mathrm{L}_{\RV_{<\omega}}$, \textit{i.e.} any relation symbol involving a sort $\RV_{\delta_n}$ or any function symbol with a domain involving a sort $\RV_{\delta_n}$ only involves such sorts. By Fact \ref{FactClosedSortResplQE}, the theory $T$ of Henselian valued fields of characteristic $0$ also eliminates quantifiers resplendently relative to $\RV_{<\omega}$. In other words, given an $\RV_{<\omega}$-enrichment $\mathrm{L}_{\RV_{<\omega},e}$, any complete $\mathrm{L}_{\RV_{<\omega},e}$-theory $T_e \supset T$ eliminates field-sorted quantifiers. A careful reading of Flenner's proof give us that Fact \ref{factfle} holds resplendently as well. 
        \end{remark}

    Now, let us discuss more specifically on the unramified mixed characteristic cases. We denote by $T$ the theory of unramified mixed characteristic Henselian valued fields with perfect residue field. We assume now that $\mathcal{K}$ is such a valued field. There is by definition a smallest positive value $\val(p)=\delta_1$, that we denote by $1$. 
    
    \textbf{Notation:}
    Notice that $\mathfrak{m}=p\mathcal{O}$ and in general that $\mathfrak{m}_{\delta_n}=\mathfrak{m}^{n+1}=p^{n+1}\mathcal{O}$ for all $n\geq 0$. We will write $\mathfrak{m}^{n+1}$ instead of $\mathfrak{m}_{\delta_n}$,  $\mathcal{O}_{n+1}$  instead of $\mathcal{O}_{\delta_n}$ and  $\RV_{n+1}$ instead of $\RV_{\delta_{n}}= K^\star/ 1+p^{n+1}\mathcal{O}$ \nomenclature[]{$\RV_n$}{$n^{\text{th}}$  $\RV$-sort}. The projection map $\res_{\delta_n}: \mathcal{O} \rightarrow \mathcal{O}_{\delta_n}$ is written $\res_{n+1}: \mathcal{O} \rightarrow \mathcal{O}_{n+1}$ etc. The idea is to denote by $\RV_n$ the $n^{\text{th}}$ $\RV$-sort, as this makes sense in unramified (or finitely ramified) mixed characteristic valued fields. The purpose is also to fit with the usual notation, and it will help to simplify the notation, although this convention contradicts the previous one ($\RV_0$ where $0$ stands for the value $0\in \Gamma$ is now $\RV_1$, the first $\RV$-sort).
    
     In this context, let us define the angular component of degree $n$:

        \begin{definition}\label{DefinitionAngularComponentsMixedCharacteristic}
        Let $n$ be an integer greater than $0$. \textit{An angular component of order} $n$ is a homomorphism $\ac_n:  K^\star \rightarrow \mathcal{O}_n^\times$ such that for all $u\in \mathcal{O}^\times$, $\ac_n(u)=res_n(u)$. A system of angular component maps $(ac_n)_{n<\omega}$ is said to be \textit{compatible} if for all $n$, $\pi_n \circ \ac_{n+1}=\ac_n$\index{Angular component! of order $n$}\nomenclature[]{$\ac_n$}{Angular component of order $n$} where $\pi_n:\mathcal{O}_{n+1}\rightarrow \mathcal{O}_n \simeq \mathcal{O}_{n+1}/p^n\mathcal{O}_{n+1}$ is the natural projection.
        

    \begin{center}
        \begin{tikzpicture}[line cap=round,line join=round,>=triangle 45,x=1cm,y=1cm,scale=0.3]\clip(0.464749615349495,-3.887764486968447) rectangle (15.782395059865305,13.80598253132582);\draw [->,line width=1pt,>=stealth] (7.873051286030098,4.136247085440175) -- (2.7509254230514903,0.27287083435127624);\draw [->,line width=1pt,>=stealth] (7.873051286030098,4.136247085440175) -- (6,-2);\draw [->,line width=1pt,>=stealth] (7.873051286030098,4.136247085440175) -- (9.964516505540304,-1.929033011080608);\draw [->,line width=1pt,>=stealth] (7.873051286030098,4.136247085440175) -- (13.130164383515968,0.45866482342147963);\draw [->,line width=1pt,>=stealth] (7.873051286030098,4.136247085440175) -- (14.28777374095427,4.251074085590988);\draw [->,line width=1pt,>=stealth] (7.873051286030098,4.136247085440175) -- (12.995177149008702,7.999623336529074);\draw [->,line width=1pt,>=stealth] (7.873051286030098,4.136247085440175) -- (9.746102572060193,10.27249417088035);\draw [->,line width=1pt,>=stealth] (2.7509254230514903,0.27287083435127624) -- (6,-2);\draw [-,line width=1pt,dotted] (6,-2) -- (9.964516505540304,-1.9290330110806082);\draw [->,line width=1pt,>=stealth] (9.964516505540304,-1.9290330110806082) -- (13.130164383515968,0.45866482342147985);\draw [->,line width=1pt,>=stealth] (13.130164383515968,0.45866482342147963) -- (14.28777374095427,4.251074085590988);\draw [->,line width=1pt,>=stealth] (14.28777374095427,4.251074085590988) -- (12.995177149008702,7.999623336529074);\draw [->,line width=1pt,>=stealth] (12.995177149008702,7.999623336529074) -- (9.746102572060193,10.27249417088035);\begin{scriptsize}\draw [fill=black] (6,-2) circle(4pt);\draw[color=black] (6.296372157630452,-2.645286355089424) node {$\mathcal{O}_n^\times$};
        \draw [fill=black] (9.964516505540304,-1.9290330110806082) circle(4pt);
        \draw[color=black] (10.734854901670144,-2.4383622178874402) node {$\mathcal{O}_5^\times$};\draw [fill=black] (13.130164383515968,0.45866482342147963) circle(4pt);\draw[color=black] (14.221539946504546,0.58344926494364) node {$\mathcal{O}_4^\times$};\draw [fill=black] (14.28777374095427,4.251074085590988) circle(4pt);\draw[color=black] (15.194110057315793,4.741932008983317) node {$\mathcal{O}_3^\times$};\draw [fill=black] (12.995177149008702,7.999623336529074) circle(4pt);\draw[color=black] (13.538716959699915,8.494263027426962) node {$\mathcal{O}_2^\times$};
        \draw [fill=black] (9.746102572060193,10.27249417088035) circle(4pt);
        \draw[color=black] (10.296905476868822,11.770428536648783) node {$\mathcal{O}_{1}^\times= k^\star$};\draw [fill=black] (2.7509254230514903,0.27287083435127624) circle(4pt);\draw[color=black] (1.881990563988111,-0.02424728386429531) node {$\mathcal{O}_{n+1}^\times$};\draw [fill=black] (7.873051286030098,4.136247085440175) circle (6pt);\draw[color=black] (6.710220432034421,5.4247549957879455) node {$K^\star$};\draw[color=black] (4.296105498011267,2.710640061764801) node {$\ac_{n+1}$};
        \draw[color=black] (5.782523883226482,1.086272251748269) node {$\ac_n$}
        ;\draw[color=black] (8.306638817249637,0.4003734021456239) node {$\ac_5$};\draw[color=black] (11.079728463673452,1.000446633542204) node {$\ac_4$};\draw[color=black] (11.814349149683377,3.54526074777472) node {$\ac_3$};\draw[color=black] (11.46947558768007,5.73527453140051) node {$\ac_2$};\draw[color=black] (8.745107777663529,7.976819354612412) node {$\ac_1=\ac$};\draw[color=black] (4.13282511206635,-1.199691202715693) node {$\pi_n$};
        \draw[color=black] (11.917944548093962,-1.29139944078709081) node {$\pi_4$};\draw[color=black] (14.495388220908516,2.4347412121621555) node {$\pi_3$};\draw[color=black] (14.435388220908516,6.080148093403816) node {$\pi_2$};\draw[color=black] (11.945374437282715,10.166833138238204) node {$\pi_1=\pi$};\end{scriptsize}\end{tikzpicture}

    \end{center}
    \end{definition}
     The convention is to contract $\ac_1$ to $\ac$ and $\pi_1$ to $\pi$. Then, let us complete the diagram given in Paragraph \ref{SusubsectionRVSort}:  
    \[\xymatrix{
        1 \ar@{->}[r] & \mathcal{O}^\times \ar@{->}[r]\ar@{->}[d]_{\res_n}& K^\star \ar@{->}[r]_\val\ar@{->}[d]_{\rv_n}\ar@/_1.0pc/@{->}[dl]^{\ac_n} & \Gamma\ar@{->}[r] \ar@{=}[d]\ar@/_1.0pc/@{-->}[l]^s &  0 \\
        1 \ar@{->}[r]& \mathcal{O}^\times /(1+\mathfrak{m}^n) \simeq \mathcal{O}_n^\times\ar@{->}[r]  & \RV_n^{\star} \ar@{->}[r]^{\val_{\rv_n}} & \Gamma \ar@{->}[r]  & 0 }\] 
    
    A section $s:\Gamma \rightarrow K^\star$ of the valuation gives immediately a compatible system of angular components (defined as $\ac_n:= a\in K^\star \mapsto \res_n(a/s(v(a)))$). As $\mathcal{O}^\times$ is a pure subgroup of $K^\star$, such a section exists when $\mathcal{K}$ is $\aleph_1$-saturated (see Fact \ref{FactSectionPureSubgroupAleph1Saturated}). As always, we assume that $\mathcal{K}$ is sufficiently saturated and we fix a compatible sequence $(\ac_n)_n$ of angular components.  
    
    We denote by $T_{ac_{<\omega}}$ the extension of $T$ to the language $\mathrm{L}_{\Gamma,k,\ac_{<\omega}}:=\mathrm{L}_{\Gamma,k} \cup \{\mathcal{O}_n, \ac_n: K \rightarrow \mathcal{O}_n, \ n\in \mathbb{N}\}$ \nomenclature[L]{$\mathrm{L}_{\Gamma,k,\ac_{<\omega}}:=\mathrm{L}_{\Gamma,k} \cup \{\mathcal{O}_n, \ac_n: K \rightarrow \mathcal{O}_n, \ n\in \mathbb{N}\}$}{} where $\ac_n$ are interpreted as compatible angular components of degree $n$ (see for instance \cite{AJ19}).

    The following proposition is well known and has been used for example in \cite[Corollary 5.2]{Bel99}. It states how the structure $\RV_n$ and the truncated Witt vectors $W_n$ are related.
    \begin{proposition}\label{KerValRV_n}
         \begin{enumerate}
             \item The residue ring $\mathcal{O}_n$ of order $n$  is isomorphic to $W_n(k)$, the set of truncated Witt vectors of length $n$ .
             \item The kernel of the valuation $\val:  \RV_n^\star \rightarrow \Gamma$ is given by $\mathcal{O}^\times/(1+\mathfrak{m}^n) \simeq (\mathcal{O}/\mathfrak{m}^n)^\times$. It is isomorphic to $W_n(k)^\times$, the set of invertible elements of $W_n(k)$.
         \end{enumerate}

    \end{proposition}
    \begin{proof}
        It is clear that (2) follows from (1) as $\mathcal{O}^\times/(1+\mathfrak{m}^n) \simeq (\mathcal{O}/\mathfrak{m}^n)^\times$.
        
        Now, we prove (1) for any discrete value group $\Gamma$. Consider 
        \[W': = \varprojlim\limits_{n< \omega} \mathcal{O}_n \subset \prod\limits_{n<\omega}\mathcal{O}_n\]
        the inverse limit of the $\mathcal{O}_n$'s. It is:
        \begin{itemize}
            \item strict, \textit{i.e.} $p^n \neq 0$ in $W'$ for every $n< \omega$, as $\pi_{n+1}(p^n) \neq 0$ in $\mathcal{O}_{n+1}=\mathcal{O}/p^{n+1}\mathcal{O}$,
            \item local, as a projective limit of the local rings $\mathcal{O}_n$,
            \item a $p$-ring. Its maximal ideal is $pW'$, it is complete as projective limit, and its residue field is the perfect field $k$.
        \end{itemize}
        By uniqueness, $W'$ is isomorphic to $W(k)$, the ring of Witt vectors over $k$. One just has to notice that $W'/p^nW'\simeq \mathcal{O}/p^n\mathcal{O}$ and it follows easily that $\mathcal{O}_n \simeq W_n(k)$ for every $n<\omega$. 
    \end{proof}
    \begin{note}
        In the above proof, one can also recover $W(k)$ by considering the coarsening $\dot{K}$ of $K$ by the convex subgroup $\mathbb{Z}\cdot 1$. Indeed, if we denote by $K^\circ$ the residue field of the coarsening, as $\mathcal{K}$ is saturated enough, one has $\varprojlim_{n< \omega}\mathcal{O}_n\simeq \mathcal{O}(K^\circ)$ (see \cite{Bas91}). 
    \end{note}

    \begin{fact}[Bélair \cite{Bel99}]\label{FactBelair}
     The theory $T_{ac_{<\omega}}$ of Henselian mixed characteristic valued fields with perfect residue field and with angular components eliminates field-sorted quantifiers in the language $\mathrm{L}_{\Gamma,k,\ac_{<\omega}}$\footnote{The Ax-Kochen-Ershov property and relative quantifier elimination for Henselian unramified mixed characteristic valued fields (with possibly imperfect residue field) has been proved in \cite{AJ19}.}  
 .
    \end{fact}
    Notice that in \cite{Bel99}, Bélair doesn't assume that the residue field $k$ is perfect, but it is indeed necessary in order to identify the ring $\mathcal{O}_n:=\mathcal{O}/\mathfrak{m}^n$ with the truncated Witt vectors $W_n(k)$.
    This implies as well that the residue field $k$ and the value group $\Gamma$ are pure sorts, and are orthogonal. This can be seen by analysing field-sorted-quantifier-free formulas, and by noticing that $\mathcal{O}_n \simeq W_n(k)$ is interpretable in $k$ (Corollary \ref{CorollaryWnBiInterpretableWithk}).

    By analogy with the previous paragraph, we name the following properties: 
    
    \begin{align}
    \tag*{$(\text{EQ})_{\Gamma,k,ac_{<\omega}}$} \text{ $T_{ac_{<\omega}}$ eliminates $K$-sorted quantifiers in the language $\mathrm{L}_{\ac_{<\omega}}$}.\\
    \tag*{$(\text{EQ})_{\RV_{<\omega}}$} \text{ $T$ has quantifier elimination (resplendently) relatively to $\RV_{<\omega}$}.
    \end{align}
    \nomenclature[P]{$(\text{EQ})_{\Gamma,k,ac_{<\omega}}$,$(\text{EQ})_{\RV_{<\omega}}$}{}
    Again, notice that $\RV_{<\omega}=\bigcup_{n<\omega}\RV_n$ is a closed set of sorts.
    
    To sum up, we have:

     \begin{fact}\label{FactUnramifiedMixedCharacteristicValuedfields}
    The theory of unramified mixed characteristic Henselian valued fields with perfect residue field satisfies $(\text{EQ})_{\RV_{<\omega}}$, $(\text{EQ})_{\Gamma,k,ac_{<\omega}}$, $(\text{Lin})_{\RV_{<\omega}}$. 
    \end{fact}

\subsubsection{Abelian groups}\label{SubsectionPreliminariesAbelianGroups}
We conclude these preliminaries with some facts on abelian groups. They will be used in Section \ref{SectionExactSequence}.  We are specifically interested in abelian groups for one main reason: we have to understand the structure of pure short exact sequences of abelian groups in order to produce our reduction principles for benign Henselian valued fields. As we obtain also reduction principles for such short exact sequences, we will take the occasion to apply it on explicit examples. 

\begin{notation}We recall some standard notation:
    \begin{itemize}
    \item $\mathbb{Z}(p^n)$ is the cyclic group of $p^n$ elements, \nomenclature[]{$\mathbb{Z}(p^n)$}{Cyclic group of $p^n$ elements}
    \item $\mathbb{Z}_{(p)}$ is the additive group of the integers localised in $(p)$, \nomenclature[]{$\mathbb{Z}_{(p)}$}{Integers localised in $(p)$}
    \item $\mathbb{Z}(p^\infty)$ is the Prüfer $p$-group, \nomenclature[]{$\mathbb{Z}(p^\infty)$}{Prüfer $p$-group}.

\end{itemize}   
    For an abelian group $A$, we denote by $A^{(\omega)}$ the direct product of $\omega$ copies of $A$.
\end{notation}

We have:

    \begin{fact} [{ \cite[Theorem 2 ${\mathbb{Z}}$ 1 ]{Pre88}} ]
    
        Let $A$ be an abelian group. Let $P_n$ be a predicate for $n$-divisibility in $A$. Then $\{A, +, - ,0 , \{P_n\}_{n\in \mathbb{N}_{>1}}\}$ eliminates quantifiers.
    \end{fact}

    The burden (or equivalently by stability, the dp-rank) of pure abelian groups has been computed in terms of their \emph{Szmielew invariants} by Halevi and Palacín in \cite{HP17}. We borrow from their work the following proposition, which says that a useful criterion to witness inp-patterns is a characterisation in the case of one-based groups (and in particular, in the case of unenriched abelian groups):

    \begin{proposition}[{\cite[Proposition 3.4]{HP17}}]\label{PropCriBdnOneBasAbeGrp}
        A stable one-based group admits an inp-pattern of depth $\kappa$ if and only if there exists $\acl^{eq}(\emptyset)$-definable subgroups $(H_{\alpha})_{\alpha<\kappa}$ such that for any $i_0<\kappa$, one has:
        \[ \left[\bigcap_{\alpha \neq i_0}H_\alpha \ : \ \bigcap_{\alpha}H_\alpha \right] = \infty.\]
        If $(b_{\alpha,j})_{j<\omega}$ are representatives of pairwise distinct classes of $ \bigcap_{\alpha \neq i_0}H_\alpha$ modulo $\bigcap_{\alpha}H_\alpha$, an inp-pattern of depth $\kappa$ is given by $\{x\in b_{\alpha,j} H_\alpha\}_{\alpha<\kappa,j<\omega}$.
    \end{proposition} 
    We will use this criterion to provide examples to Theorem \ref{ThmBdnExSeq}.

    \paragraph{\textbf{Quantifier elimination result in pure short exact sequences}}\label{SubSubSecQEExaSeqAbeGrp}
    
    \begin{definition}
        Let $B$ a group and $A$ a subgroup. We say that $A$ is a \emph{pure subgroup} \index{Pure!subgroup} of $B$ if for all $a$ in $A$, $n\in \mathbb{N}$, $a$ is $n$-divisible in $B$ if and only if $a$ is $n$-divisible in $A$.
    \end{definition}
    We recall the following fundamental fact:
    \begin{fact}\label{FactSectionPureSubgroupAleph1Saturated}
        Let $\mathcal{M}$ be an $\aleph_1$-saturated structure, and let $A,B$ be two definable abelian groups, and assume that $A$ is a pure subgroup of $B$. Then the exact sequence of abelian groups $0\rightarrow A \rightarrow B \rightarrow B/A \rightarrow 0$
        splits: there is a group homomorphism $\alpha:B \rightarrow A$ such that $\restriction{\alpha}{A}$ is the identity on $A$. In such case, $B$ is isomorphic as a group to $A \times B/A$.
    \end{fact}
    More precisely, it is an immediate corollary of a more general statement on \emph{pure-injectivity}. See \cite[Theorem 20 p.171]{Che76}.
    
    Assume that we have a pure short exact sequence of abelian groups
    \[ \xymatrix{0 \ar[r] & A \ar[r]^{\iota}& B \ar[r]^{\nu} & C\ar[r] &0}. \]
    (meaning that $\iota(A)$ is a pure subgroup of $B$)\index{Pure! short exact sequence of abelian groups}. We treat it as a three-sorted structure $(A,B,C,\iota,\nu)$, with a group structure for all sorts. In fact, in our main applications, we will consider such a sequence with more structure on $A$ and $C$. Let us explicitly state all results resplendently, by working in an enriched language. So, let $\mathcal{M}=(A,B,C,\iota,\nu,\ldots)$ be an $\{A\}$-enrichment of a $\{C\}$-enrichment  (for short: an $\{A\}$-$\{C\}$-enrichment) of the exact sequence in a language that we will denote by $\mathrm{L}$, and we denote its theory by $T$. We will always assume that $\mathcal{M}$ is sufficiently saturated ($\aleph_1$ saturated will be enough).\\
    Hypothesis of purity implies the exactness of the following sequences for $n\in \mathbb{N}$:
    \[\xymatrix{0 \ar[r] & A/nA \ar[r]^{\iota_n}& B/nB \ar[r]^{\nu_n} & C/nC\ar[r] &0}.\]
    One has indeed that 
    \[\frac{A+nB}{nB} \simeq \frac{A}{A\cap nB} = \frac{A}{nA}.\]
    
    We consider for $n\geq 0$ the following maps: \\
    \begin{itemize}
        \item the natural projections $\pi_n:A \rightarrow A/nA$,
        \item the map  \[\begin{array}{ccccc}
\rho_n  : & B & \to & A/nA \\
  & b & \mapsto & 
        \begin{cases}
            0_{A/nA} \ \text{ if }b\notin \nu^{-1}(nC) \\
            \iota_n^{-1}(b+nB) \text{ otherwise, }
        \end{cases} 
        \end{array}\]
        \end{itemize}

    where $0_{A/nA}$ is the zero element of $A/nA$ (often denoted by $0$).
    Then let us consider the language 
    \[\mathrm{L}_q= \mathrm{L} \cup \{A/nA,\pi_n,\rho_n\}_{n\geq 0},\] \nomenclature[L]{$\mathrm{L}_q$}{}
    and let $T_q$ be the natural extension of the theory $T$. 
    By $<A>$, we denote the set of sorts containing $A$, $A/nA$ and the new sorts possibly coming from the $A$-enrichment. Similarly, let $<C>$ be the set of sorts containing $C$ and the new sorts possibly coming from the $C$-enrichment. By $A$ or $A$-sort and $C$ 
    or $C$-sort, we will abusively refer to $<A>$ and $<C>$ respectively, and similarly for $A$-formulas and $C$-formulas.   Aschenbrenner, Chernikov, Gehret and Ziegler prove the following result:
    
    \begin{fact}[{\cite[Theorem 4.2]{ACGZ20}}]\label{FactACGZ}
        The theory $T_q$ (resplendently) eliminates $B$-sorted quantifiers.
        
                    \begin{center}
                \begin{tikzpicture}
                    \node{$\mathcal{M}$ }
                        child {  node {$A$}}
                        child { node {$C$}
                        };
                \end{tikzpicture}
            \end{center}

    More precisely, all $\mathrm{L}_q$-formulas $\phi(x)$ with a tuple of variables  $x\in B^{\vert x\vert}$ are equivalent to boolean combinations of formulas of the form:
    \begin{enumerate}
        \item $\phi_C(\nu(t_0(x)),\ldots,\nu(t_{s-1}(x)))$ where $t_i(x)$'s are terms in the group language, and $\phi_C$ is a $C$-formula,
        \item $\phi_{A}(\rho_{n_0}(t_0(x)),\ldots,\rho_{n_{s-1}}(t_{s-1}(x)))$ where the $t_i(x)$'s are terms in the group language, where \newline $s,n_0,n_1,\ldots,n_{s-1}\in \mathbb{N}$, and where $\phi_{A}$ is an $A$-formula.
    \end{enumerate}
    In particular there is no occurrence of the symbol $\iota$.
    
    \end{fact}
    
    In particular, notice that the formula $t(x)=0$ is equivalent to $\nu(t(x))=0 \wedge \rho_0(t(x))=0$ and $\exists y \ ny=t(x)$ is equivalent to $\exists y_C \ ny_C=\nu(t(x)) \wedge \rho_n(t(x))=0$.
    
    We have:
	    \begin{corollary}\label{CoroPureOrth}	
	        In the theory $T_q$, $<A>$ and $<C>$ are stably embedded, pure (see Definition \ref{DefinitionPureSort}) and orthogonal to each other.	
	    \end{corollary}
    
    In fact, it can be easily deduced from the existence of a section. The following proof is more technical but highlights the fact that one does not need the function $\iota$ in order to express definable sets in $\bigcup_{n<\omega} A/nA$.
    
    \begin{proof}
        In this proof, $A$ (resp.\ $C$) abusively refer to the union of the sorts $<A>$ (resp.\ $<C>$). The $C$-sort is pure and stably embedded by Fact \ref{FactACGZ} and closedness of $C$. It is also clear for the sort $A$, even if $A$ is not a closed sort\footnote{The proof shows that one does not need the function $\iota: A \rightarrow B$ in order to describe definable sets in $A$. In a certain sense, $<A>$ is a `closure' of $A$, as it describes the induced structure on $A$, with no resort to any symbol from $\mathrm{L}_q\setminus \restriction{\mathrm{L}_q}{<A>}$.}: one only needs to deal with the map $\iota:A \rightarrow B$. If $D$ is a definable set in $A^{\vert x_A\vert}$, it is given by a disjunction of formulas of the form
        \begin{align*}
         \phi(x_A)= &  \phi_{A}\left(\rho_{n_0}(k_0\iota(t_0(x_A))+b_0), \ldots, \rho_{n_{s-1}}(k_{s-1}\iota(t_{s-1}(x_A))+b_{s-1}), a\right) \\
         &\wedge 
          \phi_C\left(\nu(\iota(t_0(x_A))), \ldots, \nu(\iota(t_{s-1}(x_A))), c\right).  
        \end{align*}
        where $x_A$ is a tuple of $A$-variables, the $t_i(x_A)$'s are terms in the group language, $s,k_0, \ldots, k_{s-1}\in \mathbb{N}$, $b_0,\ldots, b_{s-1} \in B$, and  $a \in A$ and $c\in C$ are tuples of parameters (notice that we also used that $\iota$ and $\nu$ are morphisms). We apply now the following transformation in order to get a new formula $\phi'(x_A)$:
        \begin{itemize}
            \item For $l<s$, if $\nu(b_l) \notin n_lC$, then  replace $\rho_n(k_l\iota(t_l(x_A))+b_l)$ by $0_{A/n_lA}$.
            \item For $l<s$, if $\nu(b_l) \in n_lC$, replace $\rho_{n_l}(k_l\iota(t_l(x_A))+b_l)$ by $k_l\pi_{n_l}(t_l(x_A))+\rho_{n_l}(b_l)$. 
            \item Replace $\nu(\iota(t_l(x_A))$ by $0_C$.
        \end{itemize}  We obtain a pure $A$-formula $\phi'(x_A)$ such that $\phi'(A^{\vert x_A\vert})=\phi(A^{\vert x_A\vert})$. Orthogonality can also be proved similarly.
    \end{proof}

\section{Burden of pure short exact sequences of abelian groups}\label{SectionExactSequence}
    We prove in this section that the burden of a pure short exact sequence of abelian groups 
\[{0\rightarrow A\rightarrow B \rightarrow C \rightarrow 0}\]
can be computed in $A$ and $C$ (Theorem \ref{ThmBdnExSeq}). This result is motivated by valued fields, as an $\RV$-structure can be seen as an enrichment of such. It is one of the main element for our proof of Theorem \ref{ThmBdnHenValFieCha00} and Theorem \ref{theoremmixedchar}.

    \subsection{Reduction}

As in the paragraph \ref{SubSubSecQEExaSeqAbeGrp}, we consider a pure exact sequence $\mathcal{M}$ of abelian groups
    \[ \xymatrix{0 \ar[r] & A \ar[r]^{\iota}& B \ar[r]^{\nu} & C\ar[r] &0}, \]
in an $\{A\}$-$\{C\}$-enriched language $\mathrm{L}$. 
In the following paragraph, we compute the burden of the structure $\mathcal{M}$ in terms of burden of $A$ and that of $C$ (in their induced structure). 
By bi-interpretability on unary sets, one can also consider it as a one-sorted structure $A\subset B$ where $A$ is given by a predicate. It follows indeed from Fact\ref{FactBdnInterpretUnarySet} that $\bdn(A \rightarrow B \rightarrow C) = \bdn(B,A)$. We often prefer the point of view of an exact sequence as it is more relevant for the computation of the burden. We write indifferently $\bdn(\mathcal{M})$, $\bdn(A\rightarrow B \rightarrow C)$ or $\bdn(B)$, as the sort $B$ is understood as a sort of $\mathcal{M}$ with its full induced structure.

Notice that in the case where $B/nB$  and $B_{[n]}$ are finite for all $n$ and $C$ is torsion free, a straight forward generalisation of \cite[Proposition 4.1]{CS19} gives that $\bdn(A \rightarrow B \rightarrow C)= \max(\bdn(A),\bdn(C))$. We will see that one can get rid of these hypothesis and obtain a more general result using Fact \ref{FactACGZ}. We first show a trivial bound:

    \begin{fact}[Trivial bound]\label{facttrivialbound}
   Assume there is a section of the group morphism $\nu: B \rightarrow C$. Consider $\mathrm{L}_{s}$ the language $\mathrm{L}$ augmented by a symbol $s$, and interpret it by this section of $\nu$.
    \[ \xymatrix{0 \ar[r] & A \ar[r]^{\iota}& B \ar[r]^{\nu} & C\ar[r]\ar@/_1pc/[l]_{s} &0}, \]
    
    We have $\bdn_{\mathrm{L}}(C) = \bdn_{\mathrm{L}_{s}}(C)$ and $\bdn_{\mathrm{L}}(A) = \bdn_{\mathrm{L}_{s}}(A)$ as well as the following:
        $$\max\{\bdn_{\mathrm{L}}(A),\bdn_{\mathrm{L}}(C)\} \leq \bdn_{\mathrm{L}}(B) \leq \bdn_{\mathrm{L}_{s}}(B) = \bdn_{\mathrm{L}}(A)+\bdn_{\mathrm{L}}(C). $$
    \end{fact}
    
    \begin{proof}
    The two first equalities are clear since $A$ and $C$ are stably embedded (and orthogonal) in both languages. The inequality $max\{\bdn_{\mathrm{L}}(A),\bdn_{\mathrm{L}}(C)\} \leq \bdn_{\mathrm{L}}(B)$ is obvious. As the burden only grows when we add structure, the inequality $\bdn_{\mathrm{L}}(B) \leq \bdn_{\mathrm{L}_{s}}(B)$ is also clear. The last equality come from the fact that in the language $\mathrm{L}_s$, the structure $A \rightarrow B \rightarrow C$ and the structure $ \{A\times C,A, C, \pi_A:A\times C \rightarrow A, \pi_C:A\times C \rightarrow C \}$ are bi-interpretable on unary sets. We conclude by Proposition \ref{sumburden} and Fact \ref{FactBdnInterpretUnarySet}.
    \end{proof}

    \begin{theorem}\label{ThmBdnExSeq}\index{Pure! short exact sequence of abelian groups}
        Consider an $\{A\}$-$\{C\}$-enrichement of a pure exact sequence $\mathcal{M}$ of abelian groups
    \[ \xymatrix{0 \ar[r] & A \ar[r]^{\iota}& B \ar[r]^{\nu} & C\ar[r] &0}, \]
in a language $\mathrm{L}$. Let $\mathcal{D}=\mathcal{D}(x)$ be the set of formulas in the pure language of groups which are conjunction of formulas of the form $\exists y\  nx=my$ for $n,m \in \mathbb{N}$.
For $D(x)\in \mathcal{D}$ and $A$ an abelian group, $D(A)$ is an subgroup of $A$, and we have 
        \[ \bdn\mathcal{M} =\max_{D\in \mathcal{D}}(\bdn(A/D(A))+\bdn(D(C))).\]
        
        In particular:
        \begin{itemize}
            \item If $A/nA$ is finite for all $n\geq 1$, then 
        \[\bdn{\mathcal{M}}= \max_{k\in \mathbb{N}}(\bdn(kA), \bdn(C_{[k]})),\]
        where $C_{[k]}:= \{c\in C \ \vert  \ kc=0\}$ is the subgroup of $k$-torsion.
            \item If $C$ has finite $k$-torsion of all $k\geq 1$, then 
        \[\bdn{\mathcal{M}}= \max_{n\in \mathbb{N}}(\bdn(A/nA) + \bdn(nC)).\] 
            \item If $C$ has finite $n$-torsion and $A/nA$ is finite for all $n\geq 1$, then 
            \[\bdn{\mathcal{M}}= \max(\bdn(A), \bdn(C)).\]
        
        \end{itemize}

    \end{theorem}
    
    To clarify, for $D(x)\in \mathcal{D}$, $\bdn(A/D(A))$ can be computed in $A/D(A)$ endowed with its induced structure, i.e. it is the supremum of depth of patterns $P(x_{D})$ with $x_{D}$ an $A/D(A)$ variable  within the structure $\{(A,0,+,\cdots), (A/D(A),0,+), \pi_{D}: A \rightarrow A/D(A) \}$, where $\cdots$ denotes the enriched structure in $A$. 
    
    \begin{remark}\label{rmkinfiniteburdencase}
        \begin{itemize}
            \item If $\bdn(A)$ or $\bdn(C)$ is infinite (or equals $\aleph_{0-}$), then this is simply the trivial bound in Fact \ref{facttrivialbound} (as then $\max(\bdn(A),\bdn(C))=\bdn(A)+\bdn(C)$). Recall that a section exists in a $\aleph_1$-saturated model as $A$ is pure in $B$ (Fact \ref{FactSectionPureSubgroupAleph1Saturated}).
            \item The maximum is always attained by at least one $D\in \mathcal{D}$: if $\bdn(A)$ and $\bdn(C)$ are finite, this is trivial. If $\bdn(A)$ or $\bdn(C)$ is infinite,  then $\bdn(A)$ or $\bdn(C)$ (corresponding terms for resp. $D \equiv x=0$ and $D \equiv x=x$ )  realises the maximum by the previous point. 
            \item By Fact $\ref{Ft1Bdn}$, if $C$ has finite $n$-torsion for $n\in \mathbb{N}$, one has $\bdn(nC)=\bdn(C)$. 
            \item If $m\vert n$, then $A/mA$ can be seen as a quotient of $A/nA$ and naturally, one has $\bdn(A/mA) \leq \bdn(A/nA)$. 
        \end{itemize}
    \end{remark}

    In the case that the sequence is unenriched, this can gives us absolute results: 
    assume that  the induced structures on $A/D(A)$ and $D(C)$ are the structures of groups for every $D\in \mathcal{D}$. Then, Proposition \cite[Theorem 1.1.]{HP17} together with Theorem \ref{ThmBdnExSeq} gives us a computation of $\bdn(A\rightarrow B\rightarrow C)$ in term of Szmielew invariants of $A$ and $C$. We don't attempt to write a closed formula. Nonetheless, here are some examples:
    
    \begin{examples}
        We consider the following pairs of abelian groups $(A\subset B)$ with quotient $C$:
        \begin{itemize}
            \item $B= \mathbb{Z}_{(2)}^{(\omega)} \oplus \mathbb{Z}_{(3)}^{(\omega)}$ and $A=\mathbb{Z}_{(2)}^{(\omega)} \oplus \{0\}$. One has $\bdn(A/nA)+\bdn(nC)=0+1=1$  for all $2 \nmid n $, and $\bdn(A/2nA)+\bdn(2nC)=1+1=2$, which leads to $\bdn{\mathcal{M}}=2$. This can already be deduced from Halevi and Palacín's work: the sort $B$, equipped only with its group structure, is already of burden $2$. By the trivial bound, the structure $\mathcal{M}$ is also of burden $2$.
            \item $B=\mathbb{Z}_{(2)}^{(\omega)} \oplus \mathbb{Z}(2)^{(\omega)}$ and $A=\mathbb{Z}_{(2)}^{(\omega)} \oplus 0$.
            Then $\bdn(\mathcal{M})=2$ (take $D(C)=C_{[2]}$ in Theorem \ref{ThmBdnExSeq}). Again, by \cite{HP17}, the burden of $B$ as a pure abelian group is already $2$.
            \item $B= \mathbb{Z}_{(2)}^{(\omega)} \oplus \mathbb{Z}_{(3)}^{(\omega)} \oplus \mathbb{Z}_{(2)}^{(\omega)} \oplus \mathbb{Z}_{(3)}^{(\omega)}$ and $A= \mathbb{Z}_{(2)}^{(\omega)} \oplus \mathbb{Z}_{(3)}^{(\omega)} \oplus \{0\}\oplus \{0\}$. Then $\bdn(\mathcal{M})=4$ (take $D(C)=6C$ in Theorem \ref{ThmBdnExSeq}). In term of subgroups, one can consider the subgroups $A+4B$, $A+9B$, $2B$ and $3B$. The intersection is 
            \[2\mathbb{Z}_{(2)}^{(\omega)} \oplus 3\mathbb{Z}_{(3)}^{(\omega)} \oplus 4\mathbb{Z}_{(2)}^{(\omega)} \oplus 9\mathbb{Z}_{(3)}^{(\omega)}.\]
            One may see that these groups satisfy Proposition \ref{PropCriBdnOneBasAbeGrp}.
            \item $A={\mathbb{Z}(2^\infty)} $, $C=\mathbb{Z}_{(2)}^{(\omega)} \oplus \mathbb{Z}_{(3)}^{(\omega)}$ and $B=A \times C$. 
            One can see that $\bdn(C\rightarrow B \rightarrow A)=\bdn(C,B)=3$. This equality is witnessed by the subgroups $2B,3B$ and $C+B_{[2]}$.
            However, by Theorem \ref{ThmBdnExSeq}, $\bdn\left( A \rightarrow B \rightarrow C\right)=\bdn(C)=2$ as $A/nA=\{0\}$ for all $n\geq 1$ and $C_{[k]}=0$ for all $k\geq 1$.
        \end{itemize}
    \end{examples}

    \begin{proof}[Proof of Theorem \ref{ThmBdnExSeq}]
        By Fact \ref{FactBdnInterpretUnarySet}, we can work in the language $\mathrm{L}_q$ and use Fact \ref{FactACGZ}. Recall that we abusively refer the union of sorts $\{A/nA\}_{n\in N}$ as the sort $A$. 
    
        The purity of the exact sequence 
            \[ \xymatrix{0 \ar[r] & A \ar[r]^{\iota}& B \ar[r]^{\nu} & C\ar[r] &0}, \]
        implies for all $n$ the following exact sequences:
        \[ \xymatrix{0 \ar[r] & A/nA \ar[r]^{\iota_n}& B/nB \ar[r]^{\nu_n} & C/nC\ar[r] &0}, \]
        \[ \xymatrix{0 \ar[r] & A_{[n]} \ar[r]^{\iota^{[n]}}& B_{[n]} \ar[r]^{\nu^{[n]}} & C_{[n]}\ar[r] &0},\]
        and 
        \[ \xymatrix{0 \ar[r] & A/A_{[n]} \ar[r]^{\iota^n}& B/B_{[n]} \ar[r]^{\nu^n} & C/C_{[n]}\ar[r] &0}. \]
        
        All these sequences are again pure, and we can keep going. We get more generally the following:
        
        \begin{fact}
            For all $D\in \mathcal{D}$ we have the following pure exact sequences:
            \[ \xymatrix{0 \ar[r] & D(A) \ar[r]^{\iota^D}& D(B) \ar[r]^{\nu^D} & D(C)\ar[r] &0}, \]
            and
            \[ \xymatrix{0 \ar[r] & A/D(A) \ar[r]^{\iota_D}& B/D(B) \ar[r]^{\nu_D} & C/D(C)\ar[r] &0}. \]
        \end{fact}
        \begin{proof}
            By Fact \ref{FactSectionPureSubgroupAleph1Saturated}, an $\aleph_1$-saturated extension $0\rightarrow A'\rightarrow B' \rightarrow C' \rightarrow 0$ splits. 
            Then we clearly have the exactness of the following sequence: 
                        \[ \xymatrix{0 \ar[r] & D(A') \ar[r]& D(B') \ar[r] & D(C')\ar[r] &0}, \]
            and
            \[ \xymatrix{0 \ar[r] & A'/D(A') \ar[r]& B'/D(B') \ar[r] & C'/D(C')\ar[r] &0}. \]
            As they are first order properties in the language $\mathcal{L}$, we deduce the fact.
        \end{proof}
        
        It follows that if $\nu(b)\in D(C)$, there is a unique $a+D(A)\in A/D(A)$ such that $i(a) + D(B) = b + D(B)$. We denote by $\rho_D$ the following map:
            \[\begin{array}{cccc}
            \rho_D  : & B & \to & A/D(A) \\
             & b & \mapsto & 
             \begin{cases}{\iota_D}^{-1}(b+D(B)) & \text{if } \nu(b)\in D(C) \\
              0 & \text{otherwise.}
             \end{cases}
        \end{array}\]
        Notice that it is interpretable in the language $\mathcal{L}$.

        By Lemma \ref{Ft1Bdn}, if $A$ is finite, we get that $\bdn(B) = \bdn(C)$ and of course $\bdn(A/D(A))=0$ for all $D\in \mathcal{D}$. Assume that $A$ is infinite. As $A$ and $C$ are orthogonal, so are in particular $A/D(A)$ and $D(C)$ for all $D\in \mathcal{D}$. 
        It follows by Fact \ref{sumburden} that $\bdn(A/D(A) \times D(C)) = \bdn(A/D(A))+\bdn(D(C))$. The interpretable (and surjective) map 
        
        \[\begin{array}{cccc}
            \rho_D \times \nu : & \nu^{-1}(D(C)) & \to & A/D(A)\times D(C) \\
             & b & \mapsto & (\rho_{D}(b),\nu(b))
        \end{array}\]
        gives us that
        \[\bdn(B) \geq \max_{D\in \mathcal{D}}(\bdn(A/D(A))+\bdn(D(C))).\]

        
        It remains to show that $\bdn B \leq \max_{D\in \mathcal{D}}(\bdn(A/D(A))+\bdn(D(C)))$. By Remark \ref{rmkinfiniteburdencase}, we may assume that $\bdn(A)$ and $\bdn(C)$ are both finite. As $A$ is infinite, $\bdn{A}\geq 1$. If $\bdn(\mathcal{M})=1$, the equality is clear. Assume that $\bdn(\mathcal{M})>1$ and let $P(x)=\lbrace \phi_i(x,y_i),(a_{i,j})_{j< \omega},k_i \rbrace_{i < M}$ be an inp-pattern of finite depth $M \geq 2$, with $(a_{i,j})_{i,j}$ mutually indiscernible and $\vert x\vert=1$. We need to show that $M\leq \bdn(A/D(A))+ \bdn(D(C))$ for some $D \in \mathcal{D}$. If $x$ is a variable in the sort $A$ (resp.\ in the sort $C$), $P(x)$ is an inp-pattern in $A$ (resp.\ $C$) of depth bounded by $\bdn(A)$ (resp.\ $\bdn(C)$) by purity (Corollary \ref{CoroPureOrth}). Then, the inequality holds if we take $D \equiv x=x$ (respect $D \equiv x=0$ ). 
        
        Assume $x$ is a variable in the sort $B$. Consider a line $\{\phi(x,y), (a_j)_{j<\omega}\}$ of $P(x)$ (we drop the index $i<M$ for the sake of clarity). By Fact \ref{FactACGZ}, and by the fact that one can "eliminate" disjunctions in inp-patterns (see \ref{INPdisj}), we may assume that the formula $\phi(x,a_j)$ is of the form 
        \begin{eqnarray} 
             & \phi_{A}(\rho_{n_0}(t^0(x,\beta_{j})),\ldots,\rho_{n_{s-1}}(t^{s-1}(x,\beta_{j})), \alpha_{j})\\
              & \wedge \phi_C(\nu({r^0}(x,\beta_{j}),\ldots,\nu({r^{k-1}}(x,\beta_{j})),\gamma_{j})\label{FormulaPhiCPart},
        \end{eqnarray}
        where $\phi_{A}$ is an $A$-formula, $\phi_C$ a $C$-formula, and for $j<\omega$, $\alpha_{j} \in A,\  \beta_{j} \in B ,\ \gamma_{j} \in C$ are parameters, $s,k,n_0,n_1,\ldots,n_{s-1}\in \mathbb{N}$ and the $t^l$'s and $r^l$'s are terms in the group language (one needs to keep in mind that $s$,$k$, $n_l$, $t^l$, $r^l$, $\beta_j$, $\alpha_j$ and $\gamma_j$ depend on the line $i$). Also, notice that $\rho_{n_0}(t^{0}(x,\beta_{j})) \neq 0 \in A/n_0A$ implies $ \nu(t^{0}(x,\beta_{j})) \in n_0C$ (a formula of the form (\ref{FormulaPhiCPart})).
    By writing the following
    \begin{align*}
        &\phi_{A}(\rho_{n_0}(t^0(x,\beta_{j})),\ldots,\rho_{n_{s-1}}(t^{s-1}(x,\beta_{j})), \alpha_{j})  \simeq  \\
        & \qquad \qquad \left(\phi_{A}(\rho_{n_0}(t^0(x,\beta_{j})),\ldots,\rho_{n_{s-1}}(t^{s-1}(x,\beta_{j})), \alpha_{j}) \wedge \nu(t^0(x,\beta_{j})) \in n_0C\right)\\
        & \qquad \qquad \bigvee \left(\phi_{A}(0,\rho_{n_1}(t^1(x,\beta_{j})),\ldots,\rho_{n_{s-1}}(t^{s-1}(x,\beta_{j})), \alpha_{j}) \wedge \nu(t^0(x,\beta_{j})) \notin n_0C \right),&
    \end{align*}

    and by eliminating once again the disjunction, one can assume that $\phi(x,a_j)$ (or more specifically, the formula $\phi_C(\nu({t^0}(x,\beta_{j}),\ldots,\nu({t^{s-1}}(x,\beta_{j})),\gamma_{j})$ ) implies $t^0(x,\beta_{j}) \in \nu^{-1}(n_0C)$. We do the same for all terms $t^l(x,\beta_{j})$, $l<s$. This means in particular that the list of terms $\{t^l\}_{l<s}$ is included in $\{r^l\}_{l<k}$. 
    
  Let $M'\geq 0$ be the number of rows such that \[\left\{\phi_C(\nu(r^0(x,\beta_{j}),\ldots,\nu(r^{k-1}(x,\beta_{j})),\gamma_{j}) \right\}_{j<\omega}\] is consistent. Without loss, they are the $M'$ first rows of the pattern $P(x)$, and we denote by $P'(x)$ the sub-pattern consisting of these rows. For now, we work with the sub-pattern $P'(x)$.
  
  Terms $t(x,\beta_{j})$ in the group language are of the form $kx+m\cdot \beta_{j}$, with $k\in \mathbb{N}, \ m\in \mathbb{N}^{\vert \beta_{j}\vert}$.
    \begin{claim}
        Assume that, in a line $\{\phi(x,y), (a_j)_{j<\omega}\}$ of $P'(x)$, a term $\rho_n(kx+m\cdot \beta_j)$ occurs. Then $\nu(m\cdot \beta_{j}) \mod nC$ is constant for all $j<\omega$.
    \end{claim}
    \begin{proof}
        Assume not. By indiscernibility, $\nu(m\cdot \beta_{j})$ are in distinct classes modulo $nC$. As 
        \[\phi_C(\nu({r^0}(x,\beta_{j}),\ldots,\nu({r^{s-1}}(x,\beta_{j})),\gamma_{j}) \vdash \nu(kx + m \cdot \beta_j)\in nC,\]
        the $\phi_C$ part of the line  
    \[\left\{\phi_C(\nu(r^0(x,\beta_{j}),\ldots,\nu(r^{s-1}(x,\beta_{j})),\gamma_{j})\right\}_{j<\omega}\]
        is $2$-inconsistent, contradicting the fact we chose one of the first $M'$ lines.  
    \end{proof}

    \begin{claim}[Main claim]
    We may assume that all formulas in $P'(x)$ are of the form
    \begin{eqnarray*} 
         & \phi_{A}(\rho_{n_0}(k^0(x-d),\ldots,\rho_{n_{s-1}}(k^{s-1}(x-d)), \alpha_{j})\\
         & \wedge \  \phi_C(\nu(x-d),\gamma_{j}).
    \end{eqnarray*}
    for certain integer $n_0,\cdots,n_{s-1}$,$k^0,\cdots,k^{s-1}$  and a certain parameter $d\in B$.
    \end{claim}
    \begin{proof}
    Take any realisation $d$ of the first column:
    \[d \models \{\phi_i(x,a_{i,0}) \}_{i<M}.\]
    
    We fix an $i<M'$, and consider the $i^{\text{th}}$ line $\{\phi(x,a_{j}) \}_{j<\omega}$, (again, we drop the index $i$ for a simpler notation). 
    
    \textbf{Step $1$:} We may assume that all terms $t^{l}(x, \beta_j)=k^lx+m^l\cdot \beta_j$, $l<s$ are of the form $k^l(x-d)$.\\
    We change all terms one by one, starting with $t^{0}(x, \beta_j)=k^0x+m^0\cdot \beta_j$. We write
    \[\rho_{n_0}(k^0x+m^0\cdot \beta_{j}) = \rho_{n_0}(k^0(x-d)+k^0d+m^0\cdot \beta_{j}).\]
    Replace it by $\rho_{n_0}(k^0(x-d))+\rho_{n_0}(k^0d+
    m^0\cdot \beta_{j})$. This doesn't change the formula $\phi(x,a_j)$ as $\phi(x,a_j) \vdash  \nu(k^0(x-d)) \in n_0C$. Indeed, $\phi(x,a_j) \vdash \nu(k^0x+m^0\cdot \beta_{j}) \in n_0C$ and $\nu(k^0d+m^0\cdot \beta_{j}) \in n_0C$ since $d$ is a solution of the first column and $\nu(m^0\cdot \beta_{j}) \mod n_0C$ is constant by the previous claim. Then, $\rho_{n_0}(k^0d-m^0\cdot \beta_{j})$ is seen as a parameter in $A/n_0A$, and it is added to $\alpha_j$. We do the same for all terms $t^l(x,\beta_{j})$, $l<s$.

    \textbf{Step $2$:} We may assume that all terms $r^l(x,\beta_j)=k^l x+m^l\cdot\beta_j$, $l<k$, are of the form $x-d$.\\
    This is immediate, as $\nu$ is a morphism. Indeed, replace $\nu(r^l(x,\beta_j))$ by $k^l\nu(x-d)+\nu(k^ld+m^l\cdot\beta_j)$, where $\nu(k^ld+m^l\cdot\beta_j)$ is seen as a parameter in $C$, and is added to the parameters $\gamma_j$. In other words, we may assume that the formula in the $i$th line of the pattern $P'(x)$ has the required form.
    
    Let us recall finally  that $d$ has been chosen independently of the row. We can apply these steps for all rows $i<M'$.
    \end{proof}

    
    Let $D$ be the conjunction of the formulas $\exists y \ kx=ny$ for each term $\rho_n(k(x-d))$ which occurs in \textbf{any} line of $P'(x)$.
     Notice that $\{\phi_i(x,a_{i,f(i)})\}_{i<M'}$ implies that $\nu(x-d) \in D(C)$ for any choice of function $f:M'\rightarrow \omega$. If $x_{A
    }$ is a variable in $A$, notice also that the truth value of the formula
     \[\phi_{A}(\pi_{n_0}(k^0x_A),\ldots,\pi_{n_{s-1}}(k^{s-1}x_A), \alpha_{j})\]
     evaluated in $a\in A$ is independent of the class of $a$ modulo $D(A)$. In other words, it interprets a definable set in $A/D(A)$.\\
    
    

     As now the pattern $P'(x)$ is centralised, one can remark that for every line \[\left\{\phi(x,a_j)\equiv \phi_{A}(\rho_{n_0}(k^0(x-d)),\ldots,\rho_{n_{s-1}}(k^{s-1}(x-d)), \alpha_{j}) \wedge   \phi_C(\nu(x-d),\gamma_{j})\right\}_{j<\omega},\]
     at most one of the following sets:
    \begin{equation}\label{PhiAPart}
    \tag{$L_A$}
        \{\phi_{A}(\pi_{n_0}(k^0x_A),\ldots,\pi_{n_{s-1}}(k^{s-1}x_A), y_{A}), (\alpha_{j})_{j <\omega}\}
    \end{equation}
    where $\vert y_{A}\vert =\vert\alpha_{j} \vert $ and $x_{A
    }$ is a variable in $A$
    or 
     \begin{equation}\label{PhiCPart}
     \tag{$L_C$}
     \lbrace \phi_{C}(x_C,y_{C}),(\gamma_{j})_{j< \omega} \rbrace
     \end{equation}
    (where $\vert y_{C}\vert =\vert\gamma_{j} \vert $) is consistent. 
    Indeed, this follows immediately from the fact that in the monster model, the sequence splits and $B \simeq A \times C$. 
    By definition of $P'(x)$, we deduce that (\ref{PhiAPart})-- the $\phi_A$-part of the line-- is inconsistent. Now, take a path $f:M'\rightarrow \omega$ and $b$ a solution of 
    \[\{ \phi_i(x,a_{i,f(i)})\}_{i<M'}\]
    As $b-d \in D(C)$, there is an $a$ such that $\iota(a)+D(B)=b-d+D(b)$. This give us that 
    \[a\models \{\phi_{i,A}(\pi_{n_0}(k^0x_A),\ldots,\pi_{n_{s-1}}(k^{s-1}x_A), \alpha_{i,f(i)})\}_{i<M'}.\]
    This show that every path of the following pattern $P'_A(x_A)$ is consistent: 
    \[ P'_A(x_A):=\{ \phi_{i,A}(\pi_{n_0}(k^0x_A),\ldots,\pi_{n_{s-1}}(k^{s-1}x_A), y_{A}), (\alpha_{i,j})_{j <\omega}\}_{i<M'}.\]
    Thus, this is an inp-pattern in $A$ and more precisely, it interprets an inp-pattern in $A/D(A)$ by the
    remark above. This means that $M'\leq \bdn(A/D(A)).$ 

    Then, the $\phi_C$-part (\ref{PhiCPart}) of any of the $M-M'$ last lines of $P(x)$ is inconsistent by definition of $P'(x)$ and $M'$. One gets as well an inp-pattern of depth $M-M'$ in $C$. As any realisation $r$ of $P'(x)$ (in particular, of $P(x)$) satisfies $\nu(r-d) \in D(C)$, one gets actually an inp-pattern in $D(C)$. It follows that $M-M' \leq \bdn(D(C))$. At the end, we get that $M \leq  \bdn(A/D(A))+\bdn(D(C))$.
    
    To conclude, we treat the particular cases. Assume that $C$ has finite $n$-torsion for every $n \in \mathbb{N}^{\star}$. Then infinite subgroups in $\mathcal{D}(C):= \{D(C) \ \vert \ D\in \mathcal{D}\}$ are of the form $nC$. We deduce that 
            \[\bdn{\mathcal{M}}= \max_{n\in \mathbb{N}}(\bdn(A/nA) + \bdn(nC)).\] 
    
    Similarly assume that $A/nA$ is finite for all $n \in \mathbb{N}^{\star}$. 
    Then $\bdn(A/nA_{[k]})$ is equal to $\bdn(A/A_{[k]})$.
    Indeed, this can be deduced from Fact \ref{Ft1Bdn}, as we have the exact sequence
    \[0\rightarrow A_{[k]}/nA_{[k]} \rightarrow A/nA_{[k]} \rightarrow A/A_{[k]}\rightarrow 0,\]
    and as $A_{[k]}/nA_{[k]}$ is finite. Since $\bdn(C_{[k]})\geq \bdn(nC_{[k]})$ and $A/A_{[k]}\simeq kA$, we may deduce that
            \[\bdn{\mathcal{M}}= \max_{k\in \mathbb{N}}(\bdn(kA), \bdn(C_{[k]})).\]

    \end{proof}

    \subsection{Application}
    As main application of Theorem \ref{ThmBdnExSeq}, we will deduce Theorem \ref{ThmBdnHenValFieCha00}. This is the aim of the next section. For now, we want to emphasise the advantage of working resplendently by giving one straightforward generalisation of Theorem \ref{ThmBdnExSeq}.

    \begin{corollary}
    Let $\mathcal{M}$ be an exact sequence of \emph{ordered} abelian groups
    \[ \xymatrix{0 \ar[r] & A \ar[r]^{\iota}& B \ar[r]^{\nu} & C\ar[r] &0}, \]
    where $(A,<)$ is a convex subgroup of $(C,<)$. We consider it as a three sorted structure, with a structure of ordered abelian group for each sort, and function symbols for $\iota$ and $\nu$. Then, we have:  
    \begin{align*}
     \bdn\mathcal{M} &=\max_{n\in \mathbb{N}}(\bdn(A/nA)+\bdn(nC))\\
       &= \max \left(\bdn(A) , \max_{n\in \mathbb{N}^\star}(\bdn(A/nA))+\bdn(C) \right). 
     \end{align*}
\end{corollary}
    \begin{proof}
        As $C$ is torsion free, $\iota(A)$ is pure in $B$. As for $b\in B$, $b>0$ if and only if $\nu(b)>0$ or $\nu(b)=0$ and $\iota^{-1}(b)>0$, $\mathcal{M}$ is an $\{A\}$-$\{C\}$-enrichment of a short exact sequence of abelian groups. It remains to apply Theorem \ref{ThmBdnExSeq}. Notice that for all $n>0$, we have $\bdn(nC)=\bdn(C)$ (as the multiplication by $n$ in $C$ is a definable injective morphism).
    \end{proof}

\section{Burden of Henselian valued fields} \label{SectionBurden}
    We compute the burden of benign Henselian valued fields and of unramified mixed characteristic Henselian valued fields with perfect residue field in terms of the burden of the  value group and that of the residue field. The first subsection is common for both cases and treats the reduction from the valued field to the sort $\RV$  (resp.\ the sorts $\RV_{<\omega}$). 
    For the reduction to the value group and residue field, we treat (separately) the case of benign Henselian valued fields in Subsection \ref{SectionHenselianValuedFieldsEquichar0} and the case of unramified mixed characteristic henselian valued fields with perfect residue field in Subsection \ref{SectionUnrimified}.
    They are both deduced from the computation of burden in short exact sequences of abelian groups (Section \ref{SectionExactSequence}). 
    
	\subsection{Reductions to $\RV$ and $\RV_{<\omega}$}	\label{SectionReductionBurdenRV}
	We compute here the burden of Henselian valued fields of characteristc $0$ in terms of burden of  $\RV_{<\omega}$. As we explained in Paragraph \ref{Preliminaries Unramified mixed characteristic Henselian valued fields}, mixed characteristic Henselian valued fields satisfy $(\text{EQ})_{\RV_{<\omega}}$ (quantifier elimination relative to the union of sorts $\RV_{<\omega}$.), but do not satisfy in general $(\text{EQ})_{\RV}$ (elimination of quantifiers  relative to $\RV$). Our result includes naturally the case of unramified mixed Henselian valued fields, and also equicharacteristic $0$ Henselian valued fields. In the former case, we have a computation of the burden in term of the burden of $\RV$, as in equicharacteristic $0$ the structures $\RV$ and $\RV_{<\omega}$ can be identified (Remark \ref{RemarkIdentifyRVRVOmega}). In fact, the proof that we are going to present can be adapted for all benign valued fields.

        \subsubsection{Reductions}
The aim of this paragraph is to prove the following:

\begin{theorem}\label{ThmHensValuedFieldReductionRVMixedCharacteristic}\index{$\RV_{<\omega}$-sort}
 Let $K$ be a Henselian valued field of characteristic $(0,p)$, $p\geq 0$. Let $M$ be a positive integer and assume $\mathcal{K}$ is of burden $M$. Then, the sort $\RV_{<\omega}$ with the induced structure is also of burden $M$. In particular, $\mathcal{K}$ is inp-minimal if and only if $\RV_{<\omega}$ is inp-minimal.
\end{theorem} 
The demonstration (below) follows Chernikov and Simon's proof for the case of equicharacteristic $0$ and burden $1$ (see \cite{CS19}).
As we said, this statement also cover the case of equicharacteristic $0$.  One can also generalise the proof for infinite burden (see Corollary \ref{corollaryinfiniteburden} for details). A careful reading of the proof shows that one only uses properties $(\text{EQ})_{\RV_{<\omega}}$ and $(\text{Lin})_{\RV_{<\omega}}$. Of course, the proof can be written for equicaracteristic $0$ fields only (it becomes  simpler), and then it only uses Property $(\text{EQ})_{\RV}$ and $(\text{Lin})_{\RV}$. As algebraically maximal Kaplansky valued fields and algebraically closed valued fields satisfy these property, we obtain in fact:

\begin{theorem}\label{ThmHensValuedFieldReductionRV} \index{Valued field! benign Henselian} \index{Inp-minimality} \index{$\RV$-sort}
Let $\mathcal{K}$ be a benign Henselian valued field. Let $M$ be a positive integer and assume $\mathcal{K}$ is of burden $M$. Then, the sort $\RV$ with the induced structure is also of burden $M$. In particular, $\mathcal{K}$ is inp-minimal if and only if $\RV$ is inp-minimal.
\end{theorem}

\begin{proof}[Proof of Theorem \ref{ThmHensValuedFieldReductionRVMixedCharacteristic}]

We denote by $\bar{\mathbb{Z}}$ the set of natural numbers with extremal points $ \mathbb{Z}\cup \lbrace \pm \infty\rbrace$. Let $\lbrace \tilde{\phi}_i(x,y_i),(c_{i,j})_{j \in \bar{\mathbb{Z}}},k_i \rbrace_{i < M}$ be an inp-pattern in $\mathcal{K}$ of finite depth $M \geq 2$ with $\vert x \vert = 1$, where $c_{i,j}=a_{i,j}\mathbf{b}_{i,j} \in K^{k_1} \times \RV_{<\omega}^{k_2}$. Notice that the set of indices is $\bar{\mathbb{Z}}$, as we will make use of one of the extreme elements $\lbrace a_{i,-\infty},a_{i,+\infty}\rbrace$) later. We have to find an inp-pattern of depth $M$ in $\RV_{<\omega}$. Without loss of generality, we take $(c_{i,j})_{i,j}$ mutually indiscernible. By Fact \ref{factfle} and mutual indiscernibility, we can assume the formulas $\tilde{\phi}_i$ are of the form 
$$\tilde{\phi}_i(x,c_{i,j}) = \phi_i(\rv_{\delta_n}(x-a_{i,j;1}), \ldots , \rv_{\delta_n}(x-a_{i,j;k_1});\mathbf{b}_{i,j}),$$
for some integer $n$ and where $\phi_i$  are $\RV_{<\omega}$-formulas. Also recall that $\delta_n$ denotes the value $\val(p^n)$. The arguments inside symbols $\rv_{\delta_n}$ are linear terms in $x$. In some sense, difficulties coming from the field structure have been already treated and it only remains to deal with the structure coming from the valuation.

Let $d \models \lbrace \phi_i(\rv_{\delta_n}(x-a_{i,0;1}), \ldots , \rv_{\delta_n}(x-a_{i,0;k_1});\mathbf{b}_{i,0}) \rbrace_{i <M}$ be a solution of the first column. Before we give a general idea of the proof, let us reduce to the case where only one term $\rv_{\delta_n}(x-a_{i,j})$ occurs in the formula $\tilde{\phi}_i$.

\begin{claim}\label{claim1}
We may assume that for all $i<M$, $\tilde{\phi}_i(x,c_{i,j})$ is of the form $\phi_i(\rv_{\delta_n}(x-a_{i,j;1});\mathbf{b}_{i,j})$, i.e. $\vert a_{i,j} \vert = k_1=1$.
\end{claim}
\begin{proof}
We will first replace the formula $\tilde{\phi}_0(x,c_{0,j})$ by a new one with an extra parameter. 

By Lemma \ref{WDcase}, at least one of the following two cases occurs
\begin{enumerate}

    \item $\WD_{\delta_n} \left(\rv_{\delta_n}(d-a_{0,0;1}), \rv_{\delta_n}(a_{0,0;1}-a_{0,0;2}) \right)$ or
    \item $\WD_{\delta_n} \left(\rv_{\delta_n}(d-a_{0,0;2}),\rv_{\delta_n}(a_{0,0;2}-a_{0,0;1}) \right).$
\end{enumerate}
 
According to the case, we respectively define a new formula $\psi_0(x,c_{0,j}\upwedge \rv_{\delta_n}(a_{0,j;2}-a_{0,j;1}) )$  by:
\begin{enumerate}
    \item \begin{multline*}\phi_0(\rv_{\delta_n}(x-a_{0,j;1}),\rv_{\delta_n}(x-a_{0,j;1})+\rv_{\delta_n}(a_{0,j;1}-a_{0,j;2}),\rv_{\delta_n}(x-a_{0,j;3}), \ldots , \\ \rv_{\delta_n}(x-a_{0,j;k});\mathbf{b}_{0,j}) \wedge \WD_{\delta_n}(\rv_{\delta_n}(x-a_{0,j;1}), \rv_{\delta_n}(a_{0,j;1}-a_{0,j;2})),
\end{multline*}

  \item \begin{multline*}\phi_0(\rv_{\delta_n}(x-a_{0,j;2})+\rv_{\delta_n}(a_{0,j;2}-a_{0,j;1}),\rv_{\delta_n}(x-a_{0,j;2}), \ldots , \rv_{\delta_n}(x-a_{0,j;k});\mathbf{b}_{0,j}) \\ \wedge \WD_{\delta_n}(\rv_{\delta_n}(x-a_{0,j;2}), \rv_{\delta_n}(a_{0,j;2}-a_{0,j;1})).
\end{multline*}
\end{enumerate}
We will prove that the pattern where $\tilde{\phi}_0$ is replaced by $\psi_0$:
$$\lbrace \psi_0(x,y_0 \upwedge z),(c_{0,j}\upwedge \rv_{\delta_n}(a_{0,j;2}-a_{0,j;1}))_{j \in \bar{\mathbb{Z}}}, k_0 \rbrace \cup \lbrace \tilde{\phi}_i(x,y_i),(c_{i,j})_{j \in \bar{\mathbb{Z}}},k_i \rbrace_{1 \leq i <M}$$
is also an inp-pattern. First note that we have added $\rv_{\delta_n}(a_{0,j;2}-a_{0,j;1})$ to the parameters $\mathbf{b}_{0,j}$, and it still forms a mutually indiscernible array. Clearly, $d$ is still a realisation of the first column: 
$$d\models \lbrace \psi_0(x,c_{0,0}\upwedge \rv_{\delta_n}(a_{0,0;2}-a_{0,0;1}) ) \rbrace \cup \lbrace \tilde{\phi}_i(x,c_{i,0})\ \vert \ 1 \leq i <M  \rbrace.$$ By mutual indiscernibility of the parameters, every path is consistent. Since $\psi_0(\mathcal{K}) \subseteq \tilde{\phi}_0(\mathcal{K})$, inconsistency of the first row is also clear. By induction, it is clear that we may assume that $\phi_0$ is of the desired form. We can do the same for all formulas $\phi_i$, $0 < i < M$.
\end{proof}

If the array $(a_{i,j})_{i<M, j<\omega}$ is constant equal to some $a\in K$, then we obviously get an inp-pattern of depth $M$ in $\RV_{<\omega}$:
 $\lbrace \phi_i(\mathbf{x},z_i),(\mathbf{b}_{i,j})_{j \in \bar{\mathbb{Z}}},k_i \rbrace_{i <M}$, where $\mathbf{x}$ is a variable in $\RV_{\delta_n}$ (such a pattern is said to be \textit{centralised} \index{Inp-pattern! centralised}). Indeed, consistency of the path is clear. If a row is satisfied by some $\mathbf{d}\in \RV_{\delta_n}$  , any $d\in K$ such that $\rv_{\delta_n}(d-a)=\mathbf{d}$ will satisfy the corresponding row of the initial inp-pattern, which is absurd. Hence, the rows are inconsistent.  

The idea of the proof is to reduce the general case (where the $a_{i,j}$'s are distinct) to this trivial case by the same method as above: removing the parameters $a_{i,j} \in K$ and adding new parameters from $\RV_{<\omega}$ to $\mathbf{b}_{i,j}$ and specifying the formula by adding a term of the form $\WD(\rv(x-a),\rv(a-a_{i,j}))$. The main challenge is to find a suitable $a\in K$ for a center.

Recall that $d \models \lbrace \phi_i(\rv_{\delta_n}(x-a_{i,0});b_{i,0}) \rbrace_{i <M}$ is a solution to the first column.

\begin{claim}\label{claim2}
For all $j < \omega $, and $i,k < M$ with $k\neq i$, we have ${\val(d-a_{i,j}) \leq \val(d-a_{k,0}) + \delta_n}. $
\end{claim}
 \begin{proof}
 	Assume not: for some $j < \omega$, and $i,k < M$ with $k\neq i$:
 	$$\val(d-a_{i,j}) > \val(d-a_{k,0}) + \delta_n.$$
Then, $\rv_{\delta_n}(a_{i,j}-a_{k,0}) = \rv_{\delta_n}(d-a_{k,0})$. By mutual indiscernibility, we have
$$a_{i,j} \models \lbrace \phi_k(\rv_{\delta_n}(x-a_{k,l});\mathbf{b}_{k,l}) \rbrace_{l< \omega}.$$
This contradicts inconsistency of the row $k$.
\end{proof}

In particular, for all $i,k < M$, we have $ {\vert \val(d-a_{k,0})-\val(d-a_{i,0}) \vert \leq \delta_n }$. For $i<M$, let us denote $ {\gamma_i:=\val(d-a_{i,0}) }$ and let $\gamma$ be the minimum of the $\gamma_i$'s. By definition, we have the following for all $i,k<M$:
\begin{equation}
 \ \val(a_{i,0}-a_{k,0}) \geq \min \lbrace \val(d-a_{i,0}), \val(d-a_{k,0}) \rbrace \geq \gamma.\tag{$\star$}
\end{equation}

The following claim give us a correct centre $a$.

\begin{claim}\label{mainprop}
We may assume that there is $i <M $ such that for all $k < M$, the following holds:

\begin{eqnarray*}
\gamma_k = val(d-a_{k,0}) &\leq&\min \lbrace \val(d-a_{i,\infty}), \val(a_{i,\infty}-a_{k,0}) \rbrace + \delta_n. \\
\end{eqnarray*}
In particular, by Proposition \ref{WDdelta}, we have:
\begin{eqnarray*}
&&\WD_{\delta_n}(\rv_{2\delta_n}(d-a_{i,\infty}),\rv_{2\delta_n}(a_{i,\infty}-a_{k,0})).
 \end{eqnarray*}
\end{claim}

\begin{proof}
By Remark \ref{remark}, it is enough to find $i<M$ such that the following holds for all $k < M$:
$$\gamma_k  \leq \val(d-a_{i,\infty}) + \delta_n \quad \text{or} \quad \gamma_k \leq   \val(a_{i,\infty}-a_{k,0}) + \delta_n. $$ 

We will actually find $i$ such that one of the following holds:

\begin{enumerate}
    \item $\gamma_k  \leq \val(d-a_{i,\infty}) + \delta_n \quad  \text{for all $k < M$}$ 
    \item $\gamma_k \leq   \val(a_{i,\infty}-a_{k,0}) + \delta_n \quad  \text{for all $k < M$}$
\end{enumerate}
The first case will correspond to Case A, the second to Case B.  \\

	\textbf{Case A }: There are $0 \leq i, k < M$ with $i \neq k$ such that $\val(a_{i,j}-a_{k,l})$ is constant for all $j,l \in \omega$, equal to some $\epsilon$.
	Note that $(\star)$ gives $\epsilon \geq \gamma$.
	\\
	
\begin{center}
\begin{minipage}{0.70\linewidth}

	\begin{tikzpicture}[line cap=round,line join=round,>=triangle 45,x=1.0cm,y=1.0cm,scale = 0.5]
\clip(-2.226855653004253,-4.195898595071177) rectangle (17.415560940241033,5.607988346530954);
\draw [line width=1.2pt] (2.,5.)-- (6.,-3.);
\draw [line width=1.2pt] (6.,-3.)-- (10.,5.);
\draw [line width=1.2pt] (6.,-3.)-- (7.,-5.);
\draw [line width=1.2pt] (2.2477554844708085,4.504489031058383)-- (2.5,5.);
\draw [line width=1.2pt] (2.5,4.)-- (3.,5.);
\draw [line width=1.2pt] (2.741734673539459,3.5165306529210825)-- (3.5,5.);
\draw [line width=1.2pt] (3.498797557358775,2.0024048852824508)-- (5.,5.);
\draw [line width=1.2pt] (3.748010012984067,1.5039799740318665)-- (5.5,5.);
\draw [line width=1.2pt] (6.,5.)-- (4.,1.);
\draw [line width=1.2pt] (6.5,5.)-- (4.244320737961419,0.509786121281313);
\draw [line width=1.2pt] (9.,3.)-- (8.,5.);
\draw [line width=1.2pt] (8.244978810717003,4.510042378565994)-- (8.5,5.);
\draw [line width=1.2pt] (8.5,4.)-- (9.,5.);
\draw [line width=1.2pt] (8.748951172574003,3.502097654851993)-- (9.5,5.);
\draw [line width=1.2pt] (3.,3.)-- (4.,5.);
\draw [line width=1.2pt] (7.,5.)-- (4.5,0.);
\draw [line width=0.8pt,dotted] (3.,-3.)-- (10.,-3.);
\draw (7.9061370339873625,-2.966082388615786) node[anchor=north west] {$\val(a_{i,j}-a_{k,l})= \epsilon \geq \gamma $};
\begin{scriptsize}
\draw [fill=black] (2.,5.) circle (1.5pt);
\draw [fill=black] (6.,-3.) circle (1.0pt);
\draw [fill=black] (10.,5.) circle (1.5pt);
\draw[color=black] (10.521661642082634,5.174954471018492) node {$a_{k,\infty}$};
\draw [fill=black] (7.,-5.) circle (1.0pt);
\draw [fill=black] (2.2477554844708085,4.504489031058383) circle (1.0pt);
\draw [fill=black] (2.5,5.) circle (1.5pt);
\draw[color=black] (2.579820365184078,5.30959718105949) node {$a_{i,0}$};
\draw [fill=black] (2.5,4.) circle (1.0pt);
\draw [fill=black] (3.,5.) circle (1.5pt);
\draw [fill=black] (2.741734673539459,3.5165306529210825) circle (1.0pt);
\draw [fill=black] (3.5,5.) circle (1.5pt);
\draw [fill=black] (3.498797557358775,2.0024048852824508) circle (1.0pt);
\draw [fill=black] (5.,5.) circle (1.5pt);
\draw [fill=black] (3.748010012984067,1.5039799740318665) circle (1.0pt);
\draw [fill=black] (5.5,5.) circle (1.5pt);
\draw [fill=black] (6.,5.) circle (1.5pt);
\draw [fill=black] (4.,1.) circle (1.0pt);
\draw [fill=black] (6.5,5.) circle (1.5pt);
\draw [fill=black] (4.244320737961419,0.509786121281313) circle (1.0pt);
\draw [fill=black] (9.,3.) circle (1.0pt);
\draw [fill=black] (8.,5.) circle (1.5pt);
\draw [fill=black] (8.244978810717003,4.510042378565994) circle (1.0pt);
\draw [fill=black] (8.5,5.) circle (1.5pt);
\draw[color=black] (8.746222752481538,5.20959718105949) node {$a_{k,0}$};
\draw [fill=black] (8.5,4.) circle (1.0pt);
\draw [fill=black] (9.,5.) circle (1.5pt);
\draw [fill=black] (8.748951172574003,3.502097654851993) circle (1.0pt);
\draw [fill=black] (9.5,5.) circle (1.5pt);
\draw [fill=black] (3.,3.) circle (1.0pt);
\draw [fill=black] (4.,5.) circle (1.5pt);
\draw[color=black] (4.259187384908164,5.30959718105949) node {$a_{i,\infty}$};
\draw [fill=black] (7.,5.) circle (1.5pt);
\draw [fill=black] (4.5,0.) circle (1.0pt);
\end{scriptsize}
\end{tikzpicture}

\end{minipage}
\end{center}
	
	Then, we have:
	$$\val(d-a_{i,\infty})\geq \min\lbrace \val(d-a_{i,0}),\val(a_{i,0}-a_{i,\infty})\rbrace \geq \gamma$$
	Indeed, $\val(a_{i,0}-a_{i,\infty}) \geq \min\lbrace \val(a_{i,0}-a_{k,0}),\val(a_{k,0}-a_{i,\infty})\rbrace = \epsilon \geq \gamma$. 
	Hence, we have for every $0 \leq l < M$:
	$$\val(d-a_{i,\infty}) + \delta_n \geq \gamma + \delta_n \geq \val(d-a_{l,0}) = \gamma_l . $$
	
	\textbf{Case B}: For all $0\leq i,  k < M$ with $i \neq k$, $\left(\val(a_{i,j}-a_{k,l})\right)_{j,l}$ is not constant. \\
	By Lemma \ref{lemmaind1} (2) and Lemma \ref{lemmaind2}, there is $i < M$ such that for every $k< M$ and $k\neq i$,  $(a_{k,l})_{l < \omega} {\Rightarrow} a_{i,0}$  or $(a_{k,-l})_{l < \omega} {\Rightarrow} a_{i,0}$. If needed, one can flip the indices and assume that for all $k\neq i$,  $(a_{k,l})_{l < \omega} {\Rightarrow} a_{i,0}$. Note that only $(a_{i,j})_j$ could be a fan in this case.
	
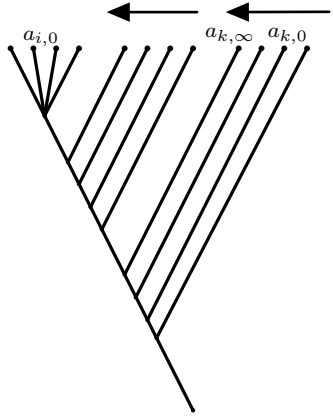
\begin{center}
\begin{minipage}{0.70\linewidth}

	\begin{tikzpicture}[line cap=round,line join=round,>=triangle 45,x=1.0cm,y=1.0cm,scale=0.6]
\clip(-2.6813490653752179,-3.244744554454015) rectangle (16.301104990983102,6.8306249304654555);
\draw [line width=1.2pt] (6.,-3.)-- (2.,5.);
\draw [line width=1.2pt] (2.745113232048317,3.509773535903366)-- (3.5,5.);
\draw [line width=1.2pt] (3.2572757384796804,2.485448523040639)-- (4.5,5.);
\draw [line width=1.2pt] (3.5,2.)-- (5.,5.);
\draw [line width=1.2pt] (3.7524098570521747,1.4951802858956507)-- (5.5,5.);
\draw [line width=1.2pt] (4.,1.)-- (6.,5.);
\draw [line width=1.2pt] (7.,5.)-- (4.5,0.);
\draw [line width=1.2pt] (7.5,5.)-- (4.74424472152211,-0.5033601209440769);
\draw [line width=1.2pt] (8.,5.)-- (5.,-1.);
\draw [line width=1.2pt] (8.5,5.)-- (5.2,-1.4);
\draw [line width=1.2pt] (2.5,5.)-- (2.745113232048317,3.509773535903366);
\draw [line width=1.2pt] (3.,5.)-- (2.745113232048317,3.509773535903366);
\draw [->,line width=1pt] (8.983574255409578,5.808506646116047) -- (6.715348226133187,5.801655878143558);
\draw [->,line width=1pt] (6.089660595481235,5.801655878143558) -- (4.1045837982323555,5.801655878143558);
\begin{scriptsize}
\draw [fill=black] (6.,-3.) circle (1.0pt);
\draw [fill=black] (2.,5.) circle (1.5pt);
\draw [fill=black] (2.5,5.) circle (1.5pt);
\draw[color=black] (2.671546559378185,5.2231835461590554) node {$a_{i,0}$};
\draw [fill=black] (3.,5.) circle (1.5pt);
\draw [fill=black] (2.745113232048317,3.509773535903366) circle (1.0pt);
\draw [fill=black] (3.5,5.) circle (1.5pt);
\draw [fill=black] (3.2572757384796804,2.485448523040639) circle (1.0pt);
\draw [fill=black] (4.5,5.) circle (1.5pt);
\draw [fill=black] (3.5,2.) circle (1.0pt);
\draw [fill=black] (5.,5.) circle (1.5pt);
\draw [fill=black] (3.7524098570521747,1.4951802858956507) circle (1.0pt);
\draw [fill=black] (5.5,5.) circle (1.5pt);
\draw [fill=black] (4.,1.) circle (1.0pt);
\draw [fill=black] (6.,5.) circle (1.5pt);
\draw [fill=black] (7.,5.) circle (1.5pt);
\draw[color=black] (6.81473270324649,5.27561091638033) node {$a_{k,\infty}$};
\draw [fill=black] (4.5,0.) circle (1.0pt);
\draw [fill=black] (7.5,5.) circle (1.5pt);
\draw [fill=black] (4.74424472152211,-0.5033601209440769) circle (1.0pt);
\draw [fill=black] (8.,5.) circle (1.5pt);
\draw[color=black] (8.07896881089334,5.27561091638033) node {$a_{k,0}$};
\draw [fill=black] (5.,-1.) circle (1.0pt);
\draw [fill=black] (8.5,5.) circle (1.5pt);
\draw [fill=black] (5.2,-1.4) circle (1.0pt);
\end{scriptsize}
\end{tikzpicture}

\end{minipage}
\end{center}

	Then we have \[\val(a_{k,0}-a_{i,\infty})=\val(a_{k,0}-a_{i,0}) \geq \gamma,\]
	since $(a_{k,l}) {\Rightarrow} a_{i,\infty}$ as well.
So 
\[\val(a_{k,0}-a_{i,\infty}) + \delta_n \geq \gamma_k,\]
It remains to prove the inequality for $k=i$. Take $l \neq i$, $l< M$. We have:
$$\val(a_{i,0}-a_{i,\infty}) \geq \min\lbrace \val(a_{i,0}-a_{l,0}),\val(a_{l,0}-a_{i,\infty}) \rbrace  \geq \gamma.$$
Hence, $\val(a_{i,0}-a_{i,\infty}) + \delta_n \geq \gamma_i .$

\end{proof}

Assume $i=0$ satisfies the conclusion of the previous claim. For every $k< M$, we have  \[\WD_{\delta_n}(\rv_{2\delta_n}(d-a_{0,\infty}),\rv_{2\delta_n}(a_{0,\infty}-a_{k,0})).
\]
Set $\tilde{b}_{i,j}:=b_{i,j} \upwedge \rv_{2\delta_n}(a_{0,\infty}-a_{i,j}) $ for $i < M, j<\omega$ and 
$$\psi_i(\tilde{x},\tilde{b}_{i,j}):=\phi_i\left(\rv_{\delta_n}(\tilde{x}+\rv_{2\delta_n}(a_{0,\infty}-a_{i,j}));b_{i,j}\right) \wedge \WD_{\delta_n}(\tilde{x},\rv_{2\delta_n}(a_{0,\infty}-a_{i,j}))$$
where $\tilde{x}$ is a variable in $\RV_{2\delta_n}$.

This is an inp-pattern. Indeed, clearly, $\rv_{2\delta_n}(d-a_{0,\infty}) \models \lbrace \psi_i(\tilde{x},\tilde{b}_{i,0})\rbrace_{i < M}$. By mutual indiscernibility of $(\tilde{b}_{i,j})_{i<M,j < \omega}$, every path is consistent. It remains to show that, for every $i < M$, $\lbrace \psi_i(\tilde{x},\tilde{b}_{i,j}) \rbrace_{j<\omega}$ is inconsistent.
Assume there is $\alpha^\star  \models \lbrace \psi_i(\tilde{x},\tilde{b}_{i,j}) \rbrace_{j<\omega}$ for some $i< M$, and let $d^\star$ be such that $\rv_{2\delta_n}(d^\star - a_{0,\infty})= \alpha^\star$. Then, since $\WD_{\delta_n}(\alpha^\star,\rv_{2\delta_n}(a_{0,\infty}-a_{i,j}))$ holds for every $j<\omega$, $d^\star$ satisfies ${\lbrace \phi_i(\rv_{\delta_n}(x-a_{i,j}), b_{i,j})\rbrace_{j<\omega}}$, which is a contradiction. All rows are inconsistent, which concludes our proof.

\end{proof}

With minor modifications, the proof goes through in the case of infinite burden $\lambda$. However, one must be careful regarding the precise statement of this generalisation. Assume we are in mixed characteristic $(0,p)$, and the burden $\lambda$ is of cofinality $\cf(\lambda)=\omega$. Then the very first argument of the proof is no longer true: one cannot necessary assume that there are $\lambda$-many formulas $\tilde{\phi}_i(x,y_i)$ in the inp-pattern of the form
$$\tilde{\phi}_i(x,c_{i,j}) = \phi_i(\rv_{\delta_n}(x-a_{i,j;1}), \ldots , \rv_{\delta_n}(x-a_{i,j;k});\mathbf{b}_{i,j}),$$
for a certain $n<\omega$. This depends of course of the cofinality of $\lambda$. Nonetheless, this is the only problem. One gets the following statement:

\begin{corollary}\label{corollaryinfiniteburden}
\begin{itemize}
Let $\lambda$ be an infinite cardinal in $\Card^\star$.
    \item Let $K$ be a mixed characteristic Henselian valued field. Assume that the union of sorts $\RV_{<\omega}$ with the induced structure is of burden $\lambda$. Then, the field $K$ is of burden $\lambda$ if $\cf(\lambda)> \omega$, and of burden $\lambda$ or $\act(\lambda)$ if $\cf(\lambda)=\omega$.
    \item Let $K$ be a benign Henselian valued field. Assume that the sort $\RV$ with the induced structure is of burden $\lambda$. Then the field $K$ is of burden $\lambda$.

\end{itemize}

\end{corollary}

\begin{proof}
        We treat the case of mixed characteristic Henselian valued field. We prove similarly the case of benign Henselian valued fields.
        Let $\kappa\geq \lambda$ be the burden of $K$, and let 
        \[ \lbrace \tilde{\phi}_i(x,y_i),(c_{i,j})_{j<\omega}\rbrace_{i<\kappa} \] be an inp-pattern of depth $\kappa$.
        If $\kappa$ is of cofinality $\cf(\kappa)>\omega$, then there are $\kappa$-many formulas $\tilde{\phi}_i(x,y_i)$ in the inp-pattern of the form
$$\tilde{\phi}_i(x,c_{i,j}) = \phi_i(\rv_{\delta_n}(x-a_{i,j;1}), \ldots , \rv_{\delta_n}(x-a_{i,j;k});\mathbf{b}_{i,j}),$$
for a certain $n<\omega$. We deduce an inp-pattern of depth $\kappa$ in the $\RV_{<\omega}$-sort. Indeed, we follow the exact same proof with few changes in Claim \ref{mainprop}:
    
    \begin{itemize}
        
    \item The minimum of $\{\gamma_k \}_{k<\lambda}$ may not exist, but one can pick $\gamma$ in an extension of the monster model, realising the cut $\{\gamma  \in\Gamma \ \vert \ \gamma<\gamma_k \ \text{for all }k<\lambda\} \cup \{\gamma  \in\Gamma \ \vert \ \gamma >\gamma_k \ \text{for some }k<\lambda\}$. By Claim \ref{claim2}, we have for all $a\in K$, $\val(a)> \gamma$ implies $\val(a)+\delta_n > \gamma_k$ for all $k<\lambda$. 
    \item Case A stays the same. 
    \item Case B is slightly different, since an $i$ such that for all $k$, $(a_{k,l})_{l<\omega} {\Rightarrow} a_{i,\infty}$  or $(a_{k,-l})_{l<\omega} {\Rightarrow} a_{i,\infty}$ does not necessarily exist either. We may distinguish three subcases:
    \begin{enumerate}
        \item there is $i$ such that for $\lambda$-many $k$, $(a_{k,l})_{l< \omega}$ or $(a_{k,-l})_{l<\omega}$ pseudo-converge to $a_{i,0}$. We conclude as in the proof.
        \item there is $i<\lambda$ such that $(a_{i,j})_{j< \omega}$  pseudo-converges to $a_{k,0}$ for $\lambda$-many $k$. For such a $k$, we have ${\val(a_{k,0}-a_{i,\infty}) > \val(a_{k,0}-a_{i,0}) \geq \gamma }$, and thus $ {\val(a_{k,0}-a_{i,\infty}) +\delta_n > \gamma_k }$. We may conclude as well.
        \item there is $i<\lambda$ such that $(a_{i,-j})_{j< \omega}$  pseudo-converges to $a_{k,0}$ for $\lambda$-many $k$ . This is an analogue to Subcase (2) just above, where $a_{i,-\infty}$ is taking the place of $a_{i,\infty}$.
    \end{enumerate}
    \end{itemize}
    Hence, we get $\lambda=\kappa$.\\
    If $\kappa$ is of cofinality $\omega$, let $(\lambda_k)_{k\in \omega}$ be a sequence of successor cardinals cofinal in $\kappa$. By the previous discussion, we find an inp-pattern in $\RV_{<\omega}$ of depth $\lambda_k$ for each $\lambda_k$. Hence, $\lambda_k \leq \lambda$ and $\kappa = \lambda$ or $\kappa = \act(\lambda)$.
\end{proof}
\begin{remark}\label{rmkreductionRV}
\begin{itemize}
    \item Consider now an enriched Henselian valued field $\mathcal{K}=(K,\RV_{<\omega}, \ldots)$ of characteristic $(0,p)$, $p\geq 0$ in an $\RV_{<\omega}$-enrichment $\mathrm{L}_{\RV_{<\omega},e}$ of $\mathrm{L}_{\RV_{<\omega}}$. Then, the above proof still holds. The burden of $K$ is equal (modulo the same subtleties when we consider the burden in $\Card^\star$) to the burden of $\RV_{<\omega} \cup \Sigma_e$ with the induced structure, where $\Sigma_e$ is the set of new sorts in $\mathrm{L}_{\RV_{<\omega},e} \setminus \mathrm{L}_{\RV_{<\omega}}$.
    \item Similarly, an $\RV$-enriched benign Henselian valued field has the same burden as  $\RV \cup \Sigma_e$ where $\Sigma_e$ is the set of new sorts in $\mathrm{L}_{\RV,e} \setminus \mathrm{L}_{\RV}$.
\end{itemize}

\end{remark}

\subsubsection{ Applications to $p$-adic fields} \label{subsubsac}

In this paragraph, $p$ is a prime number. We will deduce from Theorem \ref{ThmHensValuedFieldReductionRV}, as an application, that any finite extension of $\mathbb{Q}_p$ is dp-minimal. This is already known (in fact, all local fields of characteristic $0$ are dp-minimal). One can refer to the classification on dp-minimal fields by Will Johnson \cite{Joh15}. The fact that $\mathbb{Q}_{p}$ is dp-minimal is due to Dolich, Goodrick and Lippel \cite[Section 6]{Dol11} and Aschenbrenner, Dolich, Haskell, Macpherson and Starchenko in \cite[Corollary 7.9.]{ADHMS}. 
In Section \ref{SectionUnrimified}, we will study more generally unramified mixed characteristic Henselian valued fields.

\begin{theorem}
	The theory of any finite extension of $\mathbb{Q}_p$ in the language of rings is dp-minimal.
\end{theorem}

A characterisation of dp-minimality is the following: for any mutually indiscernible sequences $(a_i)_{i <\omega}$ and $(b_i)_{i<\omega}$ and any point $c$, one of these two sequences is indiscernible over $c$. As we already mention earlier, a theory is dp-minimal if and only if it is NIP and inp-minimal (see \cite[Lemma 1.4]{Sim11} ). Since finite extensions of $\mathbb{Q}_p$ are NIP, we have to prove that they are inp-minimal. Recall first that the valuation in a finite extension of $\mathbb{Q}_p$ is definable in the language of rings:
\[\val(x)\geq 0 \ \Leftrightarrow \ \exists y \ 1+\pi x^q= y^q, \]
where $\pi$ is an element of minimal positive valuation and $q$ is a prime with $q \neq p$. 
We can safely consider $\mathbb{Q}_p$ in the two-sorted language of valued fields $\mathrm{L} = \mathrm{L}_{\text{Mac}} \cup \mathrm{L}_{\text{Pres}} \cup \lbrace \val \rbrace$,
where $\mathrm{L}_{\text{Mac}}= \mathrm{L}_{\text{Rings}} \cup \lbrace P_n \rbrace_{n \geq 2}$ is the language of Macintyre with a predicate $P_n$ for the subgroup of $n{\text{th}}$-power of $\mathbb{Q}_p$ and where $\mathrm{L}_{\text{Pres}}$ is the language of Presburger arithmetic. We have the following well known result, that we already discussed in the example below Proposition \ref{PropositionEquivPureRelativeQuantifierElimination}:
\begin{fact}
The theory $\Th(\mathbb{Q}_p)$ eliminates quantifiers. In particular, the value group is a pure sort. 
\end{fact}

Let $\mathcal{K}=(K,\Gamma)$ be a finite extension of $\mathbb{Q}_p$ and let $\pi \in K$ be an element of minimal positive valuation. By interpretability, we obtain:
\begin{remark}
    The value group $\Gamma$ is purely stably embedded in $\mathcal{K}$. Since $\Gamma$ is a Z-group (as a finite extension of a $Z$-group), it is in particular inp-minimal.
\end{remark}

Fix some $n \in \mathbb{N}$. We have the following exact sequence
\[\xymatrix{1 \ar@{->}[r]& \mathcal{O}^\times /(1+\mathfrak{m}_{\delta_n})\ar@{->}[r]  & \RV_{\delta_n}^{\star} \ar@{->}[r]^{\val_{\rv_{\delta_n}}} & \Gamma \ar@{->}[r]  & 0 },\]
where $\delta_n=\val(p^n)$ and $\mathfrak{m}_{\delta_n}=\{ x\in K \ \vert \val(x)> \val(p^n) \}$. 
 One sees that $(\mathcal{O}/\mathfrak{m}_{\delta_n})^\times \simeq \mathcal{O}^\times /(1+\mathfrak{m}_{\delta_n})$ is finite, or in other words, that the valuation map $\val_{\rv_{\delta_n}}$ is finite to one. It follows by Lemma \ref{Ft1Bdn} that  $\RV_{\delta_n}$ is also inp-minimal. Since this holds for arbitrary $n\in \mathbb{N}$, $\RV=\bigcup_n \RV_{\delta_n}$ is inp-minimal. We conclude by using Theorem \ref{ThmHensValuedFieldReductionRV}. \\

The next application is a anticipation of the next paragraph. We provide a new proof of the non-uniform definability of an angular component. It can in fact already be deduced from \cite{CS19}. Recall that an angular component is a group homomorphism $\ac: (K^\star,\cdot) \rightarrow (k^\star,\cdot)$ such that $\restriction{ac}{\mathcal{O}^\times}=\restriction{\res}{\mathcal{O}^\times}$.
$$\xymatrix{
1 \ar@{->}[r] & O^\times \ar@{->}[r]\ar@{->}[d]_{\res}& K^\star \ar@{->}[r]_\val\ar@{->}[d]_{\rv}\ar@/_1.0pc/@{->}[dl]^\ac & \Gamma\ar@{->}[r] \ar@{=}[d] &  0 \\
1 \ar@{->}[r]& O^\times /1+\mathfrak{m} \simeq k^\star\ar@{->}[r]  & \RV^{\star} \ar@/_1.5pc/@{-->}[l]^{\ac_{\rv}}\ar@{->}[r]^{\val_{\rv}} & \Gamma \ar@{->}[r]  & 0 }$$ 
    
    Consider any theory $T_{\ac}$ of a valued field endowed with an $\ac$-map, and assume that both the value group $\Gamma$ and residue field $k$ are infinite. Then by Fact \ref{FactSumBurdenOrthSorts} and bi-interpretability on unary sets, on sees that the $\RV$-sort is of burden at least $2$. The set $\RV^\star$ is indeed in definable bijection with the direct product $\Gamma\times k^\star$.

    In the field of $p$-adics $\mathbb{Q}_p$, an angular component $\ac$ is definable in the language of rings. We can show now easily that this definition cannot be uniform:

\begin{corollary} \index{Angular component}
There is no formula which gives a uniform definition of an $\ac$-map in $\mathbb{Q}_p$ for every prime $p$.
\end{corollary}
 Notice that this has already been observed by Pas in \cite{Pas90}.  
\begin{proof}
By Chernikov-Simon \cite{CS19}, we know that the ultraproduct of $p$-adic $\mathcal{F}=\prod_{\mathcal{U}}\mathbb{Q}_p$, where $\mathcal{U}\subset \mathcal{P}$ is an ultrafilter on the set of primes, is inp-minimal in the language of rings (recall that the p-adic valuation is uniformly definable in $\mathrm{L}_{\text{Rings}}$). The residue field and the value group are infinite since they are respectively a pseudo-finite field and a $\mathbb{Z}$-group. By the above discussion, the $\ac$-map cannot be defined in the language of rings, as it would contradict inp-minimality.
\end{proof}

    \subsection{Benign Henselian Valued Fields}\label{SectionHenselianValuedFieldsEquichar0}
    Let $\mathcal{K}=(K,\Gamma,k)$ be a saturated enough benign Henselian valued field. We will compute the burden of $\RV:= K^\star/1+\mathfrak{m}$ in terms of the burden of $k$ and $\Gamma$. As the $\RV$-sort is stably embedded, we will consider it as a structure on its own. By Fact \ref{FactRVBiInterpretability}, the induced structure is given by:
 \[\left\{\RV,(k,\cdot,+,0,1), (\Gamma,+,0,<),\val_{\rv}:\RV \rightarrow \Gamma, k^\star \rightarrow \RV  \right\}. \]
Notice in particular that there is no need of the symbol $\oplus$ as we consider the sort $k$ and $\Gamma$ instead. The language is denoted by $\mathrm{L}$. In other words the sort $\RV$ is no more than an enriched exact sequence of abelian groups:
     \[1 \rightarrow k^\star \rightarrow \RV^\star \overset{\val_{\rv}}{\rightarrow} \Gamma \rightarrow 0,\]
     where $k=k^\star  \cup \{0\}$ is endowed with its field structure and $\Gamma$ is endowed with its ordered abelian group structure.  As $\Gamma$ is torsion free, $k^\star$ is a pure subgroup of $\RV^\star$. The idea to consider $\RV$ as an enrichment of abelian groups is already present in \cite{CS19} and has been  developed in \cite{ACGZ20}.
        
\subsubsection{Reduction of burden to $\Gamma$ and $k$}
Let us recall a  result of Chernikov and Simon:

    \begin{theorem}[{\cite[Theorem 1.4]{CS19}}] 
    Assume $\mathcal{K}$ is a Henselian valued field of equicharacteristic $0$. Assume the residue field $k$ satisfies
    \begin{equation}
            \tag{$H_{k}$} k^\star/(k^\star)^p \ \text{is finite for every prime p.}
    \end{equation}
    Then $\mathcal{K}$ is inp-minimal if and only if $\RV$ with the induced structure is inp-minimal if and only if $k$ and $\Gamma$ are both inp-minimal. 
    \end{theorem}
     It will now be easy to extend this theorem. We have already seen the reduction to the $\RV$-sort for any benign Henselian valued field, without the assumption $(H_k)$. For the reduction to $\Gamma$ and $k$, one can give first an easy bound, also independent of the assumption $(H_k)$. Indeed, recall that in an $\aleph_1$-saturated model, any pure exact sequence of abelian groups splits (Fact \ref{FactSectionPureSubgroupAleph1Saturated}).  In particular, there exists a section $\ac_{\rv}: \RV^{\star} \rightarrow k^{\star} $ of the valuation $\val_{\rv}$ or equivalently, there exists an angular component $\ac: K^{\star} \rightarrow k^{\star}$ (as we already discussed in Paragraph \ref{subsubsac}). 
     
     $$\xymatrix{
1 \ar@{->}[r] & \mathcal{O}^\times \ar@{->}[r]\ar@{->}[d]_{\res}& K^\star \ar@{->}[r]_\val\ar@{->}[d]_{\rv}\ar@/_1.0pc/@{->}[dl]^\ac & \Gamma\ar@{->}[r] \ar@{=}[d] &  0 \\
1 \ar@{->}[r]& \mathcal{O}^\times /1+\mathfrak{m} \simeq k^\star\ar@{->}[r]  & \RV^{\star} \ar@/_1.3pc/@{->}[l]^{\ac_{\rv}} \ar@{->}[r]^{\val_{\rv}} & \Gamma \ar@{->}[r]  & 0 }$$

     Recall that $\mathrm{L}_{\ac}$ is the language $\mathrm{L}$ extended by a unary function $\ac_{\rv}: k \rightarrow \RV$. A direct translation of Fact \ref{facttrivialbound} gives:
    
    \begin{fact}[Trivial bound]\label{factbdnRV}
    We have $\bdn_{\mathrm{L}}(\Gamma) = \bdn_{\mathrm{L}_{\ac}}(\Gamma)$ and $\bdn_{\mathrm{L}}(k) = \bdn_{\mathrm{L}_{\ac}}(k)$ as well as the following:
        $$\bdn_{\mathrm{L}}(\RV) \leq \bdn_{\mathrm{L}_{\ac}}(\RV) = \bdn_{\mathrm{L}}(\Gamma)+\bdn_{\mathrm{L}}(k). $$
    \end{fact}
    
    A valued field $\mathcal{K}_{\ac}$ together with an angular component $\ac$ can be considered as an $\RV$-enrichment of $\mathcal{K}$. Using the enriched version of Theorem \ref{ThmHensValuedFieldReductionRV} (see Remark \ref{rmkreductionRV}), we get:
    \begin{theorem}\label{ThmBdnHenselianValuedFieldequicaracteristicAngularComponent}
    Let $\mathcal{K}_{\ac}=(\mathcal{K}, \Gamma, k, \val, \ac)$ be a benign Henselian valued field endowed with an $\ac$-map. Then:
    $\bdn(\mathcal{K}_{\ac})= \bdn(k)+ \bdn(\Gamma)$.
    \end{theorem}
    
    
    If we do not want to consider an $\ac$-map, we can also compute the burden using the torsion-free case of Theorem \ref{ThmBdnExSeq} together with Theorem \ref{ThmHensValuedFieldReductionRV}. We get:
    \begin{theorem}\label{ThmBdnHenValFieCha00} \index{Burden}\index{Valued field! benign Henselian}
        Let $\mathcal{K}$ be a benign Henselian valued field. Then: 
        \[\bdn(\mathcal{K})= \max_{n\geq 0} \left(\bdn(k^\star/{k^\star}^n) + \bdn (n\Gamma) \right).\]
    \end{theorem}
    
    This gives a full answer to \cite[Problem 4.3]{CS19} and \cite[Problem 4.4]{CS19}:
    \begin{corollary}\label{thmgnrl}
        Let $\mathcal{K}=(K,\RV,k,\Gamma)$ be a benign Henselian valued field. Assume that:
        \begin{equation}
            \tag{$H_k$} k^\star/(k^\star)^p \ \text{is finite for every prime p.}
        \end{equation}
        Then we have the equalities
        $$\bdn(\mathcal{K})= \bdn(\RV)= \max(\bdn(k),\bdn(\Gamma)).$$
    \end{corollary}
    Also, in the case that $\mathcal{K}$ is not trivially valued, the value group $\Gamma$ is necessary of burden $\bdn(\Gamma)>0$. It follows that a non-trivially valued benign Henselian field $\mathcal{K}$ is inp-minimal if and only if $\Gamma$, $k$ are inp-minimal and $k$ satisfies $(H_k)$.

    Similarly to the proof of non-existence of a uniform definition of the angular component of $\mathbb{Q}_p$, we can notice the following:
    \begin{remark} \label{CorollaryNoDefinableAngularComponent}
        Let $\mathcal{K}$ be a benign Henselian valued field of finite burden. Assume that the residue field is infinite and satisfies $(H_k)$. Then, no angular component is definable in the language of valued fields $\mathrm{L}_{div}$.  
    \end{remark}
    The reason is of course that in such a case, the two terms $\max(\bdn(k),\bdn(\Gamma))$ and $\bdn(k)+\bdn(\Gamma)$ are distinct.
    
    All these results hold resplendently. In fact by definition, a benign Henselian valued field can have an enriched value group and residue field. Let us clarify by stating the previous theorem in an enriched language:
    \begin{remark}\label{rmkchar00}
        If $\mathcal{K}=(K,\RV,k,\Gamma, \ldots)$ is a $\lbrace \Gamma \rbrace$-$\{k\}$-enriched benign Henselian valued field in a $\lbrace \Gamma \rbrace$-$\{k\}$-enrichment $\mathrm{L}_{\Gamma,k,e}$ of $\mathrm{L}_{\Gamma,k}$, then
        \[\bdn(\mathcal{K})= \bdn(\RV\cup \Sigma_e)= \max_n(\bdn(k^\star/{k^\star}^n) + \bdn (n\Gamma),\bdn(\Sigma_e)),\]
        where $\Sigma_e$ is the set of new sorts in $\mathrm{L}_{\Gamma,k,e}\setminus \mathrm{L}_{\Gamma,k}.$
    \end{remark}

    

    To conclude this short subsection, let us discuss more on the hypothesis $(H_k)$ and bounded fields.
    
    \subsubsection{Bounded fields and applications}
    A \textit{bounded field}\index{Bounded field} is a field with finitely many extensions of degree $n$ for every integer $n$. The absolute Galois group is called \textit{small} if it contains finitely many open subgroups of index $n$. These two conditions are equivalent for perfect fields: a perfect field is bounded if and only if its absolute Galois group is small. Such a field $K$ satisfies in particular the following:
    \begin{equation}
        \tag{$H$} K^\star/(K^\star)^p \ \text{is finite for every prime p,}
    \end{equation}
    (see for example \cite[Proposition 2.3]{FJ16}), and it's clear that $(H)$ implies $(H_k)$. It also implies:
    \begin{equation}
        \tag{$H_\Gamma$} \Gamma/p \Gamma \ \text{is finite for every prime p,}
    \end{equation}
    
    The condition $H_k$ might be restrictive but it allows various burdens for the residue field. However, the condition $H_{\Gamma}$ implies inp-minimality \index{Inp-minimality}for the value group. Indeed, an abelian group $\Gamma$ satisfying $(H_{\Gamma})$ is called \textit{non-singular}. In the pure structure of ordered abelian groups, non-singular ordered abelian groups are exactly the dp-minimal ones (see \cite[Theorem 5.1]{JSW}). We have the following examples :
    \begin{examples}
        \begin{itemize}
            \item The Hahn field $\mathbb{F}_p^{alg}((\mathbb{Z}[1/p]))$ is algebraically maximal Kaplansky Henselian. By Jahnke, Simon and Walsberg, the value group  $\mathbb{Z}[1/p]$ is inp-minimal as it satisfies $(H_\Gamma)$. The residue field $\mathbb{F}_p^{alg}$ satisfies $(H_k)$ and is inp-minimal. By Theorem \ref{ThmBdnHenValFieCha00}, this Hahn field is inp-minimal.
            \item In general, a bounded benign Henselian valued field $\mathcal{K}$ with residue field $k$ has burden $\max(\bdn(k), 1)$. 
        \end{itemize}
    \end{examples}

    Montenegro has computed the burden of some theories of bounded fields, namely bounded pseudo real closed fields (PRC fields) and pseudo $p$-adicaly closed fields (PpC fields). We recall here these theorems (see \cite[Theorems 4.22 \& 4.23]{Mon17}):
    \begin{theorem}
        Let $k$ be a bounded PRC field. Then $\Th(k)$ is $\NTP_2$, strong and of burden the (finite) number of orders in $k$.
    \end{theorem}

    \begin{theorem}
        Let $k$ be an PpC field. Then $\Th(k)$ is $\NTP_2$ if and only if $\Th(k)$ is strong if and only if $k$ is bounded. In this case, the burden of $\Th(k)$ is the (finite) number of p-adic valuations in $k$.   
    \end{theorem}

    \subsection{Unramified mixed characteristic Henselian Valued Fields}\label{SectionUnrimified}

    Let $\mathcal{K}=(K,\RV_{<\omega},\Gamma,k)$ be an unramified Henselian valued field of characteristic $(0,p)$, $p\geq 2$ with perfect residue field $k$. 
    We denote by $1$ the valuation of $p$. The value group $\Gamma$ contains $\mathbb{Z}\cdot 1$ as a convex subgroup. Recall that in this context, it is more convenient to denote the $n^{\text{th}}$ $\RV$-sort by $\RV_n^\star:= K^\star/(1+\mathfrak{m}^n)$  where $\mathfrak{m}=\lbrace x\in K \ \vert \ \val(x)>0 \rbrace$ is the maximal ideal of the valuation ring $\mathcal{O}$. Notice that $\mathfrak{m}^n= p^n\mathcal{O}$ for every integer $n$. Similarly to the previous section, we will compute the burden of $\RV_{<\omega} = \cup_{n\in \mathbb{N}}\RV_n$ in terms of the burden of $k$ and $\Gamma$. 
     
    \subsubsection{Reduction from $\RV_{<\omega}$ to $\Gamma$ and $k$}
    Now we can look for the burden of $\RV_n$. We start with a harmless observation:
    
    \begin{observation}
        Let $m<n$ be integers. The element $p^m$ is of valuation $m$. By \cite[Proposition 2.8]{Fle11}, $\RV_{m}$ is $\emptyset$-interpretable in $\RV_{n}$, with base set $\RV_{n}$ quotiented by an equivalence relation. Hence the burden of $\RV_n$ can only grow with $n$:  for $m<n$, $\bdn(\RV_m) \leq \bdn (\RV_n)$.
    \end{observation}

    Recall that in this context of unramified mixed characteristic Henselian valued fields with perfect residue field, the $n^{\text{th}}$ residue ring $\mathcal{O}_n:= \mathcal{O}/ p^n\mathcal{O}$ is isomorphic to the $n$-truncated ring of Witt vectors (see Proposition \ref{KerValRV_n}). We work now in the following languages:
    \begin{align*}
    \mathrm{L}= \lbrace & K,\Gamma,(\RV_n)_{n<\omega},(W_n(k))_{n<\omega},
    \val: K^\star\rightarrow \Gamma,\\ 
    &(\res_n: \mathcal{O}\rightarrow W_n(k))_{n<\omega}, (\rv_n:K^\star \rightarrow \RV_n)_{n<\omega}\rbrace ,
    \end{align*}

    which is a little variation of (and bi-interpretable with) the language $\mathrm{L}_{\RV_{<\omega}}$, where the structure of the $\RV_n$'s is described with exact sequences.  We can also add the $\ac$-maps to this language:
    
    \[    \mathrm{L}_{ac_{<\omega}} = \mathrm{L} \cup \lbrace (\ac_n: K^\star \rightarrow W_n(k))_{n<\omega}\rbrace.
    \]  
    
    Here is a consequence of Corollary \ref{CorollaryWnBiInterpretableWithk}, Remark \ref{k to the n kappa n} and Fact \ref{FactBdnInterpretUnarySet}:
     \begin{corollary}\label{BurdenWn}
        We have:
        \begin{itemize}
            \item  ${\bdn(W_n(k))=\kappa_{inp}^1(W_n(k))=\kappa_{\inp}^{n}(k)}$. 
            \item $\bdn((W(k), +, \cdot, \pi:W(k) \rightarrow k))= \kappa^{\aleph_0}_{inp}(k). $
            \end{itemize}
    \end{corollary}
    
    Recall that we have the following inequalities (see Paragraph    \ref{SubsectionClassificationtheory}): 
    \[ n\cdot\kappa_{\inp}^{1}(k) \leq \kappa_{\inp}^{n}(k) \qquad \kappa_{\inp}^{n}(k)+1 \leq (\kappa_{\inp}^1(k)+1)^n. \]
    In particular, if $k$ is infinite then the burden of $(W_n(k),+,\cdot, \pi)$ is at least $n$.

    In the language $\mathrm{L}_{\ac_{<\omega}}$, a consequence of Proposition \ref{KerValRV_n} is that, for every $n<\omega$ the sort $W_n(k)$ is pure (in particular stably embedded) and orthogonal to $\Gamma$,  as it is $\emptyset$-bi-interpretable with $(k^n,+,\cdot,p_i, i<n)$, which is a pure sort orthogonal to $\Gamma$. It follows that $W_n(k)$ doesn't have more structure in $\mathrm{L}_{\ac_{<\omega}}$ than in $\mathrm{L}$. Similarly, the burden of $\Gamma$ is the same in any of the above languages.  Hence, we actually have the following equalities:
        \begin{equation}\label{eqBdnW_n2}
     \bdn_{\mathrm{L}}(W_n(k))= \bdn_{\mathrm{L}_{\ac_{<\omega}}}(W_n(k)),
    \end{equation}
        \begin{equation}
     \bdn_{\mathrm{L}}(\Gamma)= \bdn_{\mathrm{L}_{\ac_{<\omega}}}(\Gamma).
    \end{equation}
    We are now able to give a relationship between $\bdn(\RV_n)$ and $\bdn(W_n(k))$. 

    \begin{proposition}\label{ProBurdCharp} [Trivial bound]
        We have 
        \begin{align*}\max(\bdn_{\mathrm{L}}(W_n(k)),\bdn_{\mathrm{L}}(\Gamma)) \leq \bdn_{\mathrm{L}}(\RV_n) &\leq \bdn_{\mathrm{L}_{\ac_{<\omega}}}(\RV_n) \\
        &= \bdn_{\mathrm{L}}(W_n(k))+\bdn_{\mathrm{L}}(\Gamma).
        \end{align*}
    \end{proposition}
    
    \begin{proof}
    By Proposition \ref{KerValRV_n}, we have the exact sequence of abelian groups: 
    \[ 1 \rightarrow W_n(k)^\times \rightarrow \RV_n^\star \rightarrow \Gamma \rightarrow 0.\]
    
    The first inequality is clear if one shows that $\bdn_{\mathrm{L}}(W_n(k))=\bdn_{\mathrm{L}}(W_n(k)^\times)$ where $W_n(k)^\times$ is endowed with the induced structure.
    The second inequality is also clear, as adding structure can only make the burden grow. 
    Let $\lbrace \phi_i(x,y_i), (a_{i,j})_{j<\omega} \rbrace_{i \in \lambda}$ be an inp-pattern in $W_n(k)$, with $(a_{i,j})_{i<\lambda,j<\omega}$ mutually indiscernible. Let $d\models \lbrace \phi(x,a_{i,0})\rbrace_{i\in \lambda}$ be a realisation of the first column. In the case where $d\in W_n(k)^\times$, there is nothing to do. Otherwise, $1+d \in W_n(k)^\times$ and $\lbrace \phi_i(x-1,y_i), (a_{i,j})_{j<\omega} \rbrace_{i \in \lambda}$ is an inp-pattern in $W_n(k)^\times$ of depth $\lambda$. This concludes the proof of the first inequality.
   
    We work now in $\mathrm{L}_{\ac_{<\omega}}$, where we interpret $(ac_n)_n$ as a compatible sequence of angular components (it exists by $\aleph_1$-saturation). Recall that the burden may only increase. Then, the above exact sequences (definably) split in $\mathrm{L}_{\ac_{<\omega}}$, as we add a section. By the previous discussion, $W_n(k)^\times$ and $\Gamma$ are orthogonal and stably embedded. We apply now Fact \ref{FactSumBurdenOrthSorts}:  the burden  $\bdn_{\mathrm{L}_{\ac_{<\omega}}}(\RV_n^\star)$ is equal to $\bdn_{\mathrm{L}_{\ac_{<\omega}}}(W_n(k)^\times)+ \bdn_{\mathrm{L}_{\ac_{<\omega}}}(\Gamma)=\bdn_{\mathrm{L}}(W_n(k))+ \bdn_{\mathrm{L}}(\Gamma)$. 

    \end{proof}
    
    Combining Corollary \ref{corollaryinfiniteburden}, Corollary \ref{BurdenWn} and Proposition \ref{ProBurdCharp}, one gets:

    \begin{theorem}\label{theoremmixedchar}\index{Burden}\index{Angular component! of order $n$}
        Let $\mathcal{K}=(K,k,\Gamma)$ be an unramified mixed characteristic Henselian valued field. We denote by $\mathcal{K}_{\ac_{<\omega}}=(K,k,\Gamma,\ac_n,n<\omega)$ the structure $\mathcal{K}$ endowed with compatible $\ac$-maps. Assume the residue field $k$ is perfect. One has 
        \[\bdn(\mathcal{K})= \bdn(\mathcal{K}_{\ac_{<\omega}})= \max(\aleph_0 \cdot \bdn(k), \bdn(\Gamma)).\]
        
    \end{theorem}
    And its enriched version:
    \begin{remark} \label{rmkmixedchar}
        Let $\mathrm{L}_e$ be a $\lbrace \Gamma \rbrace$-$\{k\}$-enrichment of $\mathrm{L}$. Let $\mathcal{K}=(K,k,\Gamma,\ldots)$ be an enriched unramified mixed characteristic Henselian valued field in the language $\mathrm{L}_e$. Assume the residue field $k$ is perfect. We denote by $\mathcal{K}_{\ac_{<\omega}}=(K,k,\Gamma,\ac_n,n<\omega,\ldots)$ the structure $\mathcal{K}$ endowed with compatible $\ac$-maps.
        One has 
        \[\bdn(\mathcal{K})=\bdn(\mathcal{K}_{\ac_{<\omega}})= \max(\aleph_0 \cdot \bdn(k), \bdn(\Gamma), \bdn(\Sigma_e)),\]
        where $\Sigma_e$ is the set of new sorts in $\mathrm{L}_e \setminus \mathrm{L}$.
    \end{remark}

    This is a simple calculation, unless we want to consider burden in $\Card^\star$.

    \begin{remark}
        Let  $\mathcal{K}$ and  $\mathcal{K}_{\ac_{<\omega}}$ as above. We consider the second definition of burden (Definition \ref{DefOpeCarEto}).
        We have 
        \[\bdn(\mathcal{K})= \bdn(\mathcal{K}_{\ac_{<\omega}})= \max(\aleph_0 \cdot^\star \bdn(k), \bdn(\Gamma)).\]
    \end{remark}
    
    This is what we will prove now. Remark that it implies that an unramified mixed characteristic valued field of infinite perfect residue field is never strong. 
    
    \begin{proof}
        We use the same notation as before in this section. Unfortunately, due to the ambiguity in Corollary \ref{corollaryinfiniteburden} concerning $\bdn(\mathcal{K}) \in \{\bdn(\RV_{<\omega}),\act(\bdn(\RV_{<\omega}))\}$ in the case that $\cf(\bdn(\RV_{<\omega}))=\omega$, we have to go back to the proof of Theorem \ref{ThmHensValuedFieldReductionRV}.
        
        We first show that $\bdn(\mathcal{K})$ is at least $\aleph_0 \cdot^\star \bdn(k)$. Recall that $W_n(k)\simeq \mathcal{O}_n:=\mathcal{O}/\mathfrak{m}^n$ is interpretable (with one-dimensional base set $\mathcal{O}\subset K$), and so is the projective system $\{W_n(k), \pi_{n,m}: W_n(k) \rightarrow W_m(k), \ n>m\}$ and the projection maps $\chi_{n,n}: W_n(k) \rightarrow k, x=(x_1,\ldots ,x_n) \mapsto x_n$. If $\cf(\bdn(k))>\aleph_0$, there is nothing to do as $\aleph_0 \cdot \bdn(k)=\aleph_0 \cdot^\star \bdn(k)$. Assume $\cf(\bdn(k))\leq \aleph_0$. We write $\bdn(k)=\sup_{n<\omega}\lambda_n$ with $\lambda_n \in \Card$. Let $P_n(x_k)$ be an inp-pattern with $x_k \in k$, $\vert x_k \vert =1$, of depth $\lambda_n$ for every $n\in \omega$.
        Then, the pattern $P(x)=\cup_{n\in \omega} P_n(\chi_{n,n}(\pi_n(x)))$ is an inp-pattern in $K$ of depth $\aleph_0 \cdot^\star \bdn(k)$.
        One gets:
        \[\bdn(\mathcal{K})\geq \max(\aleph_0 \cdot^\star \bdn(k), \bdn(\Gamma)).\]
        We now prove that $\max(\aleph_0 \cdot^\star \bdn(k), \bdn(\Gamma))$ is an upper bound for $\bdn(\mathcal{K}_{\ac_{<\omega}})$. \\
        \textbf{Case 1:} $\aleph_0\cdot^\star \bdn(k) \geq \bdn(\Gamma)$.\\
        Subcase 1.A: $\cf(\bdn(k))> \aleph_0$. By Corollary \ref{corollaryinfiniteburden}, $\bdn(\mathcal{K}_{\ac_{<\omega}})= \bdn(\RV_{<\omega})=\sup_n(\kappa_{inp}^n(k), \bdn(\Gamma))= \bdn(k)= \aleph_0 \cdot \bdn(k)$. We used the submultiplicativity of the burden, which gives here $\kappa_{inp}^n(k)=\kappa_{inp}^1(k)=\bdn(k)$ for all $n \in \mathbb{N}$.\\
        Subcase 1.B: $\cf(\bdn(k)) \leq \aleph_0$. Then $\act(\bdn(\RV_{<\omega}))=\aleph_0\cdot \bdn(k)$. By Corollary  \ref{corollaryinfiniteburden}, we have $\bdn(\mathcal{K}_{\ac_{<\omega}})\leq \aleph_0\cdot \bdn(k)$.\\
        \textbf{Case 2:} $\bdn(\Gamma) > \aleph_0\cdot^\star \bdn(k)$. If $\bdn(\Gamma)$ is in $\Card$, this is clear by Corollary  \ref{corollaryinfiniteburden}. Assume $\bdn(\Gamma)$ is of the form $\lambda_-$ for a limit cardinal $\lambda\in \Card$. Notice that this case occurs only if the sort $\Gamma$ is enriched. 
        We work in the corresponding enrichment of language ${\mathrm{L}_{\RV_{<\omega}}}$ together with $\ac_n$-maps. We have to show that $\lambda_-$ is an upper bound for $\bdn(\mathcal{K}_{\ac_{<\omega}})$.  Let $P(x)=\lbrace \theta_i(x,y_{i,j}), (c_{i,j})_{j\in \bar{\mathbb{Z}}}\rbrace_{i\in \lambda}$ be an inp-pattern in $K$ of depth $\lambda$ with $\vert x \vert =1$ and $(c_{i,j})_{i<\lambda,j\in \bar{\mathbb{Z}}}$ be a mutually indiscernible array.
        Then, by Fact \ref{factfle}, one can assume that each formula $\theta_i(x,c_{i,j})$ in $P(x)$ ($i\leq \lambda$, $j\in \bar{\mathbb{Z}}$) is of the form 
            \[\tilde{\theta}_i(\rv_{n_i}(x-\alpha^1_{i,j}),\ldots, \rv_{n_i}(x-\alpha^m_{i,j}), \beta_{i,j}),\]
        for some integers $n_i$ and $m$, and where $\alpha^1_{i,j}, \ldots, \alpha^m_{i,j} \in K$, $\beta_{i,j} \in \RV_{n_i}$ and $\tilde{\theta}_{i}$ is an $\RV_{n_i}$-formula. As in the proof of Theorem \ref{ThmHensValuedFieldReductionRV}, we may assume with no restriction that $m=1$.
        As $\RV_{n_i}= W_{n_i}(k)^{\times} \times \Gamma$ is the direct product of the orthogonal and stably embedded sorts $W_{n_i}(k)^\times$ and $\Gamma$, we may assume $\theta_i(x,c_{i,j})$ is equivalent to a formula of the form
        \[\phi_{i}(\ac_{n_i}(x-\alpha_{i,j}),a_{i,j}) \wedge \psi_{i}(\val(x-\alpha_{i,j}),b_{i,j})\]
        where $\phi_{i}(x_{W_{n_i}},a_{i,j})$ is a $W_{n_i}$-formula and $\psi_{i}(x_\Gamma, b_{i,j})$ is a $\Gamma$-formula. By Claim \ref{mainprop} in Theorem \ref{ThmHensValuedFieldReductionRV} (or more precisely, by a generalisation of Claim \ref{mainprop} to infinite depth $M=\lambda$), one may assume that there is $k<\lambda$ such that for all $i<\lambda$, 
        \[\val(d-\alpha_{i,0}) \leq \min\lbrace \val(d-\alpha_{k,\infty}),\val(\alpha_{k,\infty}-\alpha_{i,0})\rbrace + \max(n_i,n_k).\]
        It follows that, if $\val(d-\alpha_{k,\infty}) = \val(\alpha_{k,\infty}-\alpha_{i,0})$, $\val(d-\alpha_{i,0})$ is equal to $\val(d-\alpha_{k,\infty})+n_i'$ for some $0\leq n_i'\leq \max(n_i,n_k)$. Otherwise, one has $\val(d-\alpha_{i,0}) = \min\lbrace \val(d-\alpha_{k,\infty}),\val(\alpha_{k,\infty}-\alpha_{i,0})\rbrace$.
        We can centralise $P(x)$ in $\alpha_{k,\infty}$, \textit{i.e.} we can assume that each formula in $P(x)$ is of the form 
        \[ \phi_{i}(\ac_{2n_i}(x-\alpha_{k,\infty}),a_{i,j}) \wedge \psi_{i}(\val(x-\alpha_{k,\infty}),b_{i,j})\]
        (we add new parameters $\val(\alpha_{k,\infty}-\alpha_{i,j})$ and $\ac_{2n
        _i}(\alpha_{k,\infty}-\alpha_{i,j})$. Notice that once the difference of the valuation is known, $\ac_{n_i}(d-\alpha_{i,j})$ can be computed in terms of $\ac_{2n_i}(d-\alpha_{k,\infty})$ and $\ac_{2n_i}(\alpha_{i,j}-\alpha_{k,\infty})$).
        By indiscernibility, at least one of the following sets
        \[\{\phi_{i}(x_{W_{2n_i}},a_{i,j})\}_{j<\omega}\]
        and 
        \[\{\psi_{i}(x_\Gamma,b_{i,j})\}_{j<\omega}\]
        is inconsistent.
        Since $\lambda> \sup_n \bdn(W_n(k))$, we may assume that
        \[\lbrace \psi_{i}(x_\Gamma,y_i), (b_{i,j})_{j\in \bar{\mathbb{Z}}}) \rbrace_{i<\lambda} \]
        is an inp-pattern in $\Gamma$. This is a contradiction. Hence, we have $\bdn(\mathcal{K})=\lambda_-$.
    \end{proof}

    We end now with examples:  
    \begin{examples}
        \begin{enumerate}
            \item Assume that $k$ is an algebraically closed field of characteristic $p$, and $\Gamma$ is a $\mathbb{Z}$-group. Then $\Gamma$ is inp-minimal, \textit{i.e.} of burden one (as it is quasi-o-minimal), and one has $\kappa^n_{inp}(k)=n$. By Theorem \ref{theoremmixedchar}, any Henselian mixed characteristic valued field of value group $\Gamma$ and residue field $k$ has burden $\aleph_0$. In particular, the quotient field $Q(W(k))$ of the Witt vectors $W(k)$ over $k$ is not strong.  
            \item Consider once again the field of $p$-adics $\mathbb{Q}_p$. We have $\kappa_{inp}^n(\mathbb{F}_p)=0$ for all $n$, and $\bdn(\mathbb{Z})=1$. Then Theorem \ref{theoremmixedchar} gives $\bdn(\mathbb{Q}_p)=1$. 
        \end{enumerate}
    \end{examples}
    \newpage

\appendix

\section{Reduction of burden in lexicographic products}\label{LexicographicProduct}
 Meir defined and studied the lexicographic product of relational structures in \cite{Mei16}. Using his quantifier elimination result, he notably proved a Stable and NIP transfer.  We continue here to investigate  the model theoretic  complexity of such products with respect to the burden. In this last section, we show that the burden of the lexicographic product of pure relational structures is the maximum of the burden of these structures. This is a reduction principle for pure relational structures which appears to be similar to that of pure short exact sequences of abelian groups. However, the situation is here simpler since terms are trivial and inp-patterns are automatically centralised. 

Consider a relational language $\mathrm{L}$.

\begin{definition} \index{Lexicographic product! of structures}
    Let $\mathcal{M}$ and $\mathcal{N}$ be two $\mathrm{L}$-structures. We consider the language $\mathrm{L}_{\mathbb{U},s} =\mathrm{L} \cup \{R^{\mathbb{U}} \}_{R \in \mathrm{L}}  \cup  \{s\}$ where $R^{\mathbb{U}}$ are new unary predicates and $s$ is a binary predicate. 
    The lexicographic product $\mathcal{M}[\mathcal{N}^s]^{\mathbb{U}}$ of $\mathcal{M}$ and $\mathcal{N}$ is the $\mathrm{L}_{\mathbb{U},s}$-structure of base set $M\times N$ where the relations are interpreted as follows:
     \begin{itemize}
         \item $s^{\mathcal{M}[\mathcal{N}^s]^{\mathbb{U}}} := \{ ((a,b),(a,b')) \ \vert \ a\in M, b,b' \in N\}$ 
         \item if $R\in \mathrm{L}$ is an $n$-ary predicate, 
         \begin{align*}
            R^{\mathcal{M}[\mathcal{N}^s]^{\mathbb{U}}} :=& \\
            & \left\{\left((a,b_1),\ldots, (a,b_n)\right) \ \vert \ a\in M \ \text{and } \mathcal{N} \models R(b_1,\ldots,b_n)\right\} \ \cup \\
            & \left \{\left((a_1,b_1),\ldots, (a_n,b_n)\right) \ \vert \ {\bigvee}_{0\leq i \neq j \leq n}a_i \neq a_j \ \text{and } \mathcal{M} \models R(a_1,\ldots,a_n) \right\}.
         \end{align*}
         \item if $R\in \mathrm{L}$ is an $n$-ary predicate,
         \begin{align*}
             {R^{\mathbb{U}}}^{\mathcal{M}[\mathcal{N}^s]^{\mathbb{U}}} :=& \bigg\{(a,b) \ \vert \ \mathcal{M}\models R(\underbrace{a, \ldots, a}_{n \text{ times}})\bigg \}.
         \end{align*} 
     \end{itemize}
     We denote by $\mathcal{M}[\mathcal{N}^s]$ the restriction of $\mathcal{M}[\mathcal{N}^s]^{\mathbb{U}}$ to $\mathrm{L}_s := \mathrm{L} \cup \{s\}$.
     
\end{definition}
    \begin{example}
        Consider the language of graphs $\mathrm{L}= \{R\}$, and let $\mathcal{M}$ and $\mathcal{N}$ be the graphs:
        
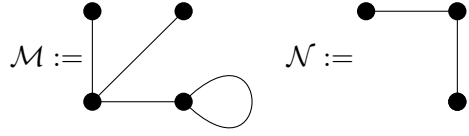
\begin{center} 
\begin{tikzpicture}[scale=1.2]\clip(-1,-0.4) rectangle (6,1.5);

\draw (-0.5,0.5) node {$\mathcal{M}:=$};
\fill (0,0) circle (0.1);
\fill (1,0) circle (0.1);
\fill (0,1) circle (0.1);
\fill (1,1) circle (0.1);
\draw (0,0) -- (1,0);
\draw (0,0) -- (0,1);
\draw (0,0) -- (1,1);
\draw (1,0) .. controls (2,1) and (2,-1) .. (1,0);
\draw (2.5,0.5) node {$\mathcal{N}:=$};
\fill (4,0) circle (0.1);
\fill (4,0) circle (0.1);
\fill (3,1) circle (0.1);
\fill (4,1) circle (0.1);
\draw (3,1) -- (4,1);
\draw (4,1) -- (4,0);
\end{tikzpicture}
\end{center}

We obtain the following graph:
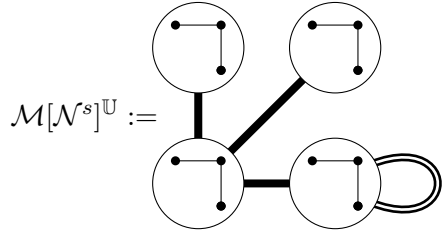
\begin{center} 
\begin{tikzpicture}[scale=0.6]\clip(-4,-1) rectangle (10,5);
\draw[line width=3pt] (0.5,1)--(0.5,2.5);
\draw[line width=3pt] (1.5,0.5)--(2.5,0.5);
\draw[line width=3pt] (0.5,0.5)--(3.5,3.5);
\draw[line width=1pt,double distance = 1pt] (3.5,0.5) .. controls (6.5,2.5) and (6.5,-1.5) .. (3.5,0.5);
\draw (-2,2) node {$\mathcal{M}[\mathcal{N}^s]^{\mathbb{U}}:=$};
\begin{scope}[shift={(0,0)}] 
\filldraw[fill=white] (0.5,0.5) circle (1);
\fill (1,0) circle (0.1);
\fill (1,0) circle (0.1);
\fill (0,1) circle (0.1);
\fill (1,1) circle (0.1);
\draw (0,1) -- (1,1);
\draw (1,1) -- (1,0);

\end{scope}
\begin{scope}[shift={(0,3)}]
\filldraw[fill=white] (0.5,0.5) circle (1);
\fill (1,0) circle (0.1);
\fill (1,0) circle (0.1);
\fill (0,1) circle (0.1);
\fill (1,1) circle (0.1);
\draw (0,1) -- (1,1);
\draw (1,1) -- (1,0);

\end{scope}
\begin{scope}[shift={(3,0)}]
\filldraw[fill=white] (0.5,0.5) circle (1);
\fill (1,0) circle (0.1);
\fill (1,0) circle (0.1);
\fill (0,1) circle (0.1);
\fill (1,1) circle (0.1);
\draw (0,1) -- (1,1);
\draw (1,1) -- (1,0);

\end{scope}
\begin{scope}[shift={(3,3)}]
\filldraw[fill=white] (0.5,0.5) circle (1);
\fill (1,0) circle (0.1);
\fill (1,0) circle (0.1);
\fill (0,1) circle (0.1);
\fill (1,1) circle (0.1);
\draw (0,1) -- (1,1);
\draw (1,1) -- (1,0);

\end{scope}
\end{tikzpicture}
\end{center}
where any point in a circle is related to any point from any linked circle. Notice that the loop is given by the predicate $R^{\mathbb{U}}$ (and not by $R$).

\end{example}

    Within $\mathcal{M}[\mathcal{N}^s]^{\mathbb{U}}$, the structure $\mathcal{M}$ can be seen as an imaginary. It is indeed the quotient of $\mathcal{M}[\mathcal{N}^s]^\mathbb{U}$ modulo the equivalence relation $s(\mathbf{x},\mathbf{y})$. We denote by $\pi_{\mathcal{M}}=\pi: \mathcal{M}[\mathcal{N}^s]^{\mathbb{U}} \rightarrow \mathcal{M}$ the projection. It is not an $\mathrm{L}$-homomorphism.  Indeed, if $R$ is an $n$-ary predicate of $\mathrm{L}$, and $D^n=\{(a,\ldots,a) \in M^n\}$ is the diagonal, one has:
    \[\pi(R^{\mathcal{M}[\mathcal{N}^s]^\mathbb{U}})= 
    \begin{cases}
      R^{\mathcal{M}}\cup D^n \ \text{ if $R^{\mathcal{N}}$ is not empty,}\\
    R^{\mathcal{M}}\setminus D^n \ \text{ if $R^{\mathcal{N}}$ is empty.}
    \end{cases}
    \]
    
  Nonetheless, we recover the $\mathrm{L}$-structure when we consider the additional symbol $R^{\mathbb{U}}$:
        \[R^{\mathcal{M}}=\pi(R^{\mathcal{M}[\mathcal{N}^s]^\mathbb{U}}) \setminus D^n \cup \{(a,\ldots,a) \ \vert \ a\in \pi({R^{\mathbb{U}}}^{\mathcal{M}[\mathcal{N}^s]^\mathbb{U}}) \}.\]

    For $a\in M$, the equivalence class of $a$ is a copy of $N$ and it is denoted by $N_a$. Here, by definition 
    \[\begin{array}{cll}
         f_a: & \mathcal{N} \rightarrow \mathcal{M}[\mathcal{N}^s]^{\mathbb{U}} \\
            & b \mapsto  (a,b)
    \end{array}\]
    is an $\mathrm{L}$-isomorphism onto $N_a$. We will see that there is no additional structure on these sets (see Corollary \ref{CoroStEm}). Here is the result of quantifier elimination obtained by Meir:
    
    \begin{theorem}
        Let $\mathrm{L}$ be a relational language. $\mathcal{M}$, $\mathcal{N}$ be two $\mathrm{L}$-structures admitting quantifier elimination. Then $\mathcal{M}[\mathcal{N}^s]^{\mathbb{U}}$ admits quantifier elimination.   
    \end{theorem}
        
        
    
    On can refer to {\cite[Theorem 2.6]{Mei16}} in the case where $\mathcal{M}$ is transitive (for the action of the automorphism group). The theorem above can be found in Meir's thesis \cite[Theorem 1.1.4]{Mei19}. We discuss here some immediate corollary. Let us denote by $\restriction{\mathcal{M}}{\mathbb{U}}$ the set $\pi(\mathcal{M}[\mathcal{N}^s]^{\mathbb{U}})$ with the full induced structure, and for $a\in M$, let $\mathcal{N}_a$ be the subset $N_a\subset M[N]$ with the full induced structure. 
    \begin{corollary}\label{CoroStEm}
        The structures $\restriction{\mathcal{M}}{\mathbb{U}}$ and $\mathcal{N}_a$ for $a\in M$ are stably embedded and setwise orthogonal. 
        The structures $\restriction{\mathcal{M}}{\mathbb{U}}$ and $\mathcal{M}$ on the one hand, and $\mathcal{N}_a$ and $\mathcal{N}$ on the other hand, have the same definable sets.
    \end{corollary}
    \begin{proof}
        By quantifier eliminations, it is enough to look at atomic formulas. It is then a simple case study.
    \end{proof}
        
    
     To simplify the notation in our next proof, we will make $\mathcal{M}$ an $\mathrm{L}_{\mathbb{U}}:= \mathrm{L} \cup \{ R^{\mathbb{U}}\}_{R \in \mathrm{L}}$- structure, where the unary predicate $R^{\mathbb{U}}(x)$ is interpreted as $R(\underbrace{x,\ldots,x}_{n \ \text{times}})$ if $R\in \mathrm{L}$ is an $n$-ary predicate. 
    \begin{theorem}\label{ThmBdnLexicoProd}
        $\bdn(\mathcal{M}[\mathcal{N}^s]^{\mathbb{U}})=\max(\bdn(\mathcal{M}), \bdn(\mathcal{N}))$.
    \end{theorem}
    
    The theorem holds also if we consider the second definition of burden (Definition \ref{DefBurden2}). The quantity $\max(\bdn(\mathcal{M}), \bdn(\mathcal{N}))$ is obviously a lower bound, as $\mathcal{M}[\mathcal{N}^s]^{\mathbb{U}}$ interprets on a unary set both $\mathcal{M}$ and $\mathcal{N}$.

    \begin{definition}
        Let $\mathbf{x}, \mathbf{y}^0, \ldots  ,\mathbf{y}^{n-1}$ be single variables. An $\{s\}$-formula $\phi_s(\textbf{x},\textbf{y}^0,\ldots,\textbf{y}^{n-1})$ is a complete $s$-diagram in $\mathbf{x}$ over $\mathbf{y}^0, \ldots  ,\mathbf{y}^{n-1}$ if it is a maximally consistent conjunction of $\neg s(\textbf{x},\textbf{y}^l)$ or $s(\textbf{x},\textbf{y}^l)$ for $l<n$. 
    \end{definition}
    
    \begin{proof}
        If $\mathcal{M}$ and $\mathcal{N}$ are finite, then so is  $\mathcal{M}[\mathcal{N}^s]^{\mathbb{U}}$ and the theorem is trivially true. Assume $\mathcal{M}$ or $\mathcal{N}$ to be infinite. We may assume that  $\mathcal{M}[\mathcal{N}^s]^{\mathbb{U}}$ is very saturated.
        Assume 
        \[\left\{\phi_i(\mathbf{x},\mathbf{y}_i), \left(\mathbf{c}_{i,j}=(a_{i,j},b_{i,j})\right)_{j<\omega}, k_i \right\}_{i<k}\]
        is an inp-pattern in $\mathcal{M}[\mathcal{N}^s]^{\mathbb{U}}$ of depth $k> \max(\bdn(\mathcal{M}), \bdn(\mathcal{N})) \geq 1$ with $\vert \mathbf{x} \vert=1$.
        By Ramsey and compactness, we can assume that the sequence of parameters $(a_{i,j},b_{i,j})_{j<\omega}$ for $i<k$ are mutually indiscernible. By quantifier elimination and usual elimination of the disjunctions, we can also assume that for $i<k$ the formula $\phi_i(\mathbf{x},\mathbf{y}_i)$ (where $\mathbf{y}_i=\mathbf{y}_i^1,\ldots, \mathbf{y}_i^{\vert\mathbf{y}_i \vert}$) is a conjunction of formulas of the form:
            \begin{itemize}
                \item $\phi_{i,s}(\mathbf{x},\mathbf{y}_i) = \bigwedge_{l<\vert \mathbf{y}_i\vert} (\neg) s(\mathbf{x},\mathbf{y}_i^l)$, a complete $s$-diagram in $\mathbf{x}$ over $\mathbf{y}_i^1,\ldots, \mathbf{y}_i^{\vert\mathbf{y}_i \vert}$ . 
                \item $\phi_{i,R}(\mathbf{x},\mathbf{y}_i)=\bigwedge (\neg)R(\mathbf{x},\mathbf{y}_i)$, for finitely many $R\in \mathrm{L}$,
                \item $\phi_{i,\mathbb{U}}(\mathbf{x})=\bigwedge(\neg) R^{\mathbb{U}}(\mathbf{x})$, for finitely many $R\in \mathrm{L}$.
            \end{itemize}
         Then, the crucial step of the proof is to remark that no inp-pattern can have a line "talking about" $\mathcal{M}$ and another one "talking about" $\mathcal{N}$.
            
        \textbf{Case 1:} Assume that the complete $s$-diagram in $\mathbf{x}$ of the first line implies $s(\mathbf{x},\mathbf{y}^0_{0,0})$. \\
        
        \begin{claim}
        By consistency of paths, and inconsistency of the lines, the same holds for every line  $i$: there is some $\mathbf{y}_i^{l_i}$, $l_i <\vert \mathbf{y}_i\vert $, such that $\phi_i(\mathbf{x},\mathbf{y}) \vdash s(\mathbf{x}, \mathbf{y}_{i}^{l_i})$.
        \end{claim}
        \begin{proof}
        Assume for example that line 1 implies $\neg s(\mathbf{x},\mathbf{y}_1^l)$ for all $l<\vert \mathbf{y}_1\vert$. By mutual indiscernibility, $a_{0,0}^0 \neq a_{1,j}^l $ for all $j<\omega$ and $l<\vert \mathbf{y}_1\vert$. Take $\textbf{d}_{j} \in \mathcal{M}[\mathcal{N}^s]^{\mathbb{U}}$ satisfying \[\phi_0(\mathbf{d}_{j},\mathbf{c}_{0,0}) \wedge \phi_1(\mathbf{d}_{j},\mathbf{c}_{1,j}),\]
        (consistency of paths). As $s(\textbf{d}_j,\mathbf{c}^0_{0,0})$ and $a_{0,0}^0 \neq a_{1,j}^l$, this implies 
         \[ \mathcal{M} \models \phi_{1,R}(a_{0,0}^0,a_{1,j}) \wedge \phi_{1,\mathbb{U}}(a_{0,0}^0)\]
         and equivalently that 
         \[ \mathcal{M}[\mathcal{N}^s]^{\mathbb{U}} \models  \phi_1(\textbf{c}_{0,0}^0,\mathbf{c}_{1,j}).\]
         Then, the line $1$ would be realised by $\mathbf{c}_{0,0}^0$, contradiction. The same argument holds for any line $i>0$.
         \end{proof}
         Without loss of generality, we may assume that for all $i$, $\phi_i(\mathbf{x},\mathbf{y}) \vdash s(\mathbf{x}, \mathbf{y}_{i}^0)$ ($l_i=0$). It follows from the claim and by consistency of paths that:
        \begin{claim}
            The parameters $\{a_{i,n}^0\}_{i<k,n<\omega}$ are all equal to some parameter $a\in \mathcal{M}$.
        \end{claim}

        Then by inconsistency of lines and by definition of $R^{\mathcal{M}[\mathcal{N}^s]^{\mathbb{U}}}$, there is a subconjunction $\phi_i'(\mathbf{x},\mathbf{y_i}')$ of $\phi_i(\mathbf{x},\mathbf{y}_i)$, where $\mathbf{y}_i'$ is a subtuple of $\mathbf{y}_i$ containing $\mathbf{y}_i^0$,  with a complete (positive) diagram $\bigwedge_{l<\vert \mathbf{y}_i'\vert} s(\mathbf{x},\mathbf{y}_i^l)$ and which already forms an inconsistent line. 
        This is simply to say that, if $\phi_i(\mathbf{x},\mathbf{y_i}) \vdash \neg s(\mathbf{x},\mathbf{y}_i^l)$ and $\phi_i(\mathbf{x},\mathbf{y_i}) \vdash (\neg) R(\mathbf{x},\mathbf{y}_i^l)$ for some $l<\vert y_i\vert$, then
        \[\mathcal{M}\models\{(\neg)R(a,a_{i,j}^l)\}_{j<\omega}\] 
        or equivalently
        \[\mathcal{M}[\mathcal{N}^s]^{\mathbb{U}} \models \{(\neg)R(\mathbf{c},\mathbf{c}_{i,j}^l)\}_{j<\omega}.\]
        whenever $\pi_{\mathcal{M}}(c)=a$. 
        
        So we may assume that the $s$-diagram in $\mathbf{x}$  in any lines is positive (with no occurrence of negation of $s(\mathbf{x},\mathbf{y}_i^{l})$). Thus, the inp-pattern translates to an inp-pattern of $\mathcal{N}$, of depth strictly bigger than $\bdn(\mathcal{N})$, namely:
        \[\{\phi_{i,R}(x,y_i), (b_{i,j})_{j<\omega}\}_{i<k}\]
        This is a contradiction.\\
        \textbf{Case 2:} The complete $s$-diagram of any line is negative.  Then we define the following pattern in $\mathcal{M}$: 
        
        \[P(x):= \left\{\bigwedge_{l<\vert \textbf{y}_i\vert} x \neq y_l \wedge \phi_{i,R}(x,y_i) \wedge \phi_{i,\mathbb{U}}(x), (a_{i,j})_{j<\omega}\right\}_{i<k}.\]
        We can show -- and this is a contradiction-- that it is an inp-pattern of depth strictly greater than $\bdn(\mathcal{M})$.
        Indeed, lines are inconsistent: if $a\in \mathcal{M}$ realises a line 
        \[\left \{\bigwedge_{l<\vert \textbf{y}_i\vert} x \neq a^l_{i,j} \wedge \phi_{i,R}(x,a_{i,j}) \wedge \phi_{i,\mathbb{U}}(x)\right\}_{j<\omega},\]
        then for any $b\in  \mathcal{N}$, $(a,b)$ satisfies the corresponding line of the original pattern 
        \[\{\phi_{i}(\mathbf{x},\mathbf{c}_{i,j})\}_{j<\omega},\]
        this is a contradiction. Paths are consistent: take $\mathbf{d}=(a,b)$ a realisation of the first column of our original inp-pattern:
        \[\{\phi_{i}(\mathbf{x},\mathbf{c}_{i,0})\}_{i<\kappa},\]
        then $a$ satisfies the first column of $P(x)$:
        \[\left \{\bigwedge_{l<\vert \textbf{y}_i\vert} x \neq a^l_{i,0} \wedge \phi_{i,R}(x,a_{i,0}) \wedge \phi_{i,\mathbb{U}}(x)\right \}_{i<\kappa},\]
        This concludes our proof.
        
    \end{proof}
    
    \begin{example}
Let $\mathrm{L}=\{R,B\}$ be the language with two binary predicates, and let $\mathcal{M}$ be a set with two cross-cutting equivalence relations with infinitely many infinite classes. 
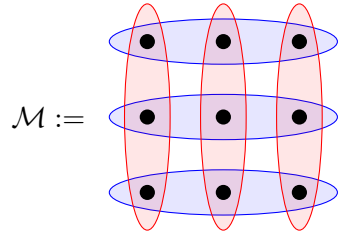
\begin{center} 
\begin{tikzpicture}\clip(-2,-0.7) rectangle (5,2.7);
\draw (-1.3,1) node {$\mathcal{M}:=$};

\draw[blue, fill=blue!10] (1,0) ellipse (1.5 and 0.3);
\draw[blue, fill=blue!10] (1,1) ellipse (1.5 and 0.3);
\draw[blue, fill=blue!10] (1,2) ellipse (1.5 and 0.3);

\fill[red,opacity=0.1] (0,1) ellipse (0.3 and 1.5);
\draw[red] (0,1) ellipse (0.3 and 1.5);
\fill[red,opacity=0.1] (1,1) ellipse (0.3 and 1.5);
\draw[red] (1,1) ellipse (0.3 and 1.5);
\fill[red,opacity=0.1] (2,1) ellipse (0.3 and 1.5);
\draw[red] (2,1) ellipse (0.3 and 1.5);

\fill (0,0) circle (0.1);
\fill (1,0) circle (0.1);
\fill (2,0) circle (0.1);
\fill (0,1) circle (0.1);
\fill (1,1) circle (0.1);
\fill (2,1) circle (0.1);
\fill (0,2) circle (0.1);
\fill (1,2) circle (0.1);
\fill (2,2) circle (0.1);
\end{tikzpicture}
\end{center}
We saw that $\bdn(\mathcal{M})=2$.
We leave it to the reader to describe $\mathcal{M}[\mathcal{M}^s]^{\mathbb{U}}$ (where $R$ and $B$ are no longer equivalence relations). Then, one can verify that it is also of burden $2$. We let the following picture for the intuition:\\
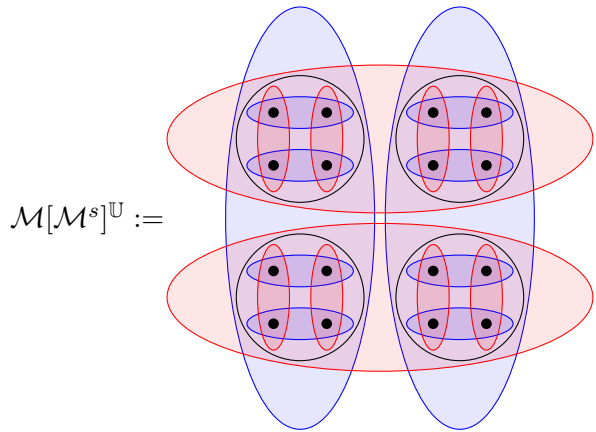
\begin{center}
    
\begin{tikzpicture}[scale=0.7]\clip(-5.5,-2) rectangle (10,6.4);
\draw (-3.5,2) node {$\mathcal{M}[\mathcal{M}^s]^{\mathbb{U}}:=$};
\fill[blue,opacity=0.1] (0.5,2) ellipse (1.4 and 4);
\draw[blue] (0.5,2) ellipse (1.4 and 4);
\fill[blue,opacity=0.1] (3.5,2) ellipse (1.4 and 4);
\draw[blue] (3.5,2) ellipse (1.4 and 4);
\fill[red,opacity=0.1] (2,0.5) ellipse (4 and 1.4);
\draw[red] (2,0.5) ellipse (4 and 1.4);
\fill[red,opacity=0.1] (2,3.5) ellipse (4 and 1.4);
\draw[red] (2,3.5) ellipse (4 and 1.4);

\begin{scope}[shift={(0,0)}]
\fill[blue,opacity=0.1] (0.5,0) ellipse (1 and 0.3);
\draw[blue] (0.5,0) ellipse (1 and 0.3);
\fill[blue,opacity=0.1] (0.5,1) ellipse (1 and 0.3);
\draw[blue] (0.5,1) ellipse (1 and 0.3);

\fill[red,opacity=0.1] (0,0.5) ellipse (0.3 and 1);
\draw[red] (0,0.5) ellipse (0.3 and 1);
\fill[red,opacity=0.1] (1,0.5) ellipse (0.3 and 1);
\draw[red] (1,0.5) ellipse (0.3 and 1);
\draw (0.5,0.5) circle (1.2);
\fill (0,0) circle (0.1);
\fill (1,0) circle (0.1);
\fill (0,1) circle (0.1);
\fill (1,1) circle (0.1);
\end{scope}
\begin{scope}[shift={(0,3)}]
\fill[blue,opacity=0.1] (0.5,0) ellipse (1 and 0.3);
\draw[blue] (0.5,0) ellipse (1 and 0.3);
\fill[blue,opacity=0.1] (0.5,1) ellipse (1 and 0.3);
\draw[blue] (0.5,1) ellipse (1 and 0.3);

\fill[red,opacity=0.1] (0,0.5) ellipse (0.3 and 1);
\draw[red] (0,0.5) ellipse (0.3 and 1);
\fill[red,opacity=0.1] (1,0.5) ellipse (0.3 and 1);
\draw[red] (1,0.5) ellipse (0.3 and 1);
\draw (0.5,0.5) circle (1.2);
\fill (0,0) circle (0.1);
\fill (1,0) circle (0.1);
\fill (0,1) circle (0.1);
\fill (1,1) circle (0.1);
\end{scope}
\begin{scope}[shift={(3,0)}]
\fill[blue,opacity=0.1] (0.5,0) ellipse (1 and 0.3);
\draw[blue] (0.5,0) ellipse (1 and 0.3);
\fill[blue,opacity=0.1] (0.5,1) ellipse (1 and 0.3);
\draw[blue] (0.5,1) ellipse (1 and 0.3);

\fill[red,opacity=0.1] (0,0.5) ellipse (0.3 and 1);
\draw[red] (0,0.5) ellipse (0.3 and 1);
\fill[red,opacity=0.1] (1,0.5) ellipse (0.3 and 1);
\draw[red] (1,0.5) ellipse (0.3 and 1);
\draw (0.5,0.5) circle (1.2);
\fill (0,0) circle (0.1);
\fill (1,0) circle (0.1);
\fill (0,1) circle (0.1);
\fill (1,1) circle (0.1);
\end{scope}
\begin{scope}[shift={(3,3)}]
\fill[blue,opacity=0.1] (0.5,0) ellipse (1 and 0.3);
\draw[blue] (0.5,0) ellipse (1 and 0.3);
\fill[blue,opacity=0.1] (0.5,1) ellipse (1 and 0.3);
\draw[blue] (0.5,1) ellipse (1 and 0.3);

\fill[red,opacity=0.1] (0,0.5) ellipse (0.3 and 1);
\draw[red] (0,0.5) ellipse (0.3 and 1);
\fill[red,opacity=0.1] (1,0.5) ellipse (0.3 and 1);
\draw[red] (1,0.5) ellipse (0.3 and 1);
\draw (0.5,0.5) circle (1.2);
\fill (0,0) circle (0.1);
\fill (1,0) circle (0.1);
\fill (0,1) circle (0.1);
\fill (1,1) circle (0.1);
\end{scope}

\end{tikzpicture}

\end{center}

\end{example}
\section{More on the $\RV$-sort}\label{AxiomatisationRV}
        We defined in Paragraph \ref{SusubsectionRVSort} the leading term structure of order $0$ 
        $(\RV,\oplus,\cdot,\mathbf{0},\mathbf{1})$ of a valued field $\mathcal{K}=(K,\Gamma,k)$ as the abelian group $K^\star/1+\mathfrak{m}$ endowed with some extra structure, where $\mathfrak{m}$ is the maximal ideal of the valuation ring. We saw that in the context of benign valued fields (Definition \ref{DefBenign}), and like the value group and residue field, the $\RV$-sort appears to be a pure imaginary sort (even with control of parameters in the sense of Definition \ref{DefinitionPureWithControlOfParameters}). This means that, like the value group and the residue field, it can be seen as a structure on its own. However, $\RV$-sorts have not  been extensively studied as algebraic structures. Let us cite Krasner, who defined in \cite{Kra44} a generalization of fields called `corpoïde' and which includes in particular the $\RV$-sorts.
        In this appendix, we give the definition of the $\RV$-sort out of the context of valued fields -- which we call an $\RV$-structure. Then we give an axiomatisation in a language $\mathrm{L}:= \{\oplus,\cdot,\mathbf{0}, \mathbf{1}\}$ with two constants and two binary functions, and show that every $\RV$-structure is the $\RV$-sort of a canonical Henselian valued field. 
        
    \subsection{Axomatisation of $\RV$-structures}

        \begin{definition}
        An $\RV$-\emph{structure} is an abelian group $(\RV^\star,\cdot)$ such that there are an ordered abelian group $(\Gamma,+,<,0) $ and a field $(k,+,\cdot,0,1)$ and the following short exact sequence:
        \[1\rightarrow k^\star \overset{\iota}{\rightarrow} \RV^\star \overset{\val_{\RV}}{\rightarrow} \Gamma \rightarrow 0.\]
        \end{definition}

        We abusively see $k^\star$ as a subset of $\RV^\star$. Outside of $k$, the addition is not defined, but one can extend it as follows. We call the group morphism $\val_{\RV}: \RV^\star \rightarrow \Gamma$ the \textit{valuation map}.

    \begin{definition}
        We first complete the embedding of $k^\star$ in $\RV^\star$ by adding a new element $\textbf{0}$, absorbing for the multiplication (\textit{i.e.} $\textbf{0}\cdot \textbf{a} =0$ for all $\textbf{a}\in \RV^\star \cup \lbrace \textbf{0}\rbrace$). We denote by $\RV$ the set $\RV^{\star} \cup \lbrace \textbf{0}\rbrace$. We also add a new element $\infty$ in $\Gamma$, with $\infty> \Gamma$, and extend the valuation map $\val_{\RV}$ in $\textbf{0}$ by setting $\val_{\RV}(\textbf{0})=\infty$. Let $\textbf{a},\textbf{b} \in \RV$. We denote by $\textbf{a} \oplus \textbf{b}$ the element:
        \[\textbf{a} \oplus \textbf{b} = \begin{cases}
            \textbf{0} \ \text{ if } \textbf{\textbf{a}}=\textbf{\textbf{b}}=0,\\
            \textbf{a} \ \text{ if } \val_{\RV}(\textbf{b})>\val_{\RV}(\textbf{a}),\\
            \textbf{b} \ \text{ if } \val_{\RV}(\textbf{a})>\val_{\RV}(\textbf{b}),\\
            (\textbf{a}/\textbf{b}+1)\cdot \textbf{b} \text{ otherwise.}
        \end{cases}
        \]

        This new law $\oplus$ extends the addition in $k$. The phenomenon ${\textbf{a} \oplus \textbf{b}=\textbf{a}}$ or $\textbf{a} \oplus \textbf{b} = \textbf{b}$, corresponding to the first two cases, will be referred to as `\textit{additive absorption}'. Notice that $\oplus$ is not associative. In particular the addition in $K$  is not compatible with $\oplus$ through the valuation maps (e.g.: in $\mathbb{R}((T)), \ 1=\val_{\RV}(1-1+t) \neq \val_{\RV}(\rv(1)\oplus \rv(-1+t))=+\infty$). 
    \end{definition}
    As the operation $\oplus$ is not (always) associative, we adopt a natural convention for the contraction $\bigoplus\limits_{i<n} \mathbf{a}_i$:
    \begin{notation}
        Let $n\in \mathbb{N}$ and $\textbf{a}_1, \ldots, \mathbf{a}_n \in \RV$. We will use the notation $\bigoplus\limits_{i<n} \textbf{a}_i$ for $((((\textbf{a}_{1} \oplus \textbf{a}_{2})\oplus \textbf{a}_{3}) + \cdots )\oplus \textbf{a}_{n})$ only when the sum $((((\textbf{a}_{\sigma(1)} \oplus \textbf{a}_{\sigma(2)})\oplus \textbf{a}_{\sigma(3)}) + \cdots )\oplus \textbf{a}_{\sigma(n)})$ is associative for any permutation $\sigma$ (in other word, when this term does not depend on the choice of parenthesis). It happens in particular when all $\textbf{a}_i$'s have the same valuation (dividing by $\textbf{a}_1$, we get a sum in the field $k$) or \textit{a contrario}, when they have pairwise distinct valuations (by additive absorption). Notice for instance that in the $\RV$-sort of $\mathbb{R}((T))$, the sum $\textbf{1} \oplus (-\textbf{1} \oplus \textbf{t})$ is not `associative' whereas the sum $\textbf{1} \oplus (\textbf{t} \oplus -\textbf{1})$ is. 
    \end{notation}

    Even if this operation $\oplus$ does not satisfy all the required properties, we refer to it as the `addition' of $\RV$. We attempt to describe its essential properties. Recall that we defined $\mathrm{L}$ as the one-sorted language with signature $\lbrace\oplus,\cdot,\textbf{0},\textbf{1}\rbrace$. We will see that $\RV$-structures are exactly $\mathrm{L}$-structures satisfying the following list of axioms (1-7): 
        \begin{enumerate}
            \item $(\RV^\star,\cdot,\textbf{1})$ is an abelian group, where $\RV^\star=\RV \setminus \lbrace \textbf{0}\rbrace$,
            \item (neutral element for $\oplus$) $\forall \textbf{a}\in \RV, \ \textbf{0} \oplus \textbf{a} = \textbf{a}, $
            \item (semi or half-associativity) $[(\textbf{a}\oplus \textbf{b}) \oplus \textbf{c} \neq \textbf{a} \oplus (\textbf{b} \oplus \textbf{c})] \Rightarrow \textbf{a}\oplus \textbf{b} =\textbf{0} \text{ or } \textbf{b} \oplus \textbf{c}=\textbf{0},$
            \item (commutativity for $\oplus$) $\forall \textbf{a},\textbf{b} \in \RV, \ \textbf{a}\oplus \textbf{b} =\textbf{b} \oplus \textbf{a}$,
            \item (distributivity) $\forall \textbf{a},\textbf{b},\textbf{c} \in \RV, \ (\textbf{a} \oplus \textbf{b}) \cdot \textbf{c} = \textbf{ac} \oplus \textbf{bc} $.
        \end{enumerate}
        We define $k^\star$ as the set:
        $$\lbrace \textbf{r}\in \RV \setminus\lbrace \textbf{0}\rbrace \ \vert \ \textbf{1} \oplus \textbf{r} \neq \textbf{1} \ \text{ and } \ \textbf{1}\oplus \textbf{r}^{-1} \neq \textbf{1} \rbrace.$$
        We write $k:=k^\star \cup \lbrace 0 \rbrace$ and we may denote its elements $a,b,c,\ldots\in k$ with the usual font and denote the restriction $\restriction{\oplus}{k}$ by the symbol $+$.
        \begin{enumerate}
            \setcounter{enumi}{5}
            \item $(k:=k^\star \cup \lbrace 0 \rbrace,\cdot,+,0,1)$ is a field,
            \item (uniform additive absorption) $\forall \textbf{a}\in \RV, \ \forall r \in k^\star, \ (\textbf{a} \oplus \textbf{1} = \textbf{1}) \Leftrightarrow (\textbf{a} \oplus r = r)$.
        \end{enumerate}
        
        The fact that all these properties are satisfied by any $\RV$-structure is clear. From the axioms, one can show the following (8-10):
        \begin{enumerate}
            \setcounter{enumi}{7}
            \item (multiplicative absorption) $\forall \textbf{a}\in \RV, \ \textbf{0} \cdot \textbf{a} = \textbf{0}. $
        \end{enumerate}
            Assume for some $\textbf{a}\in \RV$, $\textbf{0}\cdot \textbf{a} \neq \textbf{0}$. Then as $\textbf{0}$ is a neutral element for $\oplus$ and by distributivity $\textbf{0}\cdot \textbf{a} = \textbf{0}\cdot \textbf{a} \oplus \textbf{0}\cdot \textbf{a}$. We may multiply by the multiplicative inverse of $\textbf{0}\cdot \textbf{a}$ and we get a contradiction in $k$.
        \begin{enumerate}
            \setcounter{enumi}{8}
            \item (additive inverse)  $\forall \textbf{a} \in \RV, \exists ! \textbf{b} \ \textbf{a}\oplus \textbf{b} = \textbf{0}$.   
        \end{enumerate}
        This inverse is given by $-\textbf{a} := -\textbf{1}\cdot \textbf{a}$. Indeed we have $(\textbf{a} \oplus -\textbf{a})=\textbf{a}\cdot(\textbf{1}+-\textbf{1})=\textbf{0}$. Uniqueness is clear if $\textbf{a}=\textbf{0}$. Assume $\textbf{a}\neq \textbf{0}$, if $\textbf{b}\in \RV^\star$ is such that $\textbf{a}\oplus \textbf{b}=\textbf{0}$, in particular $\textbf{b} \neq \textbf{0}$. Then by distributivity and multiplicative absorption $\textbf{b}/\textbf{a} \in k^\star$ and from $\textbf{b}/\textbf{a}+\textbf{1}=\textbf{0}$ we get $\textbf{b}=-\textbf{a}$.\\
        
        We recover the value group by setting $\Gamma := \RV^\star/k^\star$. For the order in $\Gamma$, one must define it as follows:
        $$\forall [\textbf{a}],[\textbf{b}] \in \Gamma, [\textbf{a}] < [\textbf{b}] \ \Leftrightarrow \ \textbf{1}\oplus \textbf{b}/\textbf{a}= \textbf{1} \ \Leftrightarrow \ \textbf{a}\oplus \textbf{b}= \textbf{a} ,$$
        where $[\textbf{a}]$ denote the classe of $\textbf{a}$ modulo $k^\star$. By uniform additive absorption and distributivity, this definition does not depend of the representative $\textbf{a}$ and $\textbf{b}$ we have chosen. Indeed, if $r,r' \in k^\star$, one gets:
        $$1 \oplus \textbf{b}r/\textbf{a}r'=\textbf{1} \Leftrightarrow r'/r \oplus \textbf{b}/\textbf{a}= r'/r \overset{(7)}{\Leftrightarrow} \textbf{1}\oplus \textbf{b}/\textbf{a} =\textbf{1}.$$
        \begin{enumerate}
            \setcounter{enumi}{9}
            \item $(\Gamma,<)$ is an ordered group.   
        \end{enumerate}
        
        Anti-symmetry of $<$ follows from the definitions of $k^\star$ and $<$, and transitivity is given by semi-associativity:
        Assume $\textbf{a}\oplus \textbf{b}=\textbf{a}$ and $\textbf{b} \oplus \textbf{c}=\textbf{b}$, then either $\textbf{b}=\textbf{0}$, or $\textbf{a} \oplus \textbf{b} \neq \textbf{0}$ and $\textbf{b} \oplus \textbf{c} \neq \textbf{0}$. In any case, $\textbf{a} \oplus \textbf{c} = (\textbf{a} \oplus \textbf{b}) \oplus \textbf{c} = \textbf{a} \oplus (\textbf{b} \oplus \textbf{c}) = \textbf{a} \oplus \textbf{b} = \textbf{a}$. It's a total order since for all $\textbf{a},\textbf{b} \in \RV^\star$, either $\textbf{a}/\textbf{b} \in k^\star$, $\textbf{a}/\textbf{b}\oplus \textbf{1} = \textbf{1}$ or $\textbf{b}/\textbf{a} \oplus \textbf{1} =\textbf{1}$, which respectively gives $[\textbf{a}]=[\textbf{b}]$, $[\textbf{a}]>[\textbf{b}]$ or $[\textbf{b}]>[\textbf{a}]$. 
        We complete the valuation map by setting $\val_{\RV}(\textbf{0}) =\infty$ where $\infty> \Gamma$. 
        
    \begin{note} 
        \begin{itemize}
            \item To avoid the use of conventions for $\textbf{0}$, one might define the quotient $\RV/k^\star$ as the set of orbits of $\RV$ under the action of $k^\star$. This action preserve the multiplication in $\RV$. We get that $\textbf{0}$ is the unique element in its orbit $[\textbf{0}]$, which we denoted by $\infty$. the definition of $<$ gives then that $\infty> \Gamma$.
            \item In \cite{Kra44}, structures satisfying Axioms (1) and (8) are called pseudo-group for the multiplication. 
            \item One can replace (6) by:
            \[(6') \ \forall r\in k, \ 1+r \in k, 0\neq 1,  \text{ and } + \text{is associative in } k.\]
            The fact that $(k,\cdot)$ is a multiplicative group can already be deduced from (7). Let $\textbf{r},\textbf{s} \in \RV^\star$. If $\textbf{s} \in k^\star$, one has:
            \[  \textbf{1} \oplus \textbf{r}\cdot \textbf{s} = \textbf{1} \overset{(5)}{\Leftrightarrow} \textbf{1}/\textbf{s} \oplus \textbf{r} = \textbf{1}/\textbf{s} \overset{(7)}{\Leftrightarrow} \textbf{1}+\textbf{r}=\textbf{r}\]
            Hence, a product $\textbf{r}\cdot \textbf{s}$ is not in $k$ if and only if $\textbf{s}$ or $\textbf{r}$ is not in $k$. \\
        \end{itemize}
    \end{note}

    \subsection{The Hahn field associated to $\RV$}
    
        As we know, given any field $k$ and any ordered abelian group $\Gamma$, there is always a Henselian valued field with residue field $k$ and value group $\Gamma$: the Hahn field $k((\Gamma))$.
        We can ask the following question: is any $\RV$-structure the $\RV$-sort of a certain Henselian valued field? 
        We define for that the  Hahn field associated to  $\RV$.
        
    \begin{definition}[of the Hahn field $\RV^{(\Gamma)}$] \index{Hahn field associated to an $\RV$-structure} \nomenclature[]{$\RV^{(\Gamma)}$}{Hahn field associated to $\RV$}
        The Hahn field associated to the $\RV$-structure $(\RV, k, \Gamma)$ and denoted by $\RV^{(\Gamma)}$, is defined by the following set:
        \[ \lbrace (\textbf{a}_\gamma)_{\gamma \in \Gamma} \ \vert \ \forall \gamma\in \Gamma \ \mathbf{a}_\gamma \in \RV, \val_{\RV}(\mathbf{a}_\gamma) \in \lbrace \gamma,\infty \rbrace
         \text{ and } \supp(\mathbf{a}_\gamma)_{\gamma} \ \text{is well-ordered} \rbrace\]
        where $\supp(\mathbf{a}_\gamma)_\gamma = \lbrace \gamma \in \Gamma \ \vert \ \mathbf{a}_{\gamma} \neq \textbf{0}\rbrace$.\\
        It is endowed with the following laws:
        \[(\mathbf{a}_\gamma)_\gamma + (\mathbf{b}_\gamma)_\gamma = (\mathbf{a}_\gamma \oplus \mathbf{b}_\gamma)_\gamma\]
        \[(\mathbf{a}_\gamma)_\gamma \cdot (\mathbf{b}_\gamma)_\gamma = (\bigoplus\limits_{\delta+\epsilon =\gamma}\mathbf{a}_{\delta}\cdot \mathbf{b}_{\epsilon})_\gamma\]
        \[\val(\mathbf{a}_\gamma)_\gamma = \min \supp(\mathbf{a}_\gamma)_\gamma\]
    \end{definition}
    An element $(\mathbf{a}_\gamma)_{\gamma \in \Gamma} \in \RV^{(\Gamma)}$ is written $\sum\limits_{\gamma \in \Gamma} \mathbf{a}_\gamma$ where the sign $\sum$ is purely formal. 
    
    \begin{proposition}
        The Hahn field $\RV^{(\Gamma)}$ is a Henselian valued field of $\RV$-sort $\RV$. 
    \end{proposition}
    
    \begin{proof}
        The proof is straightforward. As in the Hahn field $k((t^\Gamma))$, the difficult part is to show that every non-zero element of $\RV^{(\Gamma)}$ has a multiplicative inverse. We first show that it is a spherically complete ring. Then we will deduce that it is an actual field. 
        \begin{itemize}
        
            \item (associativity for $+$): if $\textbf{a},\textbf{b},\textbf{c} \in \RV$ with $\val_{\RV}(\textbf{a})=\val_{\RV}(\textbf{b})=\val_{\RV}(\textbf{c})$, then $(\textbf{a}\oplus \textbf{b}) \oplus \textbf{c} = \textbf{a} \oplus (\textbf{b} \oplus \textbf{c})$. Then associativity for $+$ in  $\RV^{(\Gamma)}$ is clear as we sum componentwise.
            \item (commutativity for $+$): clear as $\oplus$ is commutative in $\RV$.
            \item (neutral element for $+$): $0 :=\sum_\gamma0 \in \RV^{(\Gamma)}$ is a neutral element, as $0 \in \RV$ is a neutral element for $\oplus$.
            \item (inverse for $+$): if $a=\sum\limits_{\gamma\in \Gamma}\mathbf{a}_\gamma$, the inverse of $a$ is given by $a=\sum\limits_{\gamma\in \Gamma}-\mathbf{a}_\gamma$. The support being the same, it is an element of $\RV^{(\Gamma)}$. 
        \end{itemize}
        The multiplication in $\RV^{(\Gamma)}$ is well-defined: as the supports of $a$ and $b$ are well-ordered, the sum $\bigoplus\limits_{\delta+\epsilon =\gamma}\mathbf{a}_{\delta}\cdot \mathbf{b}_{\epsilon}$ is finite. As before, it is associative since all terms have same valuation. Then $\supp(a+b) \subset \supp(a)+\supp(b)$ and it is easy to see that is a well-ordered set of $\Gamma$. We have:
        \begin{itemize}
            \item (associativity for $\cdot$): Let $a=\sum\limits_{\delta \in \Gamma} \mathbf{a}_\delta,  b=\sum\limits_{\epsilon \in \Gamma} \mathbf{b}_\epsilon ,  c=\sum\limits_{\zeta \in \Gamma}\mathbf{c}_\zeta$ be three elements of $\RV^{(\Gamma)}$, then a simple calculation gives $(a\cdot b) \cdot c = \bigoplus\limits_{\delta+\epsilon+\zeta =\gamma} (\mathbf{a}_\delta \cdot \mathbf{b}_\epsilon)\cdot \mathbf{c}_\zeta =\bigoplus\limits_{\delta+\epsilon+\zeta =\gamma} \mathbf{a}_\delta \cdot (\mathbf{b}_\epsilon\cdot \mathbf{c}_\zeta) = a\cdot (b \cdot c)$  (as $\cdot$ is associative in $\RV$).
            \item (commutativity for $\cdot$): clear as $\cdot$ is commutative in $\RV$.
            \item (neutral element for $\cdot$): Let $\textbf{1}\in \RV^{(\Gamma)}$ be the element $\sum_{\gamma}\mathbf{a}_\gamma$ where $\mathbf{a}_\gamma= 
            \begin{cases} 
                \textbf{0} \text{ if } \gamma \neq 0 \\
                \textbf{1} \text{ if } \gamma=0
            \end{cases}.$
            It is a neutral element for $\cdot$ as $\textbf{1}\in \RV$ is a neutral element for the multiplication in $\RV$.
            \item (distributivity): Let $a=\sum\limits_{\gamma \in \Gamma}\mathbf{a}_\gamma,b=\sum\limits_{\gamma \in \Gamma}\mathbf{b}_\gamma,c=\sum\limits_{\gamma \in \Gamma}\mathbf{c}_\gamma$ be three elements of $\RV^{(\Gamma)}$, then: 
        \begin{eqnarray*} (a+b)\cdot c &=& \sum_{\gamma \in \Gamma} \bigoplus_{\delta+\epsilon = \gamma} (\mathbf{a}_\delta \oplus \mathbf{b}_\delta) \cdot \mathbf{c}_\epsilon = \sum_{\gamma \in \Gamma} \bigoplus_{\delta+\epsilon = \gamma} (\mathbf{a}_\delta\cdot \mathbf{c}_\epsilon \oplus \mathbf{b}_\delta\cdot \mathbf{c}_\epsilon) \\
        &=& \sum_{\gamma \in \Gamma} \bigoplus_{\delta+\epsilon = \gamma} \mathbf{a}_\delta\cdot \mathbf{c}_\epsilon + \sum_{\gamma \in \Gamma} \bigoplus_{\delta+\epsilon = \gamma} \mathbf{b}_\delta\cdot \mathbf{c}_\epsilon\\
        &=&a\cdot c+b\cdot c.
        \end{eqnarray*}
            \item (valuation): Homomorphism of groups is clear. The ultrametric inequality of the valuation is clear from the definition.
            \item (spherically complete): We give here a usual diagonal argument. Let $(a^i)_{i<\lambda}$ be a pseudo-Cauchy sequence in $\RV^{(\Gamma)}$, where $\lambda$ is any limit ordinal. There is $i_0$ such that for all $i_0<i<j<k$, $\val_{\RV}(a^i-a^j)<\val_{\RV}(a^j-a^k)$. For $i>i_0$, we denote by $\gamma_i$ the value $\val_{\RV}(a^i-a^{i+1})$.  We define 
            \[\mathbf{a}_\gamma :=  \begin{cases} \mathbf{a}_\gamma^i \text{ if }\gamma_i>\gamma \text{ for some  }i<\lambda,\\
        0 \text{ otherwise. } 
        \end{cases}\]
        It is well defined by definition of the valuation (it does not depend on the choice of $i$). Let $a = \sum_\gamma \mathbf{a}_\gamma$, we get $\val_{\RV}(a-a^i)>\gamma_i$. We have proved that any pseudo-Cauchy sequence admits a pseudo-limit. 
            \item (multiplicative inverse) Let $a=\sum_\gamma \mathbf{a}_\gamma\in \RV^{(\Gamma)}$. Assume it has no inverse and consider the set 
            \[\Delta :=\lbrace \val(1-a\cdot b) \ \vert \ b \in \RV^{(\Gamma)} \rbrace.\]
            This set has no maximal element. Indeed, notice first if $\gamma=\val(1-a\cdot b)$, then $\gamma<\infty$ as $a$ has no inverse. It follows that if $\mathbf{c}_\gamma \in \RV$ is the coefficient of value $\gamma$ in $c=1-a\cdot b$, we have $\val(1-a\cdot(b-\mathbf{c}_\gamma\cdot \mathbf{a}_{\val(a)}^{-1}))>\gamma$. So $\Delta$ has no maximal element.
            Let $(\gamma_\nu=\val(1-a\cdot b_\nu))_{\nu\in \lambda}$ be an co-final increasing sequence in $\Delta$. Then, $\lambda$ is a limit ordinal. By definition, $(a\cdot b_\nu)_\nu$ is pseudo-Cauchy with pseudo-limit $1$.  Then, $(b_\nu)_{\nu\in \lambda}$ is a pseudo-Cauchy sequence (as multiplication by $a$ preserve pseudo-Cauchy sequences). It converges to an element $b\ in \RV^{(\Gamma)}$. Then $\val(1-a\cdot b) > \Delta$. Contradiction.
        \end{itemize}
        We have proved that $\RV^{(\Gamma)}$ is a spherically complete valued field, so in particular Henselian. Every element $\sum_\gamma \mathbf{a}_\gamma \in \RV^{(\Gamma)}$ can be written as $ {\mathbf{a}_\delta(1+\sum_\gamma \mathbf{a}_\gamma/\mathbf{a}_\delta)}$ where $\delta = \val(\sum_\gamma \textbf{a}_\gamma)$. Clearly, we have that $\RV(\RV^{(\Gamma)}) =\RV$.
    \end{proof}
    
    \begin{remarks}
        \begin{itemize}
            \item In the case where $\RV= k \times \Gamma$, one may easily show that $\RV^{(\Gamma)}$ is isomorphic to the Hahn field $k((\Gamma))$.
            \item The characteristic of $\RV^{(\Gamma)}$ is always equal to the characteristic of the residue field $k$. 
        \end{itemize}
    
    \end{remarks}

\newpage
\begin{acknowledgements}
  I would like to thank the logic team of Münster for all their support. In particular, many thanks to Martin Bays for all our discussions, Franziska Jahnke for her advice on mixed characteristic valued fields and Gonenç Onay for his enlightenment and his support. Many thanks to Yatir Halevi for his very useful comments on the previous version of this preprint and for his contribution on Theorem \ref{ThmBdnExSeq} (all errors are mine). Many thanks to Artem Chernikov and Martin Ziegler for providing me with an early version of their work in \cite{ACGZ20}. Finally, my immense gratitude to my PhD advisor Martin Hils for the numerous discussions who made the writing of this paper possible.
\end{acknowledgements}
\bibliographystyle{plain}
\bibliography{Bibtex}
\newpage

\end{document}